\documentclass[12pt]{article}

\usepackage{indentfirst}
\usepackage{amsfonts}
\usepackage{amsmath}
\usepackage{amssymb}
\usepackage{amsbsy, amsthm}
\usepackage[margin=0.8in]{geometry}

\newtheorem{dfn}{Definition}[section]
\newtheorem{theorem}[dfn]{Theorem}
\newtheorem{lemma}[dfn]{Lemma}

\newtheorem{corollary}[dfn]{Corollary}
\newtheorem{conjecture}[dfn]{Conjecture}

\newenvironment{pf}{\noindent{\bf Proof.}}
{\enspace\vrule height5pt depth0pt width5pt}
\def\deg {{\rm deg}}

\def\E {{\mathcal E}}

\def\F {{\mathcal F}}

\def\M {{\mathcal M}}

\def\C {{\mathcal C}}

\def\G {{\mathcal G}}

\def\CT {{\rm CT}}
\def\A {{\mathcal A}}

\def\td {{\rm \overline{td}}}

\begin{document}
\title{Defective coloring is perfect for minors}
\author{Chun-Hung Liu\thanks{chliu@tamu.edu. Partially supported by NSF under award DMS-1954054 and CAREER award DMS-2144042.} \\
\small Department of Mathematics, \\
\small Texas A\&M University,\\
\small College Station, TX 77843-3368, USA}

\maketitle

\begin{abstract}
\noindent
The defective chromatic number of a graph class is the infimum $k$ such that there exists an integer $d$ such that every graph in this class can be partitioned into at most $k$ induced subgraphs with maximum degree at most $d$.
Finding the defective chromatic number is a fundamental graph partitioning problem and received attention recently partially due to Hadwiger's conjecture about coloring minor-closed families.
In this paper, we prove that the defective chromatic number of any minor-closed family equals the simple lower bound obtained by the standard construction, confirming a conjecture of Ossona de Mendez, Oum, and Wood.
This result provides the optimal list of unavoidable finite minors for infinite graphs that cannot be partitioned into a fixed finite number of induced subgraphs with uniformly bounded maximum degree.
As corollaries about clustered coloring, we obtain a linear relation between the clustered chromatic number of any minor-closed family and the tree-depth of its forbidden minors, improving an earlier exponential bound proved by Norin, Scott, Seymour, and Wood and confirming the planar case of their conjecture. 
\end{abstract}

\section{Introduction}

A {\it proper coloring} of a graph is a function that maps each vertex to a color so that no adjacent vertices receive the same color.
The minimum number of required colors to properly color a graph $G$ is the {\it chromatic number} $\chi(G)$.
Clearly, $\chi(G) \geq \omega(G)$, where $\omega(G)$ is the maximum size of a set of pairwise adjacent vertices in $G$, called the {\it clique number} of $G$.
It is well-known that the gap between $\chi(G)$ and $\omega(G)$ can be arbitrarily large.
Looking for sufficient conditions for graphs $G$ to ensure a certain relationship between $\chi(G)$ and $\omega(G)$ is a very active area. 
To study such sufficient conditions, it is convenient to consider a graph class instead of just a graph.
For a graph parameter $p$ and a class $\F$ of graphs, we define $p(\F)=\sup_{G \in \F}p(G)$.

An extreme case is to study the graphs $G$ with $\chi(G)=\omega(G)$.
However, the disjoint union of an arbitrary graph $H$ and $K_{\chi(H)+1}$ is a graph $G$ satisfying $\chi(G)=\omega(G)$.
So nothing informative can be said about $G$ unless we also consider the substructures of $G$.
It leads to the notion of perfect graphs.
A graph $G$ is {\it perfect} if $\chi(H)=\omega(H)$ for every induced subgraph $H$ of $G$.
The celebrated Strong Perfect Graph Theorem \cite{crst} provides a structural characterization of perfect graphs.

The perfectness can also be defined for graph classes.
For a partial order $\preceq$ defined on graphs, we say that a class $\F$ of graphs is {\it $\preceq$-closed} if for any graphs $G$ and $H$ with $H \preceq G$, $G \in \F$ implies $H \in \F$.
And we say that a $\preceq$-closed graph class $\F$ is {\it perfect with respect to $\preceq$} if $\chi(\F')=\omega(\F')$ for every $\preceq$-closed subclass $\F'$ of $\F$. 
It is straightforward to show that a graph $G$ is perfect if and only if the class of all induced subgraphs of $G$ is perfect with respect to the induced subgraph relation.

Another example for chasing perfection is Hadwiger's conjecture about graph minors.
A graph $H$ is a {\it minor} of another graph $G$ if $H$ is isomorphic to a graph that can be obtained from a subgraph of $G$ by contracting edges.
It is easy to show that every planar graph has no $K_5$-minor.
So the case $t=4$ of the following conjecture proposed by Hadwiger implies the Four Color Theorem \cite{ah,ahk,rsst}. 

\begin{conjecture}[\cite{h}] \label{Hadwiger_conj}
For every positive integer $t$, if $\F$ is the class of $K_{t+1}$-minor free graphs, then $\chi(\F) \leq t$. 
\end{conjecture}

Conjecture \ref{Hadwiger_conj} is very difficult.
Wagner \cite{w} proved that the case $t=4$ is equivalent to the Four Color Theorem, and Robertson, Seymour, and Thomas \cite{rst} proved the case $t = 5$.
The cases $t \geq 6$ remain open.
On the other hand, Delcourt and Postle \cite{dp} proved that the chromatic number of $K_{t+1}$-minor free graphs is $O(t\log\log t)$, which is the currently best known upper bound and improves earlier results in \cite{k1,k2,nsp,t}. 

Note that Conjecture \ref{Hadwiger_conj} is equivalent to stating that $\chi(\F)=\omega(\F)$ for every minor-closed family $\F$.
In other words, Hadwiger's conjecture is equivalent to stating that every minor-closed family is perfect with respect to the minor relation.

Due to the infamous difficulty of Hadwiger's conjecture, some relaxations of proper coloring were considered.
One such relaxation is called defective coloring.
For any real number $r$, we define $[r]$ to be the set $\{x \in {\mathbb N}: 1 \leq x \leq r\}$.
For every positive integer $k$, a {\it $k$-coloring} of a graph $G$ is a function $f: V(G) \rightarrow [k]$.
For positive integers $k$ and $d$, a $k$-coloring $f$ of a graph $G$ has {\it defect} $d$ if for every $i \in [k]$, the maximum degree of the subgraph of $G$ induced by the vertices with color $i$ is at most $d$.
In other words, a graph has a $k$-coloring with defect $d$ if and only if it can be partitioned into at most $k$ induced subgraphs with maximum degree at most $d$.
Colorings with defect 0 are exactly proper colorings.

For a graph class $\F$, the {\it defective chromatic number} of $\F$, denoted by $\chi_\Delta(\F)$, is the infimum $k$ such that there exists an integer $N$ such that every graph in $\F$ has a $k$-coloring with defect $N$.
Namely, the defective chromatic number of $\F$ is the minimum number of parts required to partition every graph in $\F$ into induced subgraphs with universally bounded maximum degree (or equals $\infty$ if no such a partition exists).

Edwards, Kang, Kim, Oum, and Seymour \cite{ekkos} proved the following verbatim defective coloring analog of Hadwiger's conjecture.

\begin{theorem}[\cite{ekkos}] \label{verbatim}
For every positive integer $t$, if $\F$ is the class of $K_{t+1}$-minor free graphs, then $\chi_\Delta(\F) \leq t$.
\end{theorem}

The proof of Theorem \ref{verbatim} in \cite{ekkos} is elegant, and the bound $t$ cannot be improved for the class of $K_{t+1}$-minor free graphs.
However, Theorem \ref{verbatim} loses the essence for chasing perfection for minor-closed families.
Theorem \ref{verbatim} is equivalent to stating that $\chi_\Delta(\F) \leq \omega(\F)$ for every minor-closed family.
But unlike the chromatic number, the defective chromatic number is not lower bounded by the clique number.
In fact, the gap between $\chi_\Delta(\F)$ and $\omega(\F)$ can be arbitrarily large.
For example, for every positive integer $t$, if $\F$ is the class of graphs whose every component has at most $t$ vertices, then $\F$ is minor-closed, $\omega(\F)=t$, and $\chi_\Delta(\F)=1$.

To chase perfection for minor-closed families with respect to defective coloring, we should consider the correct analog for the clique number with respect to defective coloring.
The right answer seems to be the closure of rooted trees.

The {\it closure} of a rooted tree $T$ is the graph with vertex-set $V(T)$ such that two distinct vertices $u,v$ are adjacent if and only if one of $u$ and $v$ is an ancestor of the other.
The {\it height} of a rooted tree $T$ is the maximum number of vertices of a path from the root of $T$ to a leaf in $T$. 
For a positive integer $k$, a {\it balanced $k$-ary tree} is a rooted tree such that every path from the root to a leaf has the same length, and every non-leaf has exactly $k$ children.
For positive integers $h$ and $k$, the closure of a balanced $k$-ary tree with height $h$ is denoted by $\CT_{h,k}$.
For example, $\CT_{1,k}=K_1$ and $\CT_{2,k}=K_{1,k}$.

Like a clique, $\CT_{h,k}$ can be constructed by repeatedly adding new vertices adjacent to all existing vertices.
For graphs $G$ and $H$, we define $G \vee H$ to be the graph obtained from a disjoint union of $G$ and $H$ by adding an edge $xy$ for each pair of vertices $x \in V(G)$ and $y \in V(H)$.
For a graph $G$ and a positive integer $k$, we define $kG$ to be the union of $k$ disjoint copies of $G$.
Then it is easy to see that for every positive integer $k$, $\CT_{1,k}=K_1$, and $\CT_{h,k}=K_1 \vee k\CT_{h-1,k}$ for every $h \geq 2$.

Also, $\CT_{h,k}$ provides a lower bound for the defective chromatic number.
A simple induction on $h$ shows that there exists no $(h-1)$-coloring of $\CT_{h,k}$ with defect $k-1$.
Therefore, if a graph class $\F$ contains $\CT_{h,k}$ for infinitely many integers $k$, then $\chi_\Delta(\F) \geq h$.
So $\CT_{h,k}$ seems to be the correct analog of cliques for defective coloring.

For a graph class $\F$, we define 
$$\omega_\Delta(\F):=\sup\{h \in {\mathbb N}: \CT_{h,k} \in \F \text{ for infinitely many positive integers } k\}.$$
So for every graph class $\F$, $\chi_\Delta(\F) \geq \omega_\Delta(\F)$.

The main result of this paper (Theorem \ref{main_intro}) shows that this inequality is always an equality for minor-closed families and hence we obtain a characterization of the defective chromatic number of minor-closed families and obtain perfectness.

\begin{theorem} \label{main_intro}
For every minor-closed family $\F$, $\chi_\Delta(\F)=\omega_\Delta(\F)$.
\end{theorem}

Theorem \ref{verbatim} is an immediate corollary of both Theorem \ref{main_intro} and Hadwiger's conjecture.
And it seems that Theorem \ref{main_intro} and Hadwiger's conjecture are incomparable. 

The parameter $\omega_\Delta$ is closely related to the connected tree-depth.
The {\it connected tree-depth}\footnote{The connected tree-depth is a variant of a more commonly studied parameter, tree-depth. The {\it tree-depth} of a graph $G$ is the maximum of the connected tree-depth of the components of $G$. The difference between the tree-depth and the connected tree-depth is always at most 1.} of a graph $G$, denoted by $\td(G)$, is the minimum $h$ such that $G$ is a subgraph of a closure of a rooted tree with height $h$.
Equivalently, $\td(G)$ is the minimum $h$ such that $G$ is a subgraph of $\CT_{h,k}$ for some integer $k$.
Note that $\omega(\CT_{h,k})=h$ for any positive integers $h$ and $k$, so $\CT_{h,k}$ is not a subgraph of $\CT_{h-1,k'}$ for any $k'$.
Hence $\td(\CT_{h,k})=h$. 
So for every graph class $\F$, $\CT_{\td(\F)+1,k} \not \in \F$ for any $k$, and hence $\omega_\Delta(\F) \leq \td(\F)$.
Moreover, if $\F$ is the class of $H$-minor free graphs for some graph $H$, then $\CT_{\td(H),\ell} \not \in \F$ for all integers $\ell \geq \lvert V(H) \rvert$ (since $H$ is a subgraph of $\CT_{\td(H),\lvert V(H) \rvert}$), and $\CT_{\td(H)-1,k} \in \F$ for all integers $k$ (since every minor of $\CT_{\td(H)-1,k}$ has connected tree-depth at most $\td(H)-1$), so $\omega_\Delta(\F) = \td(H)-1$.

Ossona de Mendez, Oum, and Wood \cite{oow} proposed the following conjecture and proved the case $\td(H) \leq 3$.

\begin{conjecture}[\cite{oow}] \label{defect_conj}
Let $H$ be a graph.
If $\F$ is the class of $H$-minor free graphs, then $\chi_\Delta(\F)=\td(H)-1$.
\end{conjecture}

Note that Conjecture \ref{defect_conj} is significantly stronger than the verbatim defective coloring analog of Hadwiger's conjecture (Theorem \ref{verbatim}).
The key observation is that $K_{t,t}$ contains a $K_{t+1}$-minor.
So $K_{t+1}$-minor free graphs do not contain $K_{t,t}$ as a subgraph and hence enjoy a local condition.
Such a forbidden subgraph condition was used to prove much stronger results, including defective coloring for classes with bounded expansion \cite{oow} and clustered coloring for various graph classes \cite{lw_minor,lw_layer,lw_topo} that were used to prove the currently best results about the verbatim clustered coloring analog of Hadwiger's conjecture and Haj\'{o}s' conjecture for topological minors \cite{lw_topo}.
However, such an argument does not work for Conjecture \ref{defect_conj} because there exist no functions $f$ and $f'$ such that $K_{f(h),f'(h,k)}$ contains $\CT_{h,k}$ for all integers $k$.
So no upper bound for $\chi_\Delta$ only involving $\td(H)$ can be obtained by this argument, even when $H=\CT_{h,k}$ for some large $k$.

In fact, Conjecture \ref{defect_conj} is equivalent to our main result Theorem \ref{main_intro}.
Conjecture \ref{defect_conj} immediately follows from Theorem \ref{main_intro} because we have observed that $\omega_\Delta(\F) = \td(H)-1$ if $\F$ is the class of $H$-minor free graphs.
The converse direction holds since every minor-closed family $\F$ does not contain $\CT_{\omega_\Delta(\F)+1,k}$ for some fixed integer $k$ (so every graph in $\F$ is $\CT_{\omega_\Delta(\F)+1,k}$-minor free) and $\td(\CT_{\omega_\Delta(\F)+1,k})=\omega_\Delta(\F)+1$.

\begin{corollary} \label{defect_cor_intro}
Conjecture \ref{defect_conj} is true.
\end{corollary}

Our Theorem \ref{main_intro} can also be applied to clustered coloring.
The {\it clustered chromatic number} of a graph class $\F$, denoted by $\chi_*(\F)$, is the infimum $k$ such that there exists an integer $N$ such that every graph in $\F$ can be partitioned into at most $k$ induced subgraphs with no component on more than $N$ vertices.
Clearly, $\chi_*(\F) \geq \chi_\Delta(\F)$ for every graph class $\F$.
See \cite{w_survey} for a survey about defective coloring and clustered coloring.

A corollary of a result of Alon, Ding, Oporowski, and Vertigan \cite{adov} builds a connection between the defective chromatic number and the clustered chromatic number for minor-closed families with bounded tree-width\footnote{The {\it tree-width} of a graph $G$ is the minimum $k$ such that $G$ is a subgraph of a chordal graph with clique number at most $k+1$.} (Statement 2 in Theorem \ref{def_clu}).
The author and Oum \cite{lo} proved a tight result without the assumption for having bounded tree-width, building a connection between the defective chromatic number and the clustered chromatic number for all minor-closed families (Statement 1 in Theorem \ref{def_clu}).

\begin{theorem}[\cite{adov,lo}] \label{def_clu}
Let $\F$ be a minor-closed family.
	\begin{enumerate}
		\item Then $\chi_*(\F) \leq 3\chi_\Delta(\F)$.
		\item If $\F$ has bounded tree-width (equivalently, $\F$ does not contain all planar graphs\footnote{The equivalence follows from the Grid Minor Theorem \cite{rs}.}), then $\chi_*(\F) \leq 2\chi_\Delta(\F)$. 
	\end{enumerate}
\end{theorem}

Norin, Scott, Seymour, and Wood \cite{nssw} proposed an analog of Conjecture \ref{defect_conj} for clustered coloring.

\begin{conjecture}[\cite{nssw}] \label{clustered_conj}
For every graph $H$, if $\F$ is the class of $H$-minor free graphs, then $\chi_*(\F) \leq 2\td(H)-2$.
\end{conjecture}

Norin, Scott, Seymour, and Wood \cite{nssw} showed that the bound $2\td(H)-2$ in Conjecture \ref{clustered_conj} is tight for some graph $H$.
They \cite{nssw} also observed that Theorem \ref{def_clu} and the known $\td(H) \leq 3$ case of Conjecture \ref{defect_conj} proved in \cite{oow} imply the $\td(H) \leq 3$ case of Conjecture \ref{clustered_conj}.
In general, they \cite{nssw} showed that $\td(H)-1 \leq \chi_*(\F) \leq 2^{\td(H)+1}-4$ if $\F$ is the class of $H$-minor free graphs.
Moreover, Norin, Scott, and Wood \cite{nsw} proved that Conjecture \ref{clustered_conj} holds if $\F$ also has bounded path-width\footnote{The {\it path-width} of a graph $G$ is the minimum $k$ such that $G$ is a subgraph of an interval graph with clique number at most $k+1$. So the tree-width is at most the path-width.}.
They \cite{nsw} also showed that $\chi_*(\F)=\td(H)-1$ if $\F$ even has bounded tree-depth.

Theorem \ref{main_intro} strengths or rediscovers all results stated in the previous paragraph.
When $\F$ has bounded tree-depth, every graph in $\F$ does not have arbitrarily long paths, so $\chi_*(\F)=\chi_\Delta(\F)$, and hence Theorem \ref{main_intro} rediscovers the result $\chi_*(\F)=\td(H)-1$ in \cite{nsw} for any graph $H$ such that the class $\F$ of $H$-minor free graphs has bounded tree-depth.
And by combining Corollary \ref{defect_cor_intro} and Theorem \ref{def_clu}, we immediately obtain the following results, which improve the exponential bound in \cite{nssw} to a linear bound and generalize the bounded path-width case in \cite{nsw} to the bounded tree-width case.

\begin{corollary}
Let $H$ be a graph and $\F$ the class of $H$-minor free graphs.
	\begin{enumerate}
		\item Then $\chi_*(\F) \leq 3\td(H)-3$.
		\item If $\F$ has bounded tree-width (equivalently, $H$ is planar), then $\chi_*(\F) \leq 2\td(H)-2$. 
	\end{enumerate}
\end{corollary}

Our Theorem \ref{main_intro} is also related to known results about fractional defective and clustered coloring.
The {\it fractional defective chromatic number $\chi_\Delta^f(\F)$} (and {\it fractional clustered chromatic number $\chi_*^f(\F)$}, respectively) of a graph class $\F$ is the infimum $k$ such that for every $k'>k$, there exists an integer $d$ such that for every graph $G \in \F$, there exist a real number $p$ and at most $pk'$ induced subgraphs of $G$ with maximum degree at most $d$ (and with no component on more than $d$ vertices, respectively) such that every vertex of $G$ is contained in at least $p$ of them.
Note that the integer $d$ is allowed to be dependent on $k'$.
Clearly, $\chi_\Delta(\F) \geq \chi_\Delta^f(\F)$ and $\chi_*(\F) \geq \chi_*^f(\F) \geq \chi_\Delta^f(\F)$ for every graph class $\F$.

Norin, Scott, and Wood \cite{nsw} proved that $\chi_\Delta^f(\F)=\omega_\Delta(\F)$ for every minor-closed family $\F$.
Hence Theorem \ref{main_intro} implies that $\chi_\Delta(\F)=\chi_\Delta^f(\F)$ for every minor-closed family $\F$.
In fact, Norin, Scott, and Wood \cite{nsw} also proved that $\chi_*^f(\F)=\omega_\Delta(\F)$ for every minor-closed family $\F$.
Our Theorem \ref{main_intro} can give a proof for the same result independent from the original proof in \cite{nsw}, but we omit the details.
In summary, we have the following corollary.

\begin{corollary} \label{cor_long}
For every proper minor-closed family $\F$, $$\chi_*(\F) \geq \chi_*^f(\F) = \chi_\Delta^f(\F) = \chi_\Delta(\F) = \omega_\Delta(\F) = \min_{H \not \in \F}\td(H)-1 \geq \frac{1}{3} \cdot \chi_*(\F).$$
\end{corollary}

Note that the first inequality in Corollary \ref{cor_long} is an equality for some minor-closed families and is strict for some minor-closed families.
For example, if $\F$ is the set of $K_{t+1}$-minor free graphs, then $\chi_\Delta(\F)=t$ \cite{ekkos} and $\chi_*(\F)=t$ \cite{demw}\footnote{In fact, in an much earlier paper \cite{dn}, Dvo\v{r}\'{a}k and Norin announced that a much stronger result that implies $\chi_*(\F)=t$ will be proved in a forthcoming paper.}; if $\F$ is the set of graphs that are embeddable in a fixed surface, then $\chi_\Delta(\F)=3$ \cite{a} and $\chi_*(\F)=4$ \cite{dn}.

Hadwiger's conjecture can be restated as: if a graph does not have a proper $t$-coloring, then $K_{t+1}$ is its minor.
So it states that $K_{t+1}$ is an unavoidable minor of every non-properly $t$-colorable graph.
Theorem \ref{main_intro} implies an analogous result for partitioning an infinite graph\footnote{All graphs are finite in this paper, unless otherwise specified.} into (finite or infinite) induced subgraphs with uniformly bounded maximum degree.
Let $G$ be an infinite graph, and let $h$ be a positive integer.
Let ${\mathcal G}_{\rm ind}$ be the set of all finite induced subgraphs of $G$, and let ${\mathcal G}_{\rm min}$ be the set of all minors of graphs in ${\mathcal G}_{\rm ind}$.
So ${\mathcal G}_{\rm min}$ is a minor-closed family containing ${\mathcal G}_{\rm ind}$.
The standard compactness argument shows that if $d$ is an integer such that every graph in ${\mathcal G}_{\rm ind}$ has an $h$-coloring with defect $d$, then so does $G$, so $G$ can be partitioned into $h$ (finite or infinite) induced subgraphs with maximum degree at most $d$.
Hence, if there exists no integer $d$ such that $G$ can be partitioned into $h$ induced subgraphs with maximum degree at most $d$, then $h< \chi_\Delta({\mathcal G}_{\rm ind}) \leq \chi_\Delta({\mathcal G}_{\rm min}) = \omega_\Delta({\mathcal G}_{\rm min})$ by Theorem \ref{main_intro}, so $\CT_{h+1,k}$ is a minor of $G$ for every positive integer $k$.
As every (finite) graph with connected tree-depth at most $h+1$ is a minor of $\CT_{h+1,k}$ for some integer $k$, we obtain the following corollary.

\begin{corollary} \label{infinite_intro}
For every positive integer $h$ and every infinite graph $G$, if there exists no integer $d$ such that $G$ can be partitioned into $h$ (finite or infinite) induced subgraphs with maximum degree at most $d$, then every (finite) graph with connected tree-depth at most $h+1$ is a minor of $G$. 
\end{corollary}

Note that Corollary \ref{infinite_intro} is optimal since the closure of the balanced $\aleph_0$-ary tree with height $h+1$ cannot be partitioned into $h$ induced subgraphs with bounded maximum degree and every finite minor of it has connected tree-depth at most $h+1$.

\subsection{Proof sketch} \label{subsec:sketch}

Now we sketch the proof of Theorem \ref{main_intro}.
As discussed earlier, it suffices to prove the case $H=\CT_{h,k}$ for Conjecture \ref{defect_conj}, for any fixed integers $h$ and $k$.

Let $G$ be a $\CT_{h,k}$-minor free graph.
We shall construct a sequence $G_1,G_2,...,G_t$ of graphs with $G_1=G$ such that $G_t$ has a bounded number of vertices, and each $G_{i+1}$ is a minor of $G_i$.
Note that each vertex $v$ of $G_i$ corresponds to a connected subgraph of $G$ contracted into $v$.
When such a subgraph of $G$ corresponding to a vertex of $G_i$ has only 1 vertex, we also treat this vertex of $G_i$ as a vertex of $G$.
Then we construct a coloring by first coloring all vertices in $V(G_t) \cap V(G)$ with color 1, and then when all vertices in $V(G_{i+1}) \cap V(G)$ are colored, we extend the coloring to $V(G_i) \cap V(G)$ by coloring each vertex $v$ in $V(G_i) \cap V(G)-V(G_{i+1})$ greedily by using the smallest color that does not appear on colored neighbors of $v$ with large degree.
Since $G_1=G$, all vertices are colored eventually.
There are two remaining tasks: one is to show that this coloring has bounded defect, and the other is to show that at most $h-1$ colors are used.

To show this coloring has bounded defect, we will ensure that when a vertex is about to be colored, it only has bounded degree in the current graph.
In other words, we will only contract subgraphs of $G_i$ induced by vertices with bounded degree to obtain $G_{i+1}$.
Once we have this property, showing that the aforementioned greedy coloring has bounded defect is relatively easy.

Showing at most $h-1$ colors are used is more challenging.
The key idea is that we keep some information to ensure that when we are about to color a vertex $v \in V(G_i) \cap V(G)$, the set $S$ of colored neighbors of $v$ with large degree only use at most $h-1-j$ colors for some positive integer $j$ such that each vertex in $S$ together with $G-(V(G_i) \cap V(G))$ contains a $\CT_{j+1,k}$-minor in $G$.
As $G$ is $\CT_{h,k}$-minor free, $1 \leq j \leq h-2$, so $1 \leq h-1-j \leq h-2$, and hence there is always an available color for $v$ and the upper bound $h-1-j$ is positive.
It would be convenient if we think that we put a hyperedge on $S \cup \{v\}$ and label it with $j$, and direct the edges from $S$ to $v$ to indicate that $v$ should avoid the colors used by vertices in $S$.

So we need to ensure that such a $\CT_{j+1,k}$-minor can be constructed when we define $G_{i+1}$ from $G_i$.
Note that if there are $k$ vertices in $V(G_i) \cap V(G)$ with the same neighborhood $S$ in $G_i$, then we can delete those $k$ vertices from $G_i$ to create $G_{i+1}$ to obtain the case for $j=1$.
Similarly, if there are many distinct hyperedges whose common intersection equals their pairwise intersections, then we can assemble many $\CT_{j+1,k}$-minors corresponding to those hyperedges to construct a $\CT_{j+2,k}$-minor, put a hyperedge on the common intersection with label $j+1$, and delete the vertices not in the common intersection.
In other words, if we are able to find a large set of ``homogeneous structures'', then we can upgrade hyperedges with label $j$ to hyperedges with label $j+1$ to get flexibility for coloring and construct the corresponding $\CT_{j+2,k}$-minors.
See Figure \ref{fig_upgrade} for an example.

\begin{figure} 
	\begin{picture}(100,400) (-45,-220)
		\thinlines
		\multiput(5,120)(5,0){80}{\line(1,0){3}}
		\put(-20,100){{$V(G)-V(G_i)$}}
		\put(-20,140){{$V(G_i)$}}

		\thicklines
		\put(160,200){\circle*{5}} \put(156,190){$x$} 
		\put(240,200){\circle*{5}} \put(236,190){$y$} 
		\put(100,150){\circle*{5}} \put(96,155){$v_1$} 
		\put(200,150){\circle*{5}} \put(196,155){$v_2$} 
		\put(300,150){\circle*{5}} \put(296,155){$v_3$} 

		\thinlines
		\put(80,130){\line(0,1){90}}
		\put(80,220){\line(1,0){180}}
		\put(260,220){\line(0,-1){40}}
		\put(260,180){\line(-1,0){140}}
		\put(120,180){\line(0,-1){50}}
		\put(120,130){\line(-1,0){40}}
		\put(60,180){$S_1$}

		\put(150,130){\line(0,1){80}}
		\put(150,210){\line(1,0){100}}
		\put(250,210){\line(0,-1){80}}
		\put(250,130){\line(-1,0){100}}
		\put(160,140){$S_2$}

		\put(320,130){\line(0,1){85}}
		\put(320,215){\line(-1,0){190}}
		\put(130,215){\line(0,-1){30}}
		\put(130,185){\line(1,0){150}}
		\put(280,185){\line(0,-1){55}}
		\put(280,130){\line(1,0){40}}
		\put(330,180){$S_3$}

		\thicklines

		\put(70,50){$H_1$}

		\put(30,90){\circle*{5}} 
		\put(20,70){\circle*{5}} 
		\put(30,70){\circle*{5}} 
		\put(40,70){\circle*{5}} 
		\put(20,70){\line(1,2){10}}
		\put(30,70){\line(0,1){20}}
		\put(40,70){\line(-1,2){10}}	
	
		\put(60,90){\circle*{5}} 
		\put(50,70){\circle*{5}} 
		\put(60,70){\circle*{5}} 
		\put(70,70){\circle*{5}} 
		\put(50,70){\line(1,2){10}}
		\put(60,70){\line(0,1){20}}
		\put(70,70){\line(-1,2){10}}

		\put(90,90){\circle*{5}} 
		\put(80,70){\circle*{5}} 
		\put(90,70){\circle*{5}} 
		\put(100,70){\circle*{5}} 
		\put(80,70){\line(1,2){10}}
		\put(90,70){\line(0,1){20}}
		\put(100,70){\line(-1,2){10}}

		\put(120,90){\circle*{5}} 
		\put(110,70){\circle*{5}} 
		\put(120,70){\circle*{5}} 
		\put(130,70){\circle*{5}} 
		\put(110,70){\line(1,2){10}}
		\put(120,70){\line(0,1){20}}
		\put(130,70){\line(-1,2){10}}

		\multiput(10,60)(3,0){42}{\line(1,0){2}}
		\multiput(136,60)(0,3){12}{\line(0,1){2}}
		\multiput(10,95)(3,0){42}{\line(1,0){2}}
		\multiput(10,60)(0,3){12}{\line(0,1){2}}

		\qbezier(60,95)(90,130)(90,130)
		\qbezier(80,95)(100,130)(100,130)
		\qbezier(110,95)(110,130)(110,130)

		\put(200,50){$H_2$}
		
		\put(160,90){\circle*{5}} 
		\put(150,70){\circle*{5}} 
		\put(160,70){\circle*{5}} 
		\put(170,70){\circle*{5}} 
		\put(150,70){\line(1,2){10}}
		\put(160,70){\line(0,1){20}}
		\put(170,70){\line(-1,2){10}}
	
		\put(190,90){\circle*{5}} 
		\put(180,70){\circle*{5}} 
		\put(190,70){\circle*{5}} 
		\put(200,70){\circle*{5}} 
		\put(180,70){\line(1,2){10}}
		\put(190,70){\line(0,1){20}}
		\put(200,70){\line(-1,2){10}}

		\put(220,90){\circle*{5}} 
		\put(210,70){\circle*{5}} 
		\put(220,70){\circle*{5}} 
		\put(230,70){\circle*{5}} 
		\put(210,70){\line(1,2){10}}
		\put(220,70){\line(0,1){20}}
		\put(230,70){\line(-1,2){10}}

		\put(250,90){\circle*{5}} 
		\put(240,70){\circle*{5}} 
		\put(250,70){\circle*{5}} 
		\put(260,70){\circle*{5}} 
		\put(240,70){\line(1,2){10}}
		\put(250,70){\line(0,1){20}}
		\put(260,70){\line(-1,2){10}}

		\multiput(142,60)(3,0){42}{\line(1,0){2}}
		\multiput(268,60)(0,3){12}{\line(0,1){2}}
		\multiput(142,95)(3,0){42}{\line(1,0){2}}
		\multiput(142,60)(0,3){12}{\line(0,1){2}}

		\qbezier(170,95)(170,130)(170,130)
		\qbezier(200,95)(200,130)(200,130)
		\qbezier(230,95)(230,130)(230,130)

		\put(330,50){$H_3$}

		\put(290,90){\circle*{5}} 
		\put(280,70){\circle*{5}} 
		\put(290,70){\circle*{5}} 
		\put(300,70){\circle*{5}} 
		\put(280,70){\line(1,2){10}}
		\put(290,70){\line(0,1){20}}
		\put(300,70){\line(-1,2){10}}

		\put(320,90){\circle*{5}} 
		\put(310,70){\circle*{5}} 
		\put(320,70){\circle*{5}} 
		\put(330,70){\circle*{5}} 
		\put(310,70){\line(1,2){10}}
		\put(320,70){\line(0,1){20}}
		\put(330,70){\line(-1,2){10}}

		\put(350,90){\circle*{5}} 
		\put(340,70){\circle*{5}} 
		\put(350,70){\circle*{5}} 
		\put(360,70){\circle*{5}} 
		\put(340,70){\line(1,2){10}}
		\put(350,70){\line(0,1){20}}
		\put(360,70){\line(-1,2){10}}

		\put(380,90){\circle*{5}} 
		\put(370,70){\circle*{5}} 
		\put(380,70){\circle*{5}} 
		\put(390,70){\circle*{5}} 
		\put(370,70){\line(1,2){10}}
		\put(380,70){\line(0,1){20}}
		\put(390,70){\line(-1,2){10}}

		\multiput(273,60)(3,0){42}{\line(1,0){2}}
		\multiput(399,60)(0,3){12}{\line(0,1){2}}
		\multiput(273,95)(3,0){42}{\line(1,0){2}}
		\multiput(273,60)(0,3){12}{\line(0,1){2}}

		\qbezier(300,95)(290,130)(290,130)
		\qbezier(330,95)(300,130)(300,130)
		\qbezier(360,95)(310,130)(310,130)

		\thinlines
		\multiput(5,-100)(5,0){80}{\line(1,0){3}}
		\put(-20,-120){{$V(G)-V(G_{i+1})$}}
		\put(-20,-80){{$V(G_{i+1})$}}

		\thicklines
		\put(160,-70){\circle*{5}} \put(156,-80){$x$} 
		\put(240,-70){\circle*{5}} \put(236,-80){$y$} 
		\put(100,-120){\circle*{5}} \put(96,-115){$v_1$} 
		\put(200,-120){\circle*{5}} \put(196,-115){$v_2$} 
		\put(300,-120){\circle*{5}} \put(296,-115){$v_3$} 

		\thinlines
		\put(150,-85){\line(0,1){25}}
		\put(150,-60){\line(1,0){100}}
		\put(250,-60){\line(0,-1){25}}
		\put(250,-85){\line(-1,0){100}}
		\put(200,-55){$S'$}

		\thicklines

		\put(70,-225){$H_1$}

		\put(20,-200){\circle*{5}} \put(15,-208){$x_1$} 
		\put(30,-200){\circle*{5}} \put(28,-208){$y_1$}  
		\qbezier(20,-200)(100,-120)(100,-120)
		\qbezier(30,-200)(100,-120)(100,-120)
	
		\put(60,-180){\circle*{5}} 
		\put(50,-200){\circle*{5}} 
		\put(60,-200){\circle*{5}} 
		\put(70,-200){\circle*{5}} 
		\put(50,-200){\line(1,2){10}}
		\put(60,-200){\line(0,1){20}}
		\put(70,-200){\line(-1,2){10}}

		\put(90,-180){\circle*{5}} 
		\put(80,-200){\circle*{5}} 
		\put(90,-200){\circle*{5}} 
		\put(100,-200){\circle*{5}} 
		\put(80,-200){\line(1,2){10}}
		\put(90,-200){\line(0,1){20}}
		\put(100,-200){\line(-1,2){10}}

		\put(120,-180){\circle*{5}} 
		\put(110,-200){\circle*{5}} 
		\put(120,-200){\circle*{5}} 
		\put(130,-200){\circle*{5}} 
		\put(110,-200){\line(1,2){10}}
		\put(120,-200){\line(0,1){20}}
		\put(130,-200){\line(-1,2){10}}

		\multiput(10,-213)(3,0){42}{\line(1,0){2}}
		\multiput(136,-213)(0,3){13}{\line(0,1){2}}
		\multiput(10,-175)(3,0){42}{\line(1,0){2}}
		\multiput(10,-213)(0,3){13}{\line(0,1){2}}

		\qbezier(60,-180)(100,-120)(100,-120)
		\qbezier(50,-200)(50,-140)(100,-120)
		\qbezier(60,-200)(100,-120)(100,-120)
		\qbezier(70,-200)(100,-120)(100,-120)

		\qbezier(90,-180)(100,-120)(100,-120)
		\qbezier(80,-200)(100,-120)(100,-120)
		\qbezier(90,-200)(100,-150)(100,-120)
		\qbezier(100,-200)(100,-120)(100,-120)

		\qbezier(120,-180)(100,-120)(100,-120)
		\qbezier(110,-200)(100,-120)(100,-120)
		\qbezier(120,-200)(100,-120)(100,-120)
		\qbezier(130,-200)(120,-120)(100,-120)
		
		\qbezier(130,-175)(160,-85)(160,-85)
		\qbezier(132,-175)(170,-85)(170,-85)
		\qbezier(134,-175)(180,-85)(180,-85)

		\put(200,-225){$H_2$}
		
		\put(150,-200){\circle*{5}} \put(145,-208){$x_2$}  
		\put(160,-200){\circle*{5}} \put(158,-208){$y_2$}  
		\qbezier(150,-200)(200,-120)(200,-120)
		\qbezier(160,-200)(200,-120)(200,-120)
	
		\put(190,-180){\circle*{5}} 
		\put(180,-200){\circle*{5}} 
		\put(190,-200){\circle*{5}} 
		\put(200,-200){\circle*{5}} 
		\put(180,-200){\line(1,2){10}}
		\put(190,-200){\line(0,1){20}}
		\put(200,-200){\line(-1,2){10}}

		\put(220,-180){\circle*{5}} 
		\put(210,-200){\circle*{5}} 
		\put(220,-200){\circle*{5}} 
		\put(230,-200){\circle*{5}} 
		\put(210,-200){\line(1,2){10}}
		\put(220,-200){\line(0,1){20}}
		\put(230,-200){\line(-1,2){10}}

		\put(250,-180){\circle*{5}} 
		\put(240,-200){\circle*{5}} 
		\put(250,-200){\circle*{5}} 
		\put(260,-200){\circle*{5}} 
		\put(240,-200){\line(1,2){10}}
		\put(250,-200){\line(0,1){20}}
		\put(260,-200){\line(-1,2){10}}

		\multiput(142,-213)(3,0){42}{\line(1,0){2}}
		\multiput(268,-213)(0,3){13}{\line(0,1){2}}
		\multiput(142,-175)(3,0){42}{\line(1,0){2}}
		\multiput(142,-213)(0,3){13}{\line(0,1){2}}

		\qbezier(190,-180)(200,-120)(200,-120)
		\qbezier(180,-200)(200,-120)(200,-120)
		\qbezier(190,-200)(197,-180)(200,-120)
		\qbezier(200,-200)(200,-120)(200,-120)

		\qbezier(220,-180)(200,-120)(200,-120)
		\qbezier(210,-200)(200,-120)(200,-120)
		\qbezier(220,-200)(200,-120)(200,-120)
		\qbezier(230,-200)(210,-120)(200,-120)

		\qbezier(250,-180)(200,-120)(200,-120)
		\qbezier(240,-200)(200,-120)(200,-120)
		\qbezier(250,-200)(200,-120)(200,-120)
		\qbezier(260,-200)(245,-120)(200,-120)

		\qbezier(145,-175)(185,-85)(185,-85)
		\qbezier(150,-175)(195,-85)(195,-85)
		\qbezier(155,-175)(205,-85)(205,-85)

		\put(330,-225){$H_3$}

		\put(280,-200){\circle*{5}} \put(275,-208){$x_3$}  
		\put(290,-200){\circle*{5}} \put(288,-208){$y_3$}  
		\qbezier(280,-200)(300,-120)(300,-120)
		\qbezier(290,-200)(300,-120)(300,-120)

		\put(320,-180){\circle*{5}} 
		\put(310,-200){\circle*{5}} 
		\put(320,-200){\circle*{5}} 
		\put(330,-200){\circle*{5}} 
		\put(310,-200){\line(1,2){10}}
		\put(320,-200){\line(0,1){20}}
		\put(330,-200){\line(-1,2){10}}

		\put(350,-180){\circle*{5}} 
		\put(340,-200){\circle*{5}} 
		\put(350,-200){\circle*{5}} 
		\put(360,-200){\circle*{5}} 
		\put(340,-200){\line(1,2){10}}
		\put(350,-200){\line(0,1){20}}
		\put(360,-200){\line(-1,2){10}}

		\put(380,-180){\circle*{5}} 
		\put(370,-200){\circle*{5}} 
		\put(380,-200){\circle*{5}} 
		\put(390,-200){\circle*{5}} 
		\put(370,-200){\line(1,2){10}}
		\put(380,-200){\line(0,1){20}}
		\put(390,-200){\line(-1,2){10}}

		\multiput(273,-213)(3,0){42}{\line(1,0){2}}
		\multiput(399,-213)(0,3){13}{\line(0,1){2}}
		\multiput(273,-175)(3,0){42}{\line(1,0){2}}
		\multiput(273,-213)(0,3){13}{\line(0,1){2}}

		\qbezier(320,-180)(300,-120)(300,-120)
		\qbezier(310,-200)(300,-120)(300,-120)
		\qbezier(320,-200)(300,-120)(300,-120)
		\qbezier(330,-200)(310,-120)(300,-120)

		\qbezier(350,-180)(300,-120)(300,-120)
		\qbezier(340,-200)(300,-120)(300,-120)
		\qbezier(350,-200)(300,-120)(300,-120)
		\qbezier(360,-200)(340,-120)(300,-120)

		\qbezier(380,-180)(310,-120)(300,-120)
		\qbezier(370,-200)(310,-120)(300,-120)
		\qbezier(380,-200)(310,-120)(300,-120)
		\qbezier(390,-200)(385,-120)(300,-120)

		\qbezier(275,-175)(220,-85)(220,-85)
		\qbezier(278,-175)(230,-85)(230,-85)
		\qbezier(281,-175)(240,-85)(240,-85)

	\end{picture}
	\caption{The upper picture shows three hyperedges $S_1,S_2,S_3$ with label $2$ and three $4\CT_{2,3}$-minors $H_1,H_2,H_3$ in $G-V(G_i)$, where for every $q \in [3]$, every vertex in $H_q$ is adjacent to every vertex in $S_q=\{x,y,v_q\}$, so each vertex in $S_q$ together with $H_q$ gives a $\CT_{3,3}$-minor. The lower picture shows that, if we construct $G_{i+1}$ by ``removing'' $\bigcup_{q=1}^3S_q-\bigcap_{q=1}^3S_q=\{v_1,v_2,v_3\}$, then we obtain a hyperedge $S'=\bigcap_{q=1}^3S_q=\{x,y\}$ such that for each vertex $x$ (or $y$, respectively) in $S'$, we can contract the connected subgraph induced by $\{x,x_q:q \in [3]\}$ (or $\{y,y_q: q \in [3]\}$, respectively), so that it gives a $\CT_{4,3}$-minor together with $\bigcup_{q=1}^3(H_q \cup \{v_q\})$, where each $\{x_q,y_q\}$ is a set of two vertices contained in a $\CT_{2,3}$ in $H_q$. Note that this picture is a simplification of what we will really do in the proof. In particular, each hyperedge contains some directed edges and extra vertices with certain properties that the vertices $x,y,v_q$ in the picture do not have.} \label{fig_upgrade}
\end{figure}

However, there are some technical issues.
Since some vertex in those homogeneous structures might also belong to other hyperedges that are not in a homogeneous structure, deleting such a vertex loses too much information and cannot ensure the number of colors used in those hyperedges.
To deal with it, we should contract subgraphs instead of deleting vertices and somehow keep the information for those hyperedges that cannot be upgraded.
Figuring out what the ``right information'' (about those hyperedges, labels, directed edges, and neighborhoods) that should be kept for such a contraction is the main challenge behind the proof.

Now we give more details about the above strategy.
We will construct a sequence of graphs $G_1,G_2,...$ while keeping the information about hyperedges, labels, directed edges, and those $\CT_{j+1,k}$-minors.
This is essentially the intuition behind the strong defective elimination scheme defined in Sections \ref{sec:elimination} and \ref{sec:strong}.
It is a sequence whose each entry is a tuple of the form $(G_i,\M_i,\E_i,D_i,\A_i, \allowbreak \A_i')$, where $G_i$ is a graph mentioned above, $\M_i$ indicates what subgraphs of $G$ are contracted into vertices of $G_i$, $\E_i$ is the set of the aforementioned hyperedges with labels, $D_i$ is the set of the aforementioned directed edges, and $\A_i$ and $\A_i'$ are partial information for those $\CT_{j+1,k}$-minors.
To make sure $(G_{i+1},\M_{i+1},\E_{i+1},D_{i+1},\A_{i+1},\A_{i+1}')$ can be constructed from $(G_i,\M_i,\E_i,D_i,\A_i,\A_i')$, we should maintain many properties for $(G_i,\M_i,\E_i,D_i,\A_i,\A_i')$ (i.e.\ conditions (D1)-(D12) stated in Sections \ref{sec:elimination} and \ref{sec:strong}) during the construction, and those conditions are exactly the ``right information'' mentioned above.
Note that the precise definition of those terms might be slightly different from the ones mentioned in the above proof sketch.
But roughly speaking, (D1)-(D3) describe relationships between $G_i$ and $G_{i+1}$ by using standard language about minor models; (D4) describes properties about the set $D_i$ of directed edges; (D5) describes properties about the set $\E_i$ of hyperedges with labels; (D6) and (D7) describe further properties of $D_i$ and $\E_i$ to limit their complexities; (D9) and (D10) describe properties of the sets $\A_{i+1}$ and $\A_{i+1}'$, which are essentially the branch sets of those $\CT_{j,k}$-minors corresponding to the hyperedges in $\E_i$; (D11) describes relationships between the branch sets of the $\CT_{j,k}$-minors corresponding to different hyperedges.
All those properties are fairly straightforward to verify based on our construction from $G_i$ to $G_{i+1}$, even thought it is tedious to describe them precisely and rigorously.

Now we provide more details about how we construct $G_{i+1}$ from $G_i$ and explain what (D8) describes.
The key idea is to reduce the number of vertices in $G_i$ and try to upgrade hyperpedges to allow extra feasible colors for vertices in $V(G_i)-V(G_{i+1})$ mentioned in the above sketch.
Recall Figure \ref{fig_upgrade} for an example for upgrading a hyperedge.

Assume now we want to construct $G_{i+1}$ from $G_i$.
We can prove that there exist desired ``homogeneous structures''.
Roughly speaking, we can prove that there exist a set $U \subseteq V(G_i)$ with bounded size and pairwise disjoint sets $X_1,X_2,...,X_t$ of $V(G_i)$ (for some large integer $t$) such that each $X_j$ only contains vertices with bounded degree in $G_i$ and induces a connected subgraph with bounded size, and $U$ contains all large degree neighbors of the vertices in $\bigcup_{j=1}^tX_j$.
By a Ramsey-type argument, we may assume that the relationship between $U$ and $X_j$ is ``identical'' to the relationship between $U$ and $X_{j'}$, for any $j$ and $j'$.
That is, if we can find a hyperedge $S$ intersecting both $U$ and $X_j$ for some $j$, then for every other $j'$, we can find a hyperedge $S'$ intersecting both $U$ and $X_{j'}$ such that $S \cap U = S \cap U'$; and if a subset $U'$ of $U$ is the intersection of $U$ and the neighborhood of some vertex in $X_j$, then for every other $j'$, we can find a vertex in $X_{j'}$ such that $U'$ is also the intersection of $U$ and its neighborhood.
There are two cases for those $X_j$'s, based on whether $U$ also contains all small degree neighbors of $X_j$ or not.
Again a Ramsey-type argument allows us to assume that either $U$ actually contains all neighbors of $\bigcup_{j=1}^tX_j$, or for every $j$, some neighbor of some vertex in $X_j$ is not in $U$.

For the case that $U$ contains all neighbors of $\bigcup_{j=1}^tX_j$, we construct $G_{i+1}$ from $G_i$ by contracting $X_1$ into a single vertex $q_{i+1}$ and deleting all $X_2,X_3,...,X_t$.
Note that the ``identicalness'' of those $X_j$'s implies that $X_1$ is sufficient to represent all other $X_j$'s.
Detailed properties for this case are described in (D8h).
For the case that some neighbor of $X_j$ is not in $U$ for all $j$, we can show that $X_1$ itself is already ``homogeneous'' in the sense that if we can find a hyperedge $S$ intersecting both $U$ and $X_1$, then we can find many hyperedges $S'$ with pairwise disjoint $S'-U$ such that each $S'$ is contained in $X_1 \cup U$ and intersects both $U$ and $X_1$ with $S \cap U = S' \cap U$.
And we construct $G_{i+1}$ from $G_i$ by contracting $X_1$ into a single vertex $q_{i+1}$ and delete all edges between $X_1$ and its neighbors not in $U$.
Note that those deleted edges are between small degree vertices, so deleting them does not affect the defect of any coloring much, but it will give us convenience to simplify the proof.
Detailed properties for this case are described in (D8g).

We remark that in either case in our construction, for any hyperedge $S$ that involves with $V(G_i)-V(G_{i+1})$, there exist many hyperedges ``identical to $S$ with respect to $U$'' contained in $U \cup (V(G_i)-V(G_{i+1}))$, so we can ``remove'' them to upgrade $S$ as what we expected in Figure \ref{fig_upgrade}.

(D8a)-(D8e) describe properties about $U$ and $q_{i+1}$.
(D8f) and (D8i) describes the ``upgrade'' of the labels of the hyperedges mentioned earlier.

Properties (D8) can be verified from our construction of $G_{i+1}$ straightforward. 
The necessity of those properties naturally arise in the analysis for our procedure of greedy coloring.

\subsection{Organization and notations}

In Sections \ref{sec:elimination} and \ref{sec:color}, we show that a desired defective coloring follows from the existence of a strong defective elimination scheme.
In fact, we do not need the full strength of the strong defective elimination scheme to construct a coloring.
In particular, we do not need $\A_i$ and $\A_i'$ that record the information about the $\CT_{j+1,k}$-minors.
In Section \ref{sec:elimination}, we define the defective elimination scheme by only keeping required information.
In Section \ref{sec:color}, we show how to use it to construct a desired defective coloring.

In Section \ref{sec:homo}, we prove results that will be used for finding homogeneous structures in Section \ref{sec:strong}.
In Section \ref{sec:strong}, we define a strong elimination scheme by adding other required properties to a defective elimination scheme and prove the existence of a strong defective elimination scheme to complete the proof by inductively maintaining conditions (D1)-(D12).
The proof sketch of the two technical lemmas (Lemmas \ref{del} and \ref{contract}) was described in the previous section about how to construct $G_{i+1}$ from $G_i$.

To close this section, we introduce some notations that will be used in the rest of the paper.

All graphs in the rest of the paper are finite and simple.
Directed graphs are allowed to have different edges with the same ends but with different directions.

Let $G$ be a graph.
Let $S \subseteq V(G)$.
Let $\ell$ be a nonnegative real number.
We define $N_G^{\leq \ell}[S] = \{u \in V(G):$ there exists a path in $G$ with length at most $\ell$ between $u$ and some vertex in $S\}$.
We denote $N_G^{\leq 1}[S]$ by $N_G[S]$.
When $S$ consists of one vertex, say $v$, we denote $N_G^{\leq \ell}[S]$ by $N_G^{\leq \ell}[v]$ and denote $N_G^{\leq 1}[S]$ by $N_G[v]$.
For a set $S$ and a vertex $v$, we also define $N_G(S)=N_G[S]-S$ and $N_G(v)=N_G[v]-\{v\}$.

For every subset $S$ of $V(G)$, the subgraph of $G$ induced by $S$ is denoted by $G[S]$.

\section{Defective elimination schemes} \label{sec:elimination}

In this section, we define defective elimination schemes.
We encourage the reader to read Section \ref{subsec:sketch} in advance to get an intuition behind those terminologies.

Let $G$ be a graph.
Let $h \geq 3,k,r,d,N$ be positive integers.
Then a {\it $(G,h,k,r,d,N)$-defective elimination scheme} is a sequence $((G_i,\M_i,\E_i,D_i): i \in {\mathbb N})$ of tuples such that $(G_1,\M_1,\E_1,D_1)=(G,\{G[\{v\}]: v \in V(G)\},\emptyset, \emptyset)$, and for every $i \geq 1$, the following hold:
	\begin{enumerate}
		\item[(D1)] 
			\begin{itemize}
				\item $G_{i+1}$ is a graph.
				\item $\M_{i+1}$ is a collection $\{M_{i+1,v}: v \in V(G_{i+1})\}$ of disjoint connected induced subgraphs of $G$.
			\end{itemize}
			(If $v$ is a vertex of $G_{i+1}$ such that $M_{i+1,v}$ consists of only one vertex, then we also treat $v$ as a vertex of $G$ for simplicity.)
		\item[(D2)] $G_{i+1}$ is obtained from $G[\bigcup_{M \in \M_{i+1}}V(M)]$ by contracting each member of $\M_{i+1}$ into a vertex, deleting all resulting loops and parallel edges, and deleting possibly other edges such that 
			\begin{itemize}
				\item for every edge $uv$ of $G_{i+1}$, there exists $u'v' \in E(G_{i})$ such that $M_{i,u'} \subseteq M_{i+1,u}$ and $M_{i,v'} \subseteq M_{i+1,v}$, and
				\item for any distinct vertices $u,v$ of $G_{i+1}$, if $V(M_{i+1,u})$ is adjacent in $G$ to $V(M_{i+1,v})$ but $uv \not \in E(G_{i+1})$, then either $\lvert V(M_{i+1,u}) \rvert \geq 2$ or $\lvert V(M_{i+1,v}) \rvert \geq 2$.
			\end{itemize}

			(Note that these two statements imply that when $u,v$ are vertices in $V(G) \cap \bigcap_{\alpha=1}^{i+1}V(G_{\alpha})$, we have for every $\alpha \in [i+1]$, $uv \in E(G)$ if and only if $uv \in E(G_\alpha)$.)
		\item[(D3)] Exactly one of the following holds: 
			\begin{itemize}
				\item $\lvert V(G_{i}) \rvert \leq N$, $G_{i+1}=G_{i}$, $\M_{i+1}=\M_{i}$, $\E_{i+1}=\E_{i}$, and $D_{i+1}=D_{i}$. 
				\item $\lvert \M_{i+1} \rvert < \lvert \M_{i} \rvert$ and the following hold: 
					\begin{itemize}
						\item every member of $\M_{i}$ is either a member of $\M_{i+1}$, a subgraph of a member of $\M_{i+1}$, or disjoint from every member of $\M_{i+1}$, and 
						\item for every member $M$ of $\M_{i+1}-\M_{i}$, its vertex-set is a union of the vertex-sets of some members of $\M_{i}-\M_{i+1}$ such that $\{v \in V(G_{i}): V(M_{i,v}) \subseteq V(M)\}$ induces a connected subgraph of $G_{i}$. 
					\end{itemize} (That is, $\lvert V(G_{i+1}) \rvert < \lvert V(G_{i}) \rvert$, and $G_{i+1}$ is obtained from $G_{i}$ by contracting pairwise disjoint connected subgraphs and deleting vertices and edges. Moreover, if $v$ is a vertex in $V(G) \cap V(G_{i+1})$, then $v \in V(G) \cap \bigcap_{\alpha=1}^{i+1}V(G_\alpha)$.) 
			\end{itemize}
		\item[(D4)] $D_{i+1}$ is a set satisfying the following properties:
			\begin{itemize}
				\item $D_{i+1} \subseteq \{(u,v), (v,u): uv \in E(G_{i+1})\}$.

					(We call a vertex $v$ of $G_{i+1}$ a {\it head (with respect to $D_{i+1}$)} if $(u,v) \in D_{i+1}$ for some $u \in V(G_{i+1})$.)
				\item The digraph with vertex-set $V(G_{i+1})$ and edge-set $D_{i+1}$ has no directed cycle with length 2 and has no directed path with length 2.
				\item If there exists $(u',v') \in D_{i}$ with $M_{i,u'} \subseteq M_{i+1,u}$ and $M_{i,v'} \subseteq M_{i+1,v}$ for some edge $uv$ of $G_{i+1}$, then $(u,v) \in D_{i+1}$.
				\item If $(u,v) \in D_i$, then $\lvert V(M_{i+1,u}) \rvert=1$. 
			\end{itemize}
		\item[(D5)] $\E_{i+1}$ is a set whose each member is of the form $(S,j)$, where 
			\begin{itemize}
				\item $S$ is a subset of $V(G_{i+1})$ with $\lvert S \rvert \leq r+1$, 
				\item $j \in [h-2]$, and
				\item there exists a unique vertex $v$ in $S$ such that $(u,v) \in D_{i+1}$ for every $u \in S-\{v\}$.

					(We call this vertex $v$ the {\it sink} for $(S,j)$. Note that it is possible $S=\{v\}$, so $v$ is not necessarily a head with respect to $D_{i+1}$.)
			\end{itemize}
		\item[(D6)] For every $v \in V(G_{i+1})$, if either $\lvert V(M_{i+1,v}) \rvert \geq 2$, or $v$ is a head (with respect to $D_{i+1}$), or $v$ is the sink for some member of $\E_{i+1}$, then 
			\begin{itemize}
				\item[(D6a)] there are at most $N$ vertices $v'$ of $G_{i}$ such that $M_{i,v'} \subseteq M_{i+1,v}$, and
				\item[(D6b)] if $u$ is a vertex of $G_{i+1}$ such that either $\lvert V(M_{i+1,u}) \rvert \geq 2$, or $u$ is a head (with respect to $D_{i+1}$), or $u$ is the sink for some member of $\E_{i+1}$, then $uv \not \in E(G_{i+1})$.
			\end{itemize}
		\item[(D7)] For any $(S,j) \in \E_{i}$, if $S-\{v_S\} \subseteq V(G_{i+1})$ and $M_{i,v_S} \subseteq M_{i+1,v_{S'}}$ for some vertex $v_{S'}$ of $G_{i+1}$, where $v_S$ is the sink for $(S,j)$, then $(S',j) \in \E_{i+1}$, where $S' = \{v_{S'}\} \cup (N_{G_{i+1}}(v_{S'}) \cap S-\{v_S\})$. 
		\item[(D8)] If $G_{i+1} \neq G_{i}$, then there exist a vertex $q_{i+1} \in V(G_{i+1})$ and sets $U_{i+1}$ and $U_{i+1}^+$ with $U_{i+1} \subseteq U_{i+1}^+ \subseteq V(G_{i}) \cap V(G_{i+1}) \cap V(G)$ such that the following hold. 
			\begin{itemize}
				\item[(D8a)] For every $v \in V(G_{i})-V(G_{i+1})$ with $\lvert V(M_{i,v}) \rvert=1$, $\deg_{G_{i}}(v) \leq d$.
				\item[(D8b)] $U_{i+1} =\{u \in U_{i+1}^+: \deg_{G_{i+1}}(u)>d\}$.
				\item[(D8c)] $q_{i+1} \not \in U_{i+1}^+$ and $U_{i+1} \supseteq \{x \in N_{G_{i+1}}(q_{i+1}): \deg_{G_{i+1}}(x) >d\}$. 
				\item[(D8d)] For every $(S,j) \in \E_{i}$ with $M_{i,v_S} \subseteq M_{i+1,q_{i+1}}$, where $v_S$ is the sink for $(S,j)$, we have $\{x \in S-\{v_S\}: \deg_{G_{i}}(x) >d\} \subseteq U_{i+1}^+$. 
				\item[(D8e)] If $v \in V(G) \cap V(G_{i+1}) \cap N_{G_{i}}(v')-U_{i+1}^+$ with $\deg_{G_{i}}(v) \leq d$, for some $v' \in V(G_{i})$ with $M_{i,v'} \subseteq M_{i+1,q_{i+1}}$, then $U_{i+1} \supseteq \{x \in N_{G_{i+1}}(v) \cap V(G): \deg_{G_{i+1}}(x) >d\}$. 
				\item[(D8f)] $((U_{i+1} \cap N_{G_{i+1}}(q_{i+1})) \cup \{q_{i+1}\}, 1) \in \E_{i+1}$ and $q_{i+1}$ is its sink. 
				\item[(D8g)] If $\bigcup_{M \in \M_{i+1}}V(M) = \bigcup_{M \in \M_{i}}V(M)$, then 
					\begin{itemize}
						\item[(D8ga)] $q_{i+1} \in V(G_{i+1})-V(G_{i})$, 
						\item[(D8gb)] for every $M \in \M_{i}$, either $M \in \M_{i+1}$ or $M \subseteq M_{i+1,q_{i+1}}$, 
						\item[(D8gc)] if $uv$ is an edge of $G_{i}$ such that $u' \neq v'$ and $u'v' \not \in E(G_{i+1})$, where $u'$ and $v'$ are the vertices of $G_{i+1}$ such that $M_{i,u} \subseteq M_{i+1,u'}$ and $M_{i,v} \subseteq M_{i+1,v'}$, then $q_{i+1} \in \{u',v'\}$, and the vertex $x \in \{u',v'\}-\{q_{i+1}\}$ satisfies $x \in \{u,v\} \cap V(G)-U_{i+1}^+$, $\deg_{G_{i}}(x) \leq d$, and $x$ is not the sink for some member of $\E_{i}$, and 
						\item[(D8gd)] for every $(S,j) \in \E_{i}$ with $M_{i,v_S} \subseteq M_{i+1,q_{i+1}}$, where $v_S$ is the sink for $(S,j)$, 
							\begin{itemize}
								\item there exists $(S',j) \in \E_{i}$ with $\bigcup_{v \in S'-U_{i+1}^+}M_{i,v} \subseteq M_{i+1,q_{i+1}}$, $S \cap U_{i+1}^+=S' \cap U_{i+1}^+$ and $v_{S'} \not \in U_{i+1}^+$, where $v_{S'}$ is the sink for $(S',j)$, and 
								\item there exists a bijection $\iota$ from $S-(\{v_S\} \cup U_{i+1}^+)$ to $S'-(\{v_{S'}\} \cup U_{i+1}^+)$ such that for every $v \in S-(\{v_S\} \cup U_{i+1}^+)$, $N_{G_{i}}(v) \cap U_{i+1}^+ = N_{G_{i}}(\iota(v)) \cap U_{i+1}^+$.
							\end{itemize}
					\end{itemize}
				\item[(D8h)] If $\bigcup_{M \in \M_{i+1}}V(M) \neq \bigcup_{M \in \M_{i}}V(M)$, then 
							\begin{itemize}
								\item[(D8ha)] for every $M \in \M_{i}$, either $M \in \M_{i+1}$, or $M \subseteq M_{i+1,q_{i+1}}$, or $M$ is disjoint from all members of $\M_{i+1}$, 
								\item[(D8hb)] for every $v \in V(G_{i})$ with $M_{i,v}$ disjoint from all members of $\M_{i+1}$, we have $N_{G_{i}}(v) \subseteq (U_{i+1}^+ \cap N_{G_{i+1}}(q_{i+1})) \cup \{u \in V(G_{i}): M_{i,u}$ is disjoint from all members of $\M_{i+1}\}$, 
								\item[(D8hc)] every vertex in $N_{G_{i+1}}(q_{i+1})$ is in $V(G)$ and is not a head (with respect to $D_{i}$),
								\item[(D8hd)] if $uv$ is an edge of $G_{i}$ such that $M_{i,u} \subseteq M_{i+1,u'}$ and $M_{i,v} \subseteq M_{i+1,v'}$ for some distinct vertices $u',v'$ of $G_i$, then $u'v' \in E(G_{i+1})$,
								\item[(D8he)] $\lvert V(G_{i})-V(G_{i+1}) \rvert \leq N$,
								\item[(D8hf)] for every vertex $x \in V(G_{i}) \cap V(G)-\bigcup_{M \in \M_{i+1}}V(M)$, there exists a vertex $x' \in V(M_{i+1,q_{i+1}}) \cap V(G_{i}) \cap V(G)$ with $N_G(x) \cap U_{i+1} = N_G(x') \cap U_{i+1}$, and 
								\item[(D8hg)] for every $(S,j) \in \E_{i}$ with $V(M_{i,v_S}) \subseteq \bigcup_{M \in \M_{i}}V(M)-\bigcup_{M \in \M_{i+1}}V(M)$, where $v_S$ is the sink for $(S,j)$, there exists $(S',j) \in \E_{i}$ with $\lvert S' \rvert = \lvert S \rvert$ such that  
									\begin{itemize}
										\item $M_{i,v_{S'}} \subseteq M_{i+1,q_{i+1}}$, where $v_{S'}$ is the sink for $(S',j)$, 
										\item $S \cap V(G_{i+1}) - \{v_S\} = S' \cap V(G_{i+1}) - \{v_{S'}\}$, and
										\item there exists a bijection $\iota$ from $S -(V(G_{i+1}) \cup \{v_S\})$ to $S' \cap V(M_{i+1,q_{i+1}}) - \{v_{S'}\}$ such that for every $x \in S-(V(G_{i+1}) \cup \{v_S\})$, $N_G(x) \cap U_{i+1} = N_G(\iota(x)) \cap U_{i+1}$. 
									\end{itemize}
							\end{itemize}
				\item[(D8i)] For every $(S,j) \in \E_{i}$ with $M_{i,v} \cup M_{i,v_S} \subseteq M_{i+1,q_{i+1}}$ for some vertex $v \in S-\{v_S\}$, where $v_S$ is the sink for $(S,j)$, 
					\begin{itemize}
						\item[(D8ia)] if $\C$ is the multiset $\{N_G(x) \cap U_{i+1} \neq \emptyset: x \in S - (\{v_S\} \cup U_{i+1}^+)\}$, then for every function $f$ that maps each member $Z$ of $\C$ to a vertex in $Z$, we have $(S',j) \in \E_i$, where $S'=(S \cap U_{i+1}^+) \cup \{f(Z): Z \in \C\} \cup \{q_{i+1}\}$ and $q_{i+1}$ is the sink for $(S',j)$, and
						\item[(D8ib)] for every $u \in S-(\{v_S\} \cup U_{i+1}^+)$, if $\C_u$ is the multiset $\{N_G(x) \cap U_{i+1} \neq \emptyset: x \in S-(\{u,v_S\} \cup U_{i+1}^+)\}$, and either $S \cap U_{i+1}^+ \neq \emptyset$ or $\C_u \neq \emptyset$, then for every function $f$ that maps each member $Z$ of $\C_u$ to a vertex in $Z$, we have $(S',j+1) \in \E_i$, where $S' = (S \cap U_{i+1}^+) \cup \{f(Z): Z \in \C_u\} \cup \{q_{i+1}\}$ and $q_{i+1}$ is the sink for $(S',j+1)$.
					\end{itemize}
			\end{itemize}
	\end{enumerate}
For a positive integer $n$, an {\it $n$-$(G,h,k,r,d,N)$-defective elimination scheme} is a sequence $((G_i,\M_i, \allowbreak \E_i,D_i): i \in [n])$ such that that $(G_i,\M_i,\E_i,D_i)$ satisfies (D1)-(D8) for every $i \in [n]$, and $(G_1,\M_1,\E_1,D_1)=(G,\{G[\{v\}]: v \in V(G)\},\emptyset, \emptyset)$.

\section{Using a defective elimination scheme to color} \label{sec:color}

For a function $f$ and a subset $S$ of its domain, we define $f(S)=\{f(x): x \in S\}$.

\begin{lemma} \label{scheme_color}
For any positive integers $d$ and $N$, there exists a positive integer $d^*=d^*(d,N)$ such that the following holds.
Let $h \geq 3,k,r$ be positive integers.
Let $((G_i,\M_i,\E_i,D_i): i \in {\mathbb N})$ be a $(G,h,k,r,d,N)$-defective elimination scheme for some graph $G$.
Then for every $i \in {\mathbb N}$, there exists an $(h-1)$-coloring $c_i$ of $G[V(G_i) \cap V(G)]$ with defect $d^*$ such that 
	\begin{enumerate}
		\item for every $v \in V(G_i) \cap V(G)$ with $c_i(v) \neq 1$ and for every $\beta \in [c_i(v)-1]$, there exists $v' \in N_G(v) \cap V(G_i) \cap V(G_{i+1})$ with $\deg_{G_{i+1}}(v')>d$ and $c_i(v')=\beta$, and
		\item for every $(S,j) \in \E_i$, $\lvert c_i(S-\{v_S\}) \rvert \leq h-j-1$, where $v_S$ is the sink for $(S,j)$.
	\end{enumerate}
\end{lemma}

\begin{pf}
Let $d$ and $N$ be positive integers.
Define $d^*=2N+d$.

Let $h \geq 3,k,r$ be positive integers.
Let $((G_i,\M_i,\E_i,D_i): i \in {\mathbb N})$ be a $(G,h,k,r,d,N)$-defective elimination scheme for some graph $G$.
For every positive integer $i$ with $G_{i+1} \neq G_{i}$, let $U_{i+1},U_{i+1}^+$ be the subsets of $V(G_{i}) \cap V(G_{i+1}) \cap V(G)$ and the vertex $q_{i+1}$ of $G_{i+1}$ mentioned in (D8).
We shall prove this lemma by induction on $-i$.
Let $i$ be a fixed positive integer.

When $\lvert V(G_i) \rvert \leq N$, define $c_i$ to be the $1$-coloring of $G[V(G_i) \cap V(G)]$.  
Then $c_i$ has defect $N \leq d^*$.
Note that there exists no vertex with color non-equal to 1, so Statement 1 of this lemma holds.
And for every $(S,j) \in \E_i$, $j \in [h-2]$ by (D5), so $h-j-1 \geq 1$ and the lemma holds. 

By (D3), this proves the base case of the induction, and we may assume that $\lvert V(G_i) \rvert >N$ and there exists an $(h-1)$-coloring $c_{i+1}$ of $G[V(G_{i+1}) \cap V(G)]$ with defect $d^*$ such that the lemma holds for $c_{i+1}$. 

Since $\lvert V(G_i) \rvert >N$, by (D3), $\lvert V(G_i) \rvert > \lvert V(G_{i+1}) \rvert$.
In particular, $G_i \neq G_{i+1}$, so $q_{i+1},U_{i+1},U^+_{i+1}$ exist.

For simplicity of notations, we denote $q_{i+1},U_{i+1},U_{i+1}^+$ by $q,U,U^+$, respectively.

\medskip

\noindent{\bf Claim 1:} For every $v \in V(G) \cap V(G_i)-V(G_{i+1})$, we have $[h-1]-c_{i+1}(N_{G}(v) \cap U) \neq \emptyset$. 

\noindent{\bf Proof of Claim 1:}
Let $v \in V(G) \cap V(G_i)-V(G_{i+1})$. 
If $M_{i,v} \subseteq M_{i+1,q}$, then define $\bar{v}=v$; if $M_{i,v} \not \subseteq M_{i+1,q}$, then $\bigcup_{M \in \M_i}V(M) \neq \bigcup_{M \in \M_{i+1}}V(M)$ by (D8gb), so $M_{i,v}$ is disjoint from all members of $\M_{i+1}$ and $v \in V(G_i) \cap V(G)-\bigcup_{M \in \M_{i+1}}V(M)$ by (D8ha), and hence (D8hf) implies that there exists $\bar{v} \in V(M_{i+1,q}) \cap V(G_i) \cap V(G)$ with $N_G(v) \cap U=N_G(\bar{v}) \cap U$.
Note that since $\bar{v} \in V(G) \cap V(G_i)$, (D8) and (D2) imply $N_G(\bar{v}) \cap U = N_G(\bar{v}) \cap (U \cap V(G_{i}) \cap V(G)) = (N_G(\bar{v}) \cap V(G_{i}) \cap V(G)) \cap U = N_{G_{i}}(\bar{v}) \cap V(G) \cap U = N_{G_{i}}(\bar{v}) \cap U$.
Similarly, $N_G(v) \cap U = N_{G_{i}}(v) \cap U$.

Hence $\bar{v}$ is defined to be a vertex in $V(M_{i+1,q}) \cap V(G)$, and $N_{G_{i}}(v) \cap U = N_G(v) \cap U = N_G(\bar{v}) \cap U = N_{G_{i}}(\bar{v}) \cap U$. 

Suppose to the contrary that there exists $u \in N_G(\bar{v}) \cap U - (N_{G_{i+1}}(q) \cap U)$.
Since $N_G(\bar{v}) \cap U = N_{G_{i}}(\bar{v}) \cap U$, $u\bar{v} \in E(G_i)$.
Since $\bar{v} \in V(M_{i+1,q})$ and $u \in U \subseteq V(G_{i+1}) \cap V(G)$, (D8gc) and (D8hd) imply that $uq \in E(G_{i+1})$, contradicting to $u \not \in N_{G_{i+1}}(q) \cap U$.

Hence $N_G(\bar{v}) \cap U \subseteq N_{G_{i+1}}(q) \cap U$.
By (D8f), $((N_{G_{i+1}}(q) \cap U) \cup \{q\},1) \in \E_{i+1}$ and $q$ is its sink.
Since $c_{i+1}$ satisfies the lemma, we know $\lvert c_{i+1}(N_{G_{i+1}}(q) \cap U) \rvert \leq h-2$, so $[h-1] \not \subseteq c_{i+1}(N_{G_{i+1}}(q) \cap U)$. 
Hence $[h-1]-c_{i+1}(N_G(\bar{v}) \cap U) \supseteq [h-1]-c_{i+1}(N_{G_{i+1}}(q) \cap U) \neq \emptyset$. 
Since $N_G(v) \cap U = N_G(\bar{v}) \cap U$, we know $[h-1]-c_{i+1}(N_{G}(v) \cap U) = [h-1]-c_{i+1}(N_{G}(\bar{v}) \cap U) \neq \emptyset$.
$\Box$

\medskip

Define $c_i$ to be the $(h-1)$-coloring of $G[V(G_i) \cap V(G)]$ such that for every $v \in V(G_i) \cap V(G)$, 
	\begin{itemize}
		\item if $v \in V(G_{i+1})$, then $c_i(v)=c_{i+1}(v)$;
		\item otherwise, $c_i(v)$ is the minimum in $[h-1]-c_{i+1}(N_G(v) \cap U)$. 
	\end{itemize}
Note that $c_i$ is well-defined by Claim 1, and $c_i$ is an $(h-1)$-coloring of $G[V(G_i) \cap V(G)]$.

\medskip

\noindent{\bf Claim 2:} If $u$ is a vertex in $U$, and $u'$ is a vertex in $N_{G_i}(u) \cap V(G)$ with $c_i(u)=c_{i}(u')$, then $u' \in V(G_{i+1}) \cap V(G)$. 

\noindent{\bf Proof of Claim 2:}
Suppose to the contrary that $u' \not\in V(G_{i+1})$.
So $u' \in V(G) \cap V(G_i)-V(G_{i+1})$. 
Hence $c_i(u') \in [h-1]-c_{i+1}(N_G(u') \cap U)$.
So $u \in N_{G_i}(u') \cap U = N_G(u') \cap U$ (by (D2) and (D8)) and $c_i(u')=c_i(u)=c_{i+1}(u) \in c_{i+1}(N_G(u') \cap U)$, a contradiction.
$\Box$

\medskip

\noindent{\bf Claim 3:} $c_i$ has defect $d^*$.

\noindent{\bf Proof of Claim 3:}
Let $v \in V(G_i) \cap V(G)$. 
Let $W=\{x \in N_{G_i}(v) \cap V(G): c_i(v)=c_i(x)\}$.
To prove this claim, it suffices to show $\lvert W \rvert \leq d^*$. 
Since $d^* \geq d+2N$, we may assume $\deg_{G_i}(v) >d+2N$.
This implies $v \in V(G_{i+1})$ by (D8a).
So $c_i(v)=c_{i+1}(v)$.

We first assume $v \in U$.
By Claim 2, for every $u' \in W$, we have $u' \in V(G_{i+1}) \cap V(G)$, so $c_{i+1}(u')=c_i(u')=c_i(v)=c_{i+1}(v)$. 
So $W \subseteq \{x \in N_{G_{i+1}}(v) \cap V(G): c_{i+1}(v)=c_{i+1}(x)\}$ by (D2).
Hence $\lvert W \rvert \leq d^*$ since $c_{i+1}$ has defect $d^*$.

So we may assume $v \not \in U$. 
By (D8ga), if $\bigcup_{M \in \M_i}V(M)=\bigcup_{M \in \M_{i+1}}V(M)$, then $q \not \in V(G_i)$, so $q \neq v$.
Since $\deg_{G_i}(v)>d$ and $v \in V(G_i) \cap V(G_{i+1}) \cap V(G)$, by (D8gb) and (D8gc), if $\bigcup_{M \in \M_i}V(M)=\bigcup_{M \in \M_{i+1}}V(M)$, then $\deg_{G_{i+1}}(v) \geq \deg_{G_i}(v)- \lvert \{y \in V(G_i): M_{i,y} \subseteq M_{i+1,q}\} \rvert$. 
So by (D6a), if $\bigcup_{M \in \M_i}V(M)=\bigcup_{M \in \M_{i+1}}V(M)$, then $\deg_{G_{i+1}}(v) \geq \deg_{G_i}(v)- N>d$. 
If $\bigcup_{M \in \M_i}V(M) \neq \bigcup_{M \in \M_{i+1}}V(M)$, then by (D8ha) and (D8hd), $\deg_{G_{i+1}}(v) \geq \deg_{G_i}(v)- \lvert V(G_i)-V(G_{i+1}) \rvert$, so $\deg_{G_{i+1}}(v) \geq \deg_{G_i}(v)-N>d$ by (D8he). 

Hence $\deg_{G_{i+1}}(v)>d$. 
By (D8c), we know $U \supseteq N_{G_{i+1}}(q) \cap \{x \in V(G_{i+1}): \deg_{G_{i+1}}(x)>d\}$.
Since $v \in V(G_{i+1})-U$ and $\deg_{G_{i+1}}(v)>d$, we know $v \not \in N_{G_{i+1}}(q)$.

Suppose to the contrary that $N_{G_i}(v) \cap V(G) \not \subseteq N_{G_{i+1}}(v) \cap V(G)$.

If $\bigcup_{M \in \M_i}V(M)=\bigcup_{M \in \M_{i+1}}V(M)$, then since $v \in V(G_i) \cap V(G_{i+1}) - N_{G_{i+1}}(q)$ and $\deg_{G_{i}}(v)>d$, (D8gc) implies that $N_{G_i}(v) \cap V(G) \cap V(M_{i+1,q}) = \emptyset$, so (D8gb) and (D8gc) imply $N_{G_i}(v) \cap V(G) \subseteq N_{G_{i+1}}(v) \cap V(G)$, a contradiction.

So $\bigcup_{M \in \M_i}V(M) \neq \bigcup_{M \in \M_{i+1}}V(M)$.
Since $v \in V(G_i) \cap V(G_{i+1}) - N_{G_{i+1}}(q)$, we have $N_{G_i}(v) \cap V(G) \cap V(M_{i+1,q}) = \emptyset$ by (D8hd).
If there exists $x \in N_{G_i}(v) \cap V(G)$ with $M_{i,x}$ disjoint from all members of $\M_{i+1}$, then $v \in N_{G_i}(x)$, so (D8hb) implies either $v \in N_{G_{i+1}}(q)$ or $v \not \in V(G_{i+1})$, a contradiction. 
So $N_{G_i}(v) \cap V(G) \subseteq N_{G_{i+1}}(v) \cap V(G)$ by (D8ha) and (D8hd), a contradiction.

Therefore, $N_{G_i}(v) \cap V(G) \subseteq N_{G_{i+1}}(v) \cap V(G)$.
So $W \subseteq \{x \in N_{G_{i+1}}(v) \cap V(G) \cap V(G_i): c_{i}(x)=c_{i}(v)\} \subseteq \{x \in N_{G_{i+1}}(v) \cap V(G): c_{i+1}(x)=c_{i+1}(v)\}$ has size at most $d^*$ since $c_{i+1}$ has defect $d^*$.
$\Box$

\medskip

\noindent{\bf Claim 4:} For every $v \in V(G_i) \cap V(G)$ with $c_i(v) \neq 1$ and for every $\beta \in [c_i(v)-1]$, there exists $v' \in N_G(v) \cap V(G_i) \cap V(G_{i+1})$ with $\deg_{G_{i+1}}(v')>d$ and $c_i(v')=\beta$. 
Moreover, if $v \in V(G_i)-V(G_{i+1})$, then $v' \in U$.

\noindent{\bf Proof of Claim 4:}
Let $v \in V(G_i) \cap V(G)$ with $c_i(v) \neq 1$.
Let $\beta \in [c_i(v)-1]$.

We first assume $v \in V(G_{i+1})$.
So $c_{i+1}(v)=c_i(v) \neq 1$.
Since $c_{i+1}$ satisfies this lemma, there exists $v' \in N_G(v) \cap V(G_{i+1}) \cap V(G_{i+2})$ with $\deg_{G_{i+2}}(v')>d$ and $c_{i+1}(v')=\beta$.
In particular, $v' \in V(G) \cap V(G_i) \cap V(G_{i+1})$. 
Hence $v' \in N_G(v) \cap V(G_i) \cap V(G_{i+1})$ and $\deg_{G_{i+1}}(v') \geq \deg_{G_{i+2}}(v')>d$ by (D2).

So we may assume $v \in V(G_i)-V(G_{i+1})$.
Hence $c_i(v) = \min([h-1]-c_{i+1}(N_G(v) \cap U))$.
So there exists $v' \in N_G(v) \cap U$ such that $c_{i+1}(v')=\beta$.
Since $v' \in U$, $v' \in V(G_i) \cap V(G_{i+1})$ by (D8), so $v' \in N_G(v) \cap V(G_i) \cap V(G_{i+1})$ and $c_i(v')=c_{i+1}(v')=\beta$.
Since $v' \in U$, (D8b) implies that $\deg_{G_{i+1}}(v')>d$. 
This proves the claim.
$\Box$

\medskip

Hence by Claims 3 and 4, to prove this lemma, it suffices to prove that $\lvert c_i(S-\{v_S\}) \rvert \leq h-j-1$ for every $(S,j) \in \E_i$, where $v_S$ is the sink for $(S,j)$.

Suppose to the contrary that there exists $(S,j) \in \E_i$ such that $\lvert c_i(S-\{v_S\}) \rvert \geq h-j$. 
We further choose such an $(S,j)$ such that $(\bigcup_{x \in S}V(M_{i,x})) \cap V(M_{i+1,q}) \neq \emptyset$ if possible.
Let $v_S$ be the sink for $(S,j)$.

\medskip

\noindent{\bf Claim 5:} Either 
	\begin{itemize}
		\item $S-\{v_S\} \not \subseteq V(G_{i+1})$, or 
		\item $M_{i,v_S} \not \subseteq \bigcup_{M \in \M_{i+1}}M$, or
		\item $\bigcup_{M \in \M_i}V(M)=\bigcup_{M \in \M_{i+1}}V(M)$, $M_{i,v_S} \subseteq M_{i+1,q}$, and there exists $v \in S-(\{v_S\} \cup V(M_{i+1,q}) \cup U^+)$ with $\deg_{G_i}(v) \leq d$. 
	\end{itemize}

\noindent{\bf Proof of Claim 5:}
We assume that $S-\{v_S\} \subseteq V(G_{i+1})$ and $M_{i,v_S} \subseteq \bigcup_{M \in \M_{i+1}}M$. 
We shall prove that the third statement of this claim holds.
By (D3), there exists a vertex $v_{S'}$ of $G_{i+1}$ with $M_{i,v_S} \subseteq M_{i+1,v_{S'}}$.
Let $S'=\{v_{S'}\} \cup (N_{G_{i+1}}(v_{S'}) \cap (S-\{v_S\}))$.
By (D7), $(S',j) \in \E_{i+1}$.
Note that $v_{S'}$ is the sink for $(S',j)$ by (D4) and (D5).
Since $c_{i+1}$ satisfies the lemma, $\lvert c_{i+1}(S'-\{v_{S'}\}) \rvert \leq h-j-1$.
Hence $c_i(S-\{v_S\}) \neq c_{i+1}(S'-\{v_{S'}\})$.
Note that $S'-\{v_{S'}\} \subseteq S-\{v_S\}$.
So there exists $v \in (S-\{v_S\})-(S'-\{v_{S'}\})$ such that $c_i(v) \not \in c_{i+1}(S'-\{v_{S'}\})$.

In particular, $v \in (S-\{v_{S}\})-N_{G_{i+1}}(v_{S'})$.
Hence $S-\{v_{S}\} \not \subseteq N_{G_{i+1}}(v_{S'})$.
So $vv_S \in E(G_i)$ (by (D5)) and $vv_{S'} \not \in E(G_{i+1})$.
Hence $\bigcup_{M \in \M_i}V(M)=\bigcup_{M \in \M_{i+1}}V(M)$ by (D8hd).
By (D8gc), one of $v$ and $v_{S'}$ equals $q$.
By (D8ga), $q \in V(G_{i+1})-V(G_i)$, so $v \neq q$, and hence $q=v_{S'}$.
That is, $M_{i,v_S} \subseteq M_{i+1,v_{S'}}=M_{i+1,q}$.

Again by (D8gc), $v \in V(G)-U^+$ and $\deg_{G_i}(v) \leq d$. 
Since $v \in S-\{v_{S}\} \subseteq V(G_{i+1})$ and $v \neq q$, we know $v \not \in V(M_{i+1,q})$.
So $v \in S-(\{v_S\} \cup V(M_{i+1,q}) \cup U^+)$.
This proves the claim.
$\Box$

\medskip

\noindent{\bf Claim 6:} Either
	\begin{itemize}
		\item $(S-\{v_S\}) \cap V(M_{i+1,q}) \neq \emptyset$, or 
		\item $\bigcup_{M \in \M_i}V(M)=\bigcup_{M \in \M_{i+1}}V(M)$, $M_{i,v_S} \subseteq M_{i+1,q}$, and there exists $v \in S-(\{v_S\} \cup V(M_{i+1,q}) \cup U^+)$ with $\deg_{G_i}(v) \leq d$. 
	\end{itemize}

\noindent{\bf Proof of Claim 6:}
Suppose to the contrary that this claim does not hold.
In particular, $(S-\{v_S\}) \cap V(M_{i+1,q}) = \emptyset$.
So $M_{i,x} \not \subseteq M_{i+1,q}$ for every $x \in (S-\{v_S\})-V(G_{i+1})$.

Suppose $M_{i,v_S} \subseteq \bigcup_{M \in \M_{i+1}}M$.
By Claim 5, $S-\{v_S\} \not \subseteq V(G_{i+1})$. 
Hence there exists $v \in (S-\{v_S\})-V(G_{i+1})$.
So $M_{i,v} \not \in \M_{i+1}$ and $M_{i,v} \not \subseteq M_{i+1,q}$.
Hence $\bigcup_{M \in \M_i}V(M) \neq \bigcup_{M \in \M_{i+1}}V(M)$ by (D8gb).
So by (D8ha), $M_{i,x}$ is disjoint from all members of $\M_{i+1}$ for every $x \in (S-\{v_S\})-V(G_{i+1})$.
By the existence of $v$, $S-\{v_S\} \neq \emptyset$.
So by (D5), $v_S$ is a head with respect to $D_i$.
Hence by (D8hc), $v_S \not \in N_{G_{i+1}}(q)$.
Since $M_{i,v}$ is disjoint from all members of $\M_{i+1}$, by (D8hb), $v_S \in N_{G_i}(v) \subseteq N_{G_{i+1}}(q) \cup \{u \in V(G_i): M_{i,u}$ is disjoint from all members of $\M_{i+1}\}$.
Hence $M_{i,v_S}$ is disjoint from all members of $\M_{i+1}$, contradicting $M_{i,v_S} \subseteq \bigcup_{M \in \M_{i+1}}M$.

So $M_{i,v_S} \not \subseteq \bigcup_{M \in \M_{i+1}}M$.
Hence $\bigcup_{M \in \M_i}V(M) \neq \bigcup_{M \in \M_{i+1}}V(M)$ by (D8gb).
By (D8ha), $V(M_{i,v_S}) \subseteq \bigcup_{M \in \M_i}V(M) - \bigcup_{M \in \M_{i+1}}V(M)$.
So by (D8hg), there exists $(S',j) \in \E_i$ with $\lvert S' \rvert = \lvert S \rvert$ such that 
	\begin{itemize}
		\item[(i)] $M_{i,v_{S'}} \subseteq M_{i+1,q}$, where $v_{S'}$ is the sink for $(S',j)$, 
		\item[(ii)] $S \cap V(G_{i+1}) - \{v_S\} = S' \cap V(G_{i+1}) - \{v_{S'}\}$, and 
		\item[(iii)] there exists a bijection $\iota$ from $S-(V(G_{i+1}) \cup \{v_S\})$ to $S' \cap V(M_{i+1,q}) - \{v_{S'}\}$ such that for every $x \in S-(V(G_{i+1}) \cup \{v_S\})$, $N_G(x) \cap U = N_G(\iota(x)) \cap U$.
	\end{itemize}
Since $(S-\{v_S\}) \cap V(M_{i+1,q}) = \emptyset$ and $V(M_{i,v_S}) \subseteq \bigcup_{M \in \M_i}V(M) - \bigcup_{M \in \M_{i+1}}V(M)$, we know $(\bigcup_{x \in S}V(M_{i,x})) \cap V(M_{i+1,q}) = \emptyset$.
Since $(S,j)$ is chosen such that $(\bigcup_{x \in S}V(M_{i,x})) \cap V(M_{i+1,q}) \neq \emptyset$ if possible, we know that $(S',j)$ satisfies the lemma by (i). 
That is, $\lvert c_i(S'-\{v_{S'}\}) \rvert \leq h-j-1$.

By (i), either 
	\begin{itemize}
		\item[(iv)] $v_{S'}=q \in V(G_i) \cap V(G_{i+1})$ and $S' \cap V(M_{i+1,q})-\{v_{S'}\} = \emptyset$, or 
		\item[(v)] $S' \cap V(M_{i+1,q}) \subseteq S'-V(G_{i+1})$.
	\end{itemize}

If (iv) holds and $S'-V(G_{i+1}) \neq \emptyset$, then there exists $u \in S'-\{v_{S'}\}$ with $M_{i,u}$ disjoint from all members of $\M_{i+1}$ by (D8ha), so (D8hb) implies that $q=v_{S'} \in N_{G_i}(u) \subseteq N_{G_{i+1}}(q) \cup (V(G_i)-V(G_{i+1}))$, a contradiction. 
So $S'-V(G_{i+1})=\emptyset = S' \cap V(M_{i+1,q})-\{v_{S'}\}$ when (iv) holds. 
Hence when (iv) holds, $\iota$ is a bijection from $S-(V(G_{i+1}) \cup \{v_S\})$ to $S'-(V(G_{i+1}) \cup \{v_{S'}\})$ by (iii).

Since $\lvert S' \rvert = \lvert S \rvert$, (ii) implies that $\lvert S-(V(G_{i+1}) \cup \{v_S\}) \rvert = \lvert S'-(V(G_{i+1}) \cup \{v_{S'}\}) \rvert$.
By (iii), $\lvert S-(V(G_{i+1}) \cup \{v_S\}) \rvert = \lvert S' \cap V(M_{i+1,q}) - \{v_{S'}\} \rvert$.
Hence $\lvert S'-(V(G_{i+1}) \cup \{v_{S'}\}) \rvert = \lvert S' \cap V(M_{i+1,q}) - \{v_{S'}\} \rvert$.
So if (v) holds, then $S' \cap V(M_{i+1,q})-\{v_{S'}\} = S'-(V(G_{i+1}) \cup \{v_{S'}\})$, so (iii) imply that $\iota$ is a bijection from $S-(V(G_{i+1}) \cup \{v_S\})$ to $S'-(V(G_{i+1}) \cup \{v_{S'}\})$.

Hence $\iota$ is a bijection from $S-(V(G_{i+1}) \cup \{v_S\})$ to $S'-(V(G_{i+1}) \cup \{v_{S'}\})$ in either case.
So by (iii), for every $x \in S-(V(G_{i+1}) \cup \{v_S\})$, $c_i(x) = \min([h-1]-c_{i+1}(N_G(x) \cap U))=\min([h-1]-c_{i+1}(N_G(\iota(x)) \cap U)) = c_i(\iota(x))$.
So $c_i(S-(V(G_{i+1}) \cup \{v_S\})) = c_i(S'-(V(G_{i+1}) \cup \{v_{S'}\}))$.
By (ii), $c_i(S \cap V(G_{i+1}) - \{v_S\}) = c_i(S' \cap V(G_{i+1}) - \{v_{S'}\})$.
Hence $c_i(S-\{v_S\}) = c_i(S'-\{v_{S'}\})$ has size at most $h-j-1$, a contradiction.
$\Box$

\medskip

\noindent{\bf Claim 7:} $M_{i,v_S} \subseteq M_{i+1,q}$, $\lvert V(M_{i+1,q}) \rvert \geq 2$, and $\{x \in S-\{v_S\}: \deg_{G_{i}}(x) >d\} \subseteq U^+$.

\noindent{\bf Proof of Claim 7:}
Since $M_{i,v_S} \subseteq M_{i+1,q}$ implies $\{x \in S-\{v_S\}: \deg_{G_{i}}(x) >d\} \subseteq U^+$ by (D8d), it suffices to prove that $M_{i,v_S} \subseteq M_{i+1,q}$ and $\lvert V(M_{i+1,q}) \rvert \geq 2$.

We first assume the second statement of Claim 6 holds.
Then $M_{i,v_S} \subseteq M_{i+1,q}$ and $\bigcup_{M \in \M_i}V(M)=\bigcup_{M \in \M_{i+1}}V(M)$.
Since $\lvert V(G_i) \rvert \neq \lvert V(G_{i+1}) \rvert$, (D8gb) implies that $\lvert V(M_{i+1,q}) \rvert \geq 2$.

Hence we may assume the first statement of Claim 6 holds.
That is, $(S-\{v_S\}) \cap V(M_{i+1,q}) \neq \emptyset$.

Suppose to the contrary that $M_{i,v_S} \subseteq M_{i+1,v}$ for some $v \in V(G_{i+1})-\{q\}$.
Since $(S-\{v_S\}) \cap V(M_{i+1,q}) \neq \emptyset$, there exists $(u',v_S) \in D_i$ with $M_{i,u'} \subseteq M_{i+1,q}$ by (D5). 
Since $M_{i,v_S} \subseteq M_{i+1,v}$, if $qv \in E(G_{i+1})$, then $(q,v) \in D_{i+1}$ by (D4).
Since (D8f) implies that $q$ is a sink for some member of $\E_{i+1}$, if $(q,v) \in D_{i+1}$, then $v$ is a head with respect to $D_{i+1}$, so $qv \not \in E(G_{i+1})$ by (D6b).
So $qv \not \in E(G_{i+1})$.
Since $u'v_S \in E(G_i)$ and $qv \not \in E(G_{i+1})$, we have $\bigcup_{M \in \M_{i+1}}V(M) = \bigcup_{M \in \M_i}V(M)$ by (D8hd).
Again since $u'v_S \in E(G_i)$ and $qv \not \in E(G_{i+1})$, (D8gc) implies that $v=v_S$ and $v_S$ is not the sink of a member of $\E_i$, a contradiction.

So $M_{i,v_S} \not \subseteq M_{i+1,v}$ for every $v \in V(G_{i+1})-\{q\}$.
Suppose to the contrary that $M_{i,v_S} \not \subseteq M_{i+1,q}$.
Then $\bigcup_{M \in \M_{i+1}}V(M) \neq \bigcup_{M \in \M_i}V(M)$ by (D8gb), and hence $M_{i,v_S}$ is disjoint from every member of $\M_{i+1}$ by (D8ha).
So (D8hb) implies that $N_{G_i}(v_S) \subseteq N_{G_{i+1}}(q) \cup \{x \in V(G_i): M_{i,x}$ is disjoint from all members of $\M_{i+1}\}$, contradicting $(S-\{v_S\}) \cap V(M_{i+1,q}) \neq \emptyset$.

Hence $M_{i,v_S} \subseteq M_{i+1,q}$.
Since $(S-\{v_S\}) \cap V(M_{i+1,q}) \neq \emptyset$, $V(M_{i+1,q})-V(M_{i,v_S}) \neq \emptyset$, so $\lvert V(M_{i+1,q}) \rvert \geq 2$.
$\Box$

\medskip

\noindent{\bf Claim 8:} For every $v \in S-(\{v_S\} \cup U^+)$ with $c_i(v) \neq 1$ and for every $\beta \in [c_i(v)-1]$, there exists $u_{v,\beta} \in N_{G}(v) \cap U$ such that $c_i(u_{v,\beta})=\beta$. 

\noindent{\bf Proof of Claim 8:}
Let $v \in S-(\{v_S\} \cup U^+)$ with $c_i(v) \neq 1$.
Let $\beta \in [c_i(v)-1]$.
Note that $v \in V(G) \cap V(G_i)$ by (D4) and (D5).
By Claim 4, there exists $u \in N_G(v) \cap V(G_i) \cap V(G_{i+1})$ with $\deg_{G_{i+1}}(u)>d$ and $c_i(u)=\beta$.
So to prove this claim, it suffices to show $u \in U$.
By Claim 4, we may assume $v \in V(G_i) \cap V(G_{i+1})$.

Since both $u$ and $v$ are in $V(G) \cap V(G_i) \cap V(G_{i+1})$ and $u \in N_G(v)$, $u \in N_{G_i}(v) \cap N_{G_{i+1}}(v)$ by (D2).
By Claim 7, $\{x \in S-\{v_S\}: \deg_{G_{i}}(x) >d\} \subseteq U^+$.
Since $v \in S-\{v_S\}$ but $v \not \in U^+$, $\deg_{G_i}(v) \leq d$.
Note that $v \in V(G) \cap V(G_{i+1}) \cap N_{G_i}(v_S)-U^+$.
Since $M_{i,v_S} \subseteq M_{i+1,q}$ by Claim 7, (D8e) implies that $u \in \{x \in N_{G_{i+1}}(v) \cap V(G): \deg_{G_{i+1}}(x)>d\} \subseteq U$.
$\Box$

\medskip

For every $\alpha \in c_i(S-\{v_S\})$, let $v_\alpha$ be a vertex in $S-\{v_S\}$ such that $c_i(v_\alpha)=\alpha$, and subject to this, $v_\alpha \in U^+$ if possible.

Let $A=\{\alpha \in c_i(S-\{v_S\}): v_\alpha \not \in U^+\}$.
Note that for every $\alpha \in A$, by the choice of $v_\alpha$, there exists no vertex $x$ in $(S-\{v_S\}) \cap U^+$ such that $c_i(x)=\alpha$, so $\alpha \not \in c_i((S-\{v_S\}) \cap U^+)$.

A {\it segment} is a maximal set of consecutive integers in $A$.
We say that a segment $I$ is {\it tilted} if either $\min I=1$, or $\min I \geq 2$ and $\min I-1 \in c_i(S-\{v_S\})$ and $v_{\min I-1} \in U^+$.
For any two distinct tilted segments $I_1,I_2$, we define $I_1 \prec I_2$ if and only if $\min I_1 < \min I_2$.
So $\prec$ is a total order of all tilted segments.
A tilted segment is {\it $\prec$-minimum} if it is the minimum respect to $\prec$.
For a non-$\prec$-minimum tilted segment $I$, we define the {\it precursor} of $I$ to be the tilted segment $I'$ such that $I' \prec I$ and there exists no tilted segment $I''$ with $I' \prec I'' \prec I$; and we define $g(I)$ to be $\max I'$, where $I'$ is the precursor of $I$.

Define $g$ to be the function with domain $c_i(S-\{v_S\})$ such that for every $x \in c_i(S-\{v_S\})$, 
	\begin{itemize}
		\item if $x \in c_i(S-\{v_S\})-A$, then $g(x)=x$, 
		\item if $x \in A$ and $x$ is not the minimum of any tilted segment, then $g(x)=x-1$, 
		\item if $x \in A$ and $x$ is the minimum of some non-$\prec$-minimum tilted segment $I$, then $g(x)=g(I)$, 
		\item if $x \in A$ and $x$ is the minimum of the $\prec$-minimum tilted segment, then $g(x)=0$. 
	\end{itemize}
Note that $g$ maps each element in $A$ to $0$ or an element in $A \cup ({\mathbb N}-c_i(S-\{v_S\}))$.

\medskip

\noindent{\bf Claim 9:} $\lvert g(c_i(S-\{v_S\})) \rvert \geq h-j$ and $\lvert g(c_i(S-\{v_S\}))-\{0\} \rvert \geq h-j-1$. 

\noindent{\bf Proof of Claim 9:}
We first show that $g$ is injective.
Suppose to the contrary that there exist distinct $x,y \in c_i(S-\{v_S\})$ such that $g(x)=g(y)$.
If one of $x,y$ is not in $A$, say $x \not \in A$, then $g(x)=x \in c_i(S-\{v_S\})-A$, and since $g$ maps each element in $A$ to $0$ or an element in $A \cup ({\mathbb N}-c_i(S-\{v_S\}))$, $g(x)=g(y)$ implies that $y \in c_i(S-\{v_S\})-A$, so $x=g(x)=g(y)=y$, a contradiction.
So $x,y \in A$.
Note that the minimum of the $\prec$-minimum tilted segment is the only element mapped to 0 by $g$.
So none of $x,y$ is the minimum of the $\prec$-minimum tilted segment.
If each of $x,y$ is the minimum of some non-$\prec$-minimum tilted segment, then $g(x)=g(y)$ implies that $x=y$, a contradiction.
Hence by symmetry, we may assume that $x \in A$ is not the minimum of any tilted segment and $y$ is the minimum of some non-$\prec$-minimum tilted segment $I$.
So $g(y)$ is the maximum element of some tilted segment, but $g(x)$ is not, a contradiction.

Hence $g$ is injective.
So $\lvert g(c_i(S-\{v_S\})) \rvert \geq \lvert c_i(S-\{v_S\}) \rvert \geq h-j$.
Hence $\lvert g(c_i(S-\{v_S\}))-\{0\} \rvert \geq h-j-1$.
$\Box$

\medskip

If $(S-\{v_S\}) \cap V(M_{i+1,q}) \neq \emptyset$, then define $S^*=S$ and $v_{S^*}=v_S$, and define $\iota$ to be the identity function from $S^*$ to $S^*$; otherwise, $\bigcup_{M \in \M_i}V(M)=\bigcup_{M \in \M_{i+1}}V(M)$ and $M_{i,v_S} \subseteq M_{i+1,q}$ by Claim 6, so (D8gd) implies that there exists $(S^*,j) \in \E_i$ with $\bigcup_{v \in S^*-U^+}M_{i,v} \subseteq M_{i+1,q}$, $S^* \cap U^+ = S \cap U^+$ and $v_{S^*} \not \in U^+$, where $v_{S^*}$ is the sink for $(S^*,j)$, and there exists a bijection $\iota: S-(\{v_S\} \cup U^+) \rightarrow S^*-(\{v_{S^*}\} \cup U^+)$ such that for every $v \in S-(\{v_S\} \cup U^+)$, $N_{G_i}(v) \cap U^+ =  N_{G_i}(\iota(v)) \cap U^+$.

\medskip

\noindent{\bf Claim 10:} $M_{i,v} \cup M_{i,v_{S^*}} \subseteq M_{i+1,q}$ for some $v \in S^*-\{v_{S^*}\}$.

\noindent{\bf Proof of Claim 10:}
If $(S-\{v_S\}) \cap V(M_{i+1,q}) \neq \emptyset$, then Claim 7 implies that $M_{i,v} \cup M_{i,v_{S^*}} \subseteq M_{i+1,q}$ for some $v \in S^*-\{v_{S^*}\}$.
Otherwise, Claim 6 implies $M_{i,\iota(v)} \cup M_{i,v_{S^*}} \subseteq M_{i+1,q}$ for some $v \in S-(\{v_S\} \cup U^+)$. 
$\Box$

\medskip

Note that for every $\alpha \in c_i(S-\{v_S\})$, $\alpha \not \in A$ if and only if $v_\alpha \in U^+$.
And if $\alpha=1$ and $v_\alpha \not \in U^+$, then $\alpha = 1 \in A$ is the minimum of the $\prec$-minimum tilted segment, so $g(\alpha)=0$. 
So if $v_\alpha \not \in U^+$ and $g(\alpha) \neq 0$, then $c_i(v_\alpha) \neq 1$.

Define $f$ to be the function with domain $\{v_\alpha: \alpha \in c_i(S-\{v_S\}), g(\alpha) \neq 0\}$ such that for every $\alpha \in c_i(S-\{v_S\})$ with $g(\alpha) \neq 0$,
	\begin{itemize}
		\item if $v_\alpha \in U^+$, then we define $f(v_\alpha)=v_\alpha$,
		\item if $v_\alpha \not \in U^+$, then since $g(\alpha) \neq 0$, we know $c_i(v_\alpha) \neq 1$, so there exists a vertex $u_{v_\alpha,g(\alpha)}$ in $N_G(v_\alpha) \cap U$ with $c_i(u_{v_\alpha,g(\alpha)})=g(\alpha)$ by Claim 8, and we define $f(v_\alpha)=u_{v_\alpha,g(\alpha)}$. 
	\end{itemize}

\noindent{\bf Claim 11:} For every $\alpha \in c_i(S-\{v_S\})$ with $g(\alpha) \neq 0$, if $v_\alpha \not \in U^+$, then $f(v_\alpha)=u_{v_\alpha,g(\alpha)} \in N_G(v_\alpha) \cap U = N_{G_i}(v_\alpha) \cap U$ and $c_i(f(v_\alpha))=c_{i+1}(f(v_\alpha))=g(\alpha)$.

\noindent{\bf Proof of Claim 11:}
Note that $f(v_\alpha)=u_{v_\alpha,g(\alpha)} \in N_G(v_\alpha) \cap U = N_{G_i}(v_\alpha) \cap U$ by (D2) since $\{v_\alpha\} \cup U \subseteq V(G) \cap V(G_i)$.
And since $U \subseteq V(G_i) \cup V(G_{i+1})$ by (D8), $c_i(f(v_\alpha))=c_{i+1}(f(v_\alpha))=g(\alpha)$.
$\Box$

\medskip

\noindent{\bf Claim 12:} There exists a tilted segment.

\noindent{\bf Proof of Claim 12:}
Let $\C$ be the multiset $\{N_G(x) \cap U \neq \emptyset: x \in S^*-(\{v_{S^*}\} \cup U^+)\}$.
Note that $\C$ equals the multiset $\{N_G(\iota(x)) \cap U \neq \emptyset: x \in S-(\{v_{S}\} \cup U^+)\}$.
For each $Z \in \C$, let $x_Z$ be the element in $S-(\{v_{S}\} \cup U^+)$ such that $Z=N_G(\iota(x_Z)) \cap U$. 
For every $\alpha \in A-\{1\}$, let $Z_\alpha$ be the element in the multiset $\C$ such that $x_{Z_\alpha}=v_\alpha$.

Suppose to the contrary that there exists no tilted segment.
So $1 \not \in A$ and $g(\alpha) \neq 0$ for every $\alpha \in c_i(S-\{v_S\})$.
Let $f'$ be a function such that for each member $Z$ of the multiset $\C$,
	\begin{itemize}
		\item if $x_Z \in \{v_\alpha: \alpha \in A\}$, then $f'(Z)=f(x_Z)$, 
		\item otherwise, $f'(Z)$ is an arbitrary element in $Z$.
	\end{itemize}
So $f'$ maps each member $Z$ of $\C$ to a vertex in $Z$.

Recall that $(S^*,j)$ is a member of $\E_i$ such that $M_{i,v} \cup M_{i,v_{S^*}} \subseteq M_{i+1,q}$ for some $v \in S^*-\{v_{S^*}\}$ by Claim 10.
So by (D8ia), $(S',j) \in \E_{i+1}$, where $S'=(S^* \cap U^+) \cup \{f'(Z): Z \in \C\} \cup \{q\}$, and $q$ is the sink for $(S',j)$.
Since $c_{i+1}$ satisfies this lemma, $\lvert c_{i+1}(S'-\{q\}) \rvert \leq h-j-1$.

Note that for every $\alpha \in c_i(S-\{v_S\})$, if $\alpha \not \in A$, then $v_\alpha \in S \cap U^+ = S^* \cap U^+ \subseteq S'-\{q\}$ and $g(\alpha)=\alpha=c_i(v_\alpha)=c_{i+1}(v_\alpha) \in c_{i+1}(S'-\{q\})$; if $\alpha \in A$, then $v_\alpha \not \in U^+$, and since $1 \not \in A$, we know $v_\alpha=x_{Z_\alpha}$, so $f(v_\alpha)=f'(Z_\alpha) \in \{f'(Z): Z \in \C\} \subseteq S'-\{q\}$, and hence by Claim 11, $g(\alpha)=c_{i+1}(f(v_\alpha)) \in c_{i+1}(S'-\{q\})$. 
Therefore, $c_{i+1}(S'-\{q\}) \supseteq g(c_i(S-\{v_S\}))$.
So $h-j-1 \geq \lvert c_{i+1}(S'-\{q\}) \rvert \geq \lvert g(c_i(S-\{v_S\})) \rvert$.
But by Claim 9, $h-j-1 \geq \lvert g(c_i(S-\{v_S\})) \rvert \geq h-j$, a contradiction.
$\Box$

\medskip

By Claim 12, there exists a tilted segment.
So there exists a unique element $\alpha^* \in c_i(S-\{v_S\})$ such that $g(\alpha^*)=0$.
Note that $\alpha^* \in A$.
So $v_{\alpha^*} \not \in U^+$ and hence $\iota(v_{\alpha^*}) \in S^*-(\{v_{S^*}\} \cup U^+)$.

Let $\C^*$ be the multiset $\{N_G(x) \cap U \neq \emptyset: x \in S^*-(\{\iota(v_{\alpha^*}),v_{S^*}\} \cup U^+)\}$.
Note that $\C^*$ equals the multiset $\{N_G(\iota(x)) \cap U \neq \emptyset: x \in S-(\{v_{\alpha^*},v_{S}\} \cup U^+)\}$ by the definition of $\iota$.

\medskip

\noindent{\bf Claim 13:} $S \cap U^+ = S^* \cap U^+ \neq \emptyset$ or $\C^* \neq \emptyset$.

\noindent{\bf Proof of Claim 13:}
Recall $S^* \cap U^+=S \cap U^+$ by the definition of $S^*$.
Suppose to the contrary that $S \cap U^+=\emptyset=\C^*$.
If there exists $y \in S-(\{v_S,v_{\alpha^*}\} \cup U^+)$ with $c_i(y) \neq 1$, then by Claim 8, $N_G(y) \cap U \neq \emptyset$, so $\iota(y) \in S^*-(\{\iota(v_{\alpha^*}),v_{S^*}\} \cup U^+)$ and $N_G(\iota(y)) \cap U = N_G(y) \cap U \neq \emptyset$ (by (D2) and the property of $\iota$), and hence $\C^* \neq \emptyset$, a contradiction.
Therefore, $c_i(S-(\{v_S,v_{\alpha^*}\} \cup U^+)) \subseteq \{1\}$.
Since $S \cap U^+=\emptyset$, $c_i(S-\{v_S\})=c_i(S-(\{v_S\} \cup U^+))$.
If $c_i(S-(\{v_S,v_{\alpha^*}\} \cup U^+)) = \emptyset$, then $c_i(S-\{v_S\})=c_i(S-(\{v_S\} \cup U^+)) = \{c_i(v_{\alpha^*})\}$ has size at most 1, contradicting $\lvert c_i(S-\{v_S\}) \rvert \geq h-j \geq 2$ by (D5).
So $c_i(S-(\{v_S,v_{\alpha^*}\} \cup U^+)) = \{1\}$.
Since $S \cap U^+=\emptyset$, $1 \in A$.
So $\alpha^*=1$.
Hence $c_i(S-\{v_S\})=c_i(S-(\{v_S\} \cup U^+)) = \{1,c_i(v_{\alpha^*})\} = \{1\}$ has size smaller than $h-j$ by (D5), a contradiction.
$\Box$

\medskip

For each member $Z$ of the multiset $\C^*$, let $x_Z$ be an element in $S-(\{v_{\alpha^*},v_S\} \cup U^+)$ such that $Z=N_G(\iota(x_Z)) \cap U$. 
For every $\alpha \in A-\{\alpha^*\}$, let $Z_\alpha$ be the member of the multiset $\C^*$ such that $x_{Z_\alpha}=v_\alpha$. 

Let $f'$ be a function such that for each member $Z$ of the multiset $\C^*$, 
	\begin{itemize}
		\item if $x_Z \in \{v_\alpha: \alpha \in A-\{\alpha^*\}\}$, then $f'(Z)=f(x_Z)$, (note that $f(x_Z)$ is defined since $\alpha^*$ is the unique element such that $g(\alpha^*)=0$)
		\item otherwise, $f'(Z)$ is an arbitrary element in $Z$.
	\end{itemize}
So $f'$ maps each member $Z$ of $\C^*$ to a vertex in $Z$.
Recall that $(S^*,j)$ is a member of $\E_i$ such that $M_{i,v} \cup M_{i,v_{S^*}} \subseteq M_{i+1,q}$ for some $v \in S^*-\{v_{S^*}\}$ by Claim 10.
Hence by Claim 13 and (D8ib), $(S',j+1) \in \E_{i+1}$, where $S'=(S^* \cap U^+) \cup \{f'(Z): Z \in \C^*\} \cup \{q\} =(S \cap U^+) \cup \{f'(Z): Z \in \C^*\} \cup \{q\}$, and $q$ is the sink for $(S',j)$.
Since $c_{i+1}$ satisfies this lemma, $\lvert c_{i+1}(S'-\{q\}) \rvert \leq h-j-2$.

Note that for every $\alpha \in c_i(S-\{v_S\})$, if $\alpha \not \in A$, then $v_\alpha \in S \cap U^+ \subseteq S'-\{q\}$ and $g(\alpha)=\alpha=c_i(v_\alpha)=c_{i+1}(v_\alpha) \in c_{i+1}(S'-\{q\})$; if $\alpha \in A-\{\alpha^*\}$, then $v_\alpha \not \in U^+$ and $g(\alpha) \neq 0$ and $v_\alpha=x_{Z_\alpha}$, so $f(v_\alpha)=f'(Z_\alpha) \in \{f'(Z): Z \in \C^*\} \subseteq S'-\{q\}$, and hence by Claim 11, $g(\alpha)=c_{i+1}(f(v_\alpha)) \in c_{i+1}(S'-\{q\})$. 
Therefore, $c_{i+1}(S'-\{q\}) \supseteq g(c_i(S-\{v_S\})-\{\alpha^*\}) = g(c_i(S-\{v_S\}))-\{0\}$.
So $h-j-2 \geq \lvert c_{i+1}(S'-\{q\}) \rvert \geq \lvert g(c_i(S-\{v_S\}))-\{0\} \rvert$, contradicting Claim 9.

This proves the lemma.
\end{pf}

\section{Homogeneous structures} \label{sec:homo}

A {\it geodesic} in a graph $G$ is a path $P$ in $G$ such that its length equals the distance between its ends in $G$. 
Note that every subpath of a geodesic is a geodesic.

\begin{lemma} \label{geodesic}
For any positive integers $t,k,\ell$, there exists a positive integer $n=n(t,k,\ell)$ with $n > t\ell$ such that for every graph $G$ having a vertex $v^*$ with $V(G)=N_G^{\leq n}[v^*]$, every function $f: V(G) \rightarrow [t]$ and every geodesic $P$ in $G$ with an end $v^*$ on $n-t\ell$ vertices, there exist a subpath $Q$ of $P$ and a nonempty set $S \subseteq [t]$ such that 
	\begin{enumerate}
		\item for every $x \in [t]-S$, there exists no vertex $u \in N_G^{\leq \lvert S \rvert\ell}[V(Q)]$ with $f(u)=x$, and
		\item $Q$ can be partitioned into $k$ disjoint subpaths $Q_1,Q_2,...,Q_k$ of $Q$ with the same length such that for any $i \in [k]$ and $x \in S$, $N_G^{\leq (\lvert S \rvert-1)\ell}[V(Q_i)]$ contains a vertex $v_{i,x}$ with $f(v_{i,x})=x$. 
	\end{enumerate}
\end{lemma}

\begin{pf}
For any positive integers $y,z$, define $n(1,y,z)=y+z$, and for every positive integer $x \geq 2$, $n(x,y,z)=y \cdot (n(x-1,y,z)-(x-1)z) + xz$. 
Clearly, $n(x,y,z) > xz$ for any positive integers $x,y,z$.

Let $t,k,\ell$ be positive integers.
Let $n=n(t,k,\ell)$.
We shall prove this lemma by induction on $t$.

Let $G,v^*,f,P$ be as stated in the lemma.
When $t=1$, $P$ has $n(1,k,\ell)-\ell \geq k$ vertices, so every subpath $Q$ of $P$ on $k$ vertices and the set $[1]$ satisfy the conclusion of this lemma (with $S=[1]$).
So we may assume that $t \geq 2$ and the lemma holds when $t$ is smaller.

Let $W$ be a subpath of $P$ on $n(t-1,k,\ell)-(t-1)\ell$ vertices.
Since $W$ is a subpath of a geodesic $P$, $W$ is a geodesic.

We first assume that there exists $x_W \in [t]$ such that there exists no vertex $u \in N_G^{\leq (t-1)\ell}[V(W)]$ with $f(u)=x_W$.
By symmetry, we may assume $x_W=t$.
Let $H=G[N_G^{\leq (t-1)\ell}[V(W)]]$.
Let $v_W$ be the end of $W$ closest to $v^*$ in $G$.
Hence $V(W) \subseteq N_H^{\leq n(t-1,k,\ell)-(t-1)\ell-1}[v_W]$.
So $V(H)=N_H^{\leq n(t-1,k,\ell)-1}[v_W]$.
Hence $V(H)=N_H^{\leq n(t-1,k,\ell)}[v_W]$.
Moreover, $f|_{V(H)}$ is a function from $V(H)$ to $[t-1]$, and $W$ is a geodesic in $H$ with an end $v_W$ on $n(t-1,k,\ell)-(t-1)\ell$ vertices.
Hence by the induction hypothesis, there exist a subpath $W'$ of $W$ and a nonempty set $S_W \subseteq [t-1]$ such that 
	\begin{itemize}
		\item[(i)] for every $x \in [t-1]-S_W$, there exists no $u \in N_H^{\leq \lvert S_W \rvert\ell}[V(W')]$ with $f(u)=x$, and
		\item[(ii)] $W'$ can be partitioned into $k$ disjoint subpaths $W_1,W_2,...,W_k$ of $W'$ with the same length such that for every $i \in [k]$ and $x \in S_W$, $N_H^{\leq (\lvert S_W \rvert-1)\ell}[V(W_i)]$ contains a vertex $w_{i,x}$ with $f(w_{i,x})=x$. 
	\end{itemize}
Since $S_W \subseteq [t-1]$, $\lvert S_W \rvert \leq t-1$. 
So $G[N_G^{\leq \lvert S_W \rvert\ell}[V(W')]]$ is a subgraph of $H$.
Hence for every $i \in [\lvert S_W \rvert]$, $N_G^{\leq i\ell}[V(W')] = N_H^{\leq i\ell}[V(W')]$.
Recall that there exists no vertex $u \in N_G^{\leq (t-1)\ell}[V(W)]$ with $f(u)=x_W=t$.
So (i) implies that for every $x \in [t]-S_W$, there exists no $u \in N_G^{\leq \lvert S_W \rvert\ell}[V(W')]$ with $f(u)=x$.
Hence this lemma follows from taking $Q=W'$ and $S=S_W$.

So we may assume that for every subpath $P'$ of $P$ on $n(t-1,k,\ell)-(t-1)\ell$ vertices and for every $x \in [t]$, $N_G^{\leq (t-1)\ell}[V(P')]$ contains a vertex $v_{P',x}$ with $f(v_{P',x})=x$.
Since $P$ has $n(t,k,\ell)-t\ell \geq k \cdot (n(t-1,k,\ell)-(t-1)\ell)$ vertices, there exists a subpath $P^*$ of $P$ on $k \cdot (n(t-1,k,\ell)-(t-1)\ell)$ vertices.
So $P^*$ can be partitioned into $k$ disjoint subpaths $P_1,P_2,...,P_k$ of $P^*$ each having $n(t-1,k,\ell)-(t-1)\ell$ vertices.
So for any $i \in [k]$ and $x \in [t]$, $N_G^{\leq (t-1)\ell}[V(P_i)]$ contains a vertex $v_{i,x}$ with $f(v_{i,x})=x$.
Hence we are done by choosing $Q=P^*$ and $S=[t]$.
This proves the lemma.
\end{pf}

\bigskip

Note that for every proper minor-closed family $\G$, there exist positive integers $r$ and $k$ such that every graph in $\G$ has edge-density at most $k$ \cite{m} and no $K_{r,r}$-minor.
By taking $\epsilon=\epsilon'=0$ in \cite[Lemma 4.4]{l_homo}, we obtain the following immediate corollary for proper minor-closed families.

\begin{lemma}[{special case of {\cite[Lemma 4.4]{l_homo}}}] \label{homo} 
For any proper minor-closed family $\G$, there exists a positive integer $r=r(\G)$ such that for any integers $t \geq 1,\ell \geq 2$, there exist positive integers $d=d(\G,t,\ell),N=N(\G,t,\ell)$ such that for any graph $G \in \G$ with $\lvert V(G) \rvert >N$, there exist $X,Z,W \subseteq V(G)$ with $Z \subseteq X$, $\lvert Z \rvert =t$, $W \subseteq V(G)-X$ and $\lvert W \rvert \leq r-1$ such that
	\begin{enumerate}
		\item every vertex in $X$ has degree at most $d$ in $G$,
		\item for any distinct $z,z' \in Z$, the distance in $G[X]$ between $z,z'$ is at least $2\ell-1$, and
		\item $N_G(N_{G[X]}^{\leq \ell-1}[z])-X =W$ for every $z \in Z$.
	\end{enumerate}
\end{lemma}

Note that Lemma \ref{homo} can also be derived from the machinery developed in \cite{lw_phase}, which is a simpler version of the machinery developed in \cite{l_homo}.

\section{Strong elimination schemes} \label{sec:strong}

Let $G$ be a graph.
Let $h \geq 3,k,r,d,N$ be positive integers.
Then a {\it strong $(G,h,k,r,d,N)$-defective elimination scheme} is a sequence $((G_i,\M_i,\E_i,D_i,\A_i,\A_i'): i \in {\mathbb N})$ of tuples such that $((G_i,\M_i,\E_i,D_i): i \in {\mathbb N})$ is a $(G,h,k,r,d,N)$-defective elimination scheme such that $\A_1=\A'_1=\emptyset$, and for every $i \geq 2$, the following hold: 
	\begin{enumerate}
		\item[] (Recall that we denote the member of $\M_i$ corresponding to a vertex $v$ of $G_i$ by $M_{i,v}$.)
		\item[(D9)] $\A_i$ is a collection $\{\A_{i,Q}: Q \in \E_i\}$, and $\A_i'$ is a collection $\{\A'_{i,Q}: Q \in \E_i\}$.
		\item[(D10)] For every $Q=(S,j) \in \E_i$, $\A_{i,Q}$ is a set of pairwise disjoint connected subgraphs of $G[V(M_{i,v_S}) \cup (V(G)-\bigcup_{v \in V(G_i)}V(M_{i,v}))]$ (where $v_S$ is the sink for $Q$), and $\A'_{i,Q}$ is a set of $\lvert S \rvert-1$ pairwise disjoint connected subgraphs of $G[\bigcup_{v \in S}V(M_{i,v}) \cup (V(G)-\bigcup_{v \in V(G_i)}V(M_{i,v}))]$ such that the following hold.  
			\begin{itemize}
				\item For each member $A$ of $\A'_{i,Q}$, 
					\begin{itemize}
						\item $V(A) \cap S-\{v_S\} \neq \emptyset$, 
						\item $A$ is disjoint from all members of $\A_{i,Q}$, and 
						\item $A$ is adjacent in $G$ to all members of $\A_{i,Q}$. 
					\end{itemize}
					(Note that $\lvert V(A) \cap S-\{v_S\} \rvert=1$ for each $A \in \A'_{i,Q}$ since $\lvert \A'_{i,Q} \rvert = \lvert S \rvert-1 = \lvert S -\{v_S\} \rvert$.)
				\item Contracting each member of $\A_{i,Q}$ into a single vertex creates a $(k+h-j)\CT_{j,k}$-minor in $G$. 
			\end{itemize}
		\item[(D11)] For any $Q_1=(S_1,j_1),Q_2=(S_2,j_2) \in \E_i$ with distinct sinks, 
			\begin{itemize}
						\item every member of $\A_{i,Q_1}$ is disjoint from every member of $\A_{i,Q_2}$, 
						\item for $\alpha \in [2]$, every member of $\A_{i,Q_\alpha}$ is disjoint from every member of $\A'_{i,Q_{3-\alpha}}$, and 
						\item if $A_1'$ is a member of $\A'_{i,Q_1}$ and $A'_2$ is a member of $\A'_{i,Q_2}$, then $V(A_1') \cap V(A_2') \subseteq V(A_1') \cap V(S_1) \cap V(A_2') \cap V(S_2)$. 
			\end{itemize}
		\item[(D12)] For every $v \in V(G_i)$, if $\lvert V(M_{i,v}) \rvert \geq 2$ or $v$ is a head (with respect to $D_i$) or the sink for some member of $\E_i$, then $\deg_{G_i}(v) \leq r$.
	\end{enumerate}
For any positive integer $n$, a {\it strong $n$-$(G,h,k,r,d,N)$-defective elimination scheme} is a sequence $((G_i,\M_i,\E_i,D_i,\A_i,\A_i'): i \in [n])$ of tuples such that $((G_i,\M_i,\E_i,D_i): i \in [n])$ is an $n$-$(G,h,k,r,d, \allowbreak N)$-defective elimination scheme and $(G_i,\M_i,\E_i,D_i,\A_i,\A_i')$ satisfies (D9)-(D12) for every $i \in [n]$, and $\A_1=\A'_1=\emptyset$.

\begin{lemma} \label{del}
For any positive integers $h \geq 3,k,r$, there exists a positive integer $t=t(h,k,r)$ such that for any positive integers $d \geq 2,\ell_0$, there exists a positive integer $N^*=N^*(h,k,r,d,\ell_0)$ such that the following hold.
Let $G$ be a graph with no $\CT_{h,k}$-minor.
Let $N$ be an integer with $N \geq N^*$.
Let $i$ be a positive integer.
Let $((G_\alpha,\M_\alpha,\E_\alpha,D_\alpha,\A_\alpha,\A_\alpha'): \alpha \in [i])$ be a strong $i$-$(G,h,k,r,d,N)$-defective elimination scheme.
If $\lvert V(G_i) \rvert>N$, and there exist $X,Z,W \subseteq V(G_i)$ with $Z \subseteq X$, $\lvert Z \rvert =t$, $W \subseteq V(G_i)-X$ and $\lvert W \rvert \leq r-1$ such that
	\begin{enumerate}
		\item every vertex in $X$ has degree at most $d$ in $G_i$,
		\item for any distinct $z,z' \in Z$, the distance in $G_i[X]$ between $z,z'$ is at least $2\ell_0-1$, and
		\item $N_{G_i}(N_{G_i[X]}^{\leq \ell_0-1}[z]) =W$ for every $z \in Z$, 
	\end{enumerate}
then there exists a graph $G_{i+1}$ with $\lvert V(G_{i+1}) \rvert < \lvert V(G_i) \rvert$ and a tuple $(G_{i+1},\M_{i+1},\E_{i+1},D_{i+1}, \allowbreak \A_{i+1}, \allowbreak \A_{i+1}')$ such that $((G_\alpha,\M_\alpha,\E_\alpha,D_\alpha,\A_\alpha,\A'_\alpha): \alpha \in [i+1])$ is a strong $(i+1)$-$(G,h,k,r,d,N)$-defective elimination scheme.
\end{lemma}

\begin{pf}
Let $h \geq 3,k,r$ be positive integers.
Define $t=2^{(h-2)(r+1)^{2^{r-1}}(k+h)2^{r-1}} \cdot 2^{2^{r-1}}$.
Let $d \geq 2,\ell_0$ be positive integers.
Define $N^* = (k+h)d^{\ell_0}$.

Let $G$ be a graph with no $\CT_{h,k}$-minor.
Let $N$ be an integer with $N \geq N^*$.
Let $i$ be a positive integer.
Let $((G_\alpha,\M_\alpha,\E_\alpha,D_\alpha,\A_\alpha,\A_\alpha'): \alpha \in [i])$ be a strong $i$-$(G,h,k,r,d,N)$-defective elimination scheme. 
Assume $\lvert V(G_i) \rvert > N$.

By assumption, there exist $X,Z,W \subseteq V(G_i)$ with $Z \subseteq X$, $\lvert Z \rvert =t$, $W \subseteq V(G_i)-X$ and $\lvert W \rvert \leq r-1$ such that
	\begin{itemize}
		\item[(i)] every vertex in $X$ has degree at most $d$ in $G_i$,
		\item[(ii)] for any distinct $z,z' \in Z$, the distance in $G_i[X]$ between $z,z'$ is at least $2\ell_0-1$, and
		\item[(iii)] $N_{G_i}(N_{G_i[X]}^{\leq \ell_0-1}[z]) =W$ for every $z \in Z$. 
	\end{itemize}

\noindent{\bf Claim 1:} For every $v \in W$, $v \in V(G) \cap V(G_i)$, $v$ is not a head (with respect to $D_i$), and $v$ is not the sink for a member of $\E_i$.

\noindent{\bf Proof of Claim 1:}
Suppose to the contrary that there exists a vertex $v \in W$ such that $\lvert V(M_{i,v}) \rvert \geq 2$, or $v$ is a head (with respect to $D_i$), or $v$ is the sink for some member of $\E_i$.
By (D12), $\deg_{G_i}(v) \leq r$.
By (iii), for every $z \in Z$, $v$ is adjacent in $G_i$ to a vertex in $N_{G_i[X]}^{\leq \ell_0-1}[z]$.
So (ii) implies that $\deg_{G_i}(v) \geq \lvert Z \rvert$.
Hence $t=\lvert Z \rvert \leq \deg_{G_i}(v) \leq r$, a contradiction.
$\Box$

\medskip

Let $\sigma$ be a linear ordering of the subsets of $W$.
For any $z \in Z$, define the following: 
	\begin{itemize}
		\item For any subset $T$ of $W$ and nonnegative integers $j_0,j_1,...,j_{2^{\lvert W \rvert}}$ with $j_0 \in [h-2]$ and $0 \leq j_\alpha \leq r$ for every $\alpha \in [2^{\lvert W \rvert}]$, 
			\begin{itemize}
				\item $a_{z,T,j_0,...,j_{2^{\lvert W \rvert}}}=1$ if there exists $(S,j_0) \in \E_i$ such that 
					\begin{itemize}
						\item the sink $v_S$ for $(S,j_0)$ is in $N_{G_i[X]}^{\leq \ell_0-1}[z]$, 
						\item $S \cap W=T$, and
						\item for every $\alpha \in [2^{\lvert W \rvert}]$, there are exactly $j_\alpha$ vertices $x$ in $S-\{v_S\}$ such that $N_{G}(x) \cap W$ equals the $\alpha$-th subset of $W$ based on $\sigma$; 
					\end{itemize}
				\item otherwise, $a_{z,T,j_0,...,j_{2^{\lvert W \rvert}}}=0$. 
			\end{itemize}
		\item For any $T \subseteq W$,
			\begin{itemize}
				\item $a_{z,T}=1$ if there exists $v \in N_{G_i[X]}^{\leq \ell_0-1}[z] \cap V(G)$ such that $N_{G}(v) \cap W = T$;
				\item $a_{z,T}=0$ otherwise.
			\end{itemize}
	\end{itemize}
Since $\lvert W \rvert \leq r-1$, there exists $Z_1 \subseteq Z$ with $$\lvert Z_1 \rvert \geq \frac{\lvert Z \rvert}{2^{2^{\lvert W \rvert}(h-2)(r+1)^{2^{\lvert W \rvert}}}} \geq \frac{t}{2^{2^{r-1}(h-2)(r+1)^{2^{r-1}}}} \geq (k+h)2^{2^{r-1}}$$ such that $a_{z_1,T,j_0,j_1,...,j_{2^{\lvert W \rvert}}}=a_{z_2,T,j_0,j_1,...,j_{2^{\lvert W \rvert}}}$ for any $z_1,z_2 \in Z_1$, $T \subseteq W$ and nonnegative integers $j_0 \in [h-2], j_1,...,j_{2^{\lvert W \rvert}}$ with $0 \leq j_\alpha \leq r$ for every $\alpha \in [2^{\lvert W \rvert}]$.
Similarly, there exists $Z_2 \subseteq Z_1$ with $\lvert Z_2 \rvert \geq \frac{\lvert Z_1 \rvert}{2^{2^{r-1}}} \geq \frac{(k+h)2^{2^{r-1}}}{2^{2^{r-1}}} \geq k+h$ such that $a_{z_1,T}=a_{z_2,T}$ for any $z_1,z_2 \in Z_2$ and $T \subseteq W$.
Let $Z^*$ be a subset of $Z_2$ with $\lvert Z^* \rvert = k+h$.

Let $z^*$ be a vertex in $Z^*$.
Define $G_{i+1}$ to be the graph obtained from $G_i-\bigcup_{z \in Z^*-\{z^*\}}N_{G_i[X]}^{\leq \ell_0-1}[z]$ by contracting $N_{G_i[X]}^{\leq \ell_0-1}[z^*]$ into a new vertex $v^*$, and deleting resulting parallel edges and loops.
Note that $\lvert V(G_{i+1}) \rvert < \lvert V(G_i) \rvert$ since $\lvert Z^* \rvert \geq k+h \geq 3$.
Define $\M_{i+1} = (\M_i-\{M_{i,v}: v \in \bigcup_{z \in Z^*}N_{G_i[X]}^{\leq \ell_0-1}[z]\}) \cup \{G[\bigcup_{v \in N_{G_i[X]}^{\leq \ell_0-1}[z^*]}V(M_{i,v})]\}$.
Then $(G_{i+1},\M_{i+1})$ satisfies (D1)-(D3) by (ii).
Note that $M_{i+1,v^*}=G[\bigcup_{v \in N_{G_i[X]}^{\leq \ell_0-1}[z^*]}V(M_{i,v})]$.

Define $D_{i+1} = \{(u,v) \in D_i: u,v \in V(G_i) \cap V(G_{i+1})\} \cup \{(w,v^*): w \in W\}$.
By (iii) and Claim 1, $(G_{i+1},\M_{i+1},D_{i+1})$ satisfies (D4).

Let $U^+_{i+1}=W$ and $q_{i+1}=v^*$.
Let $U_{i+1} = \{v \in U_{i+1}^+: \deg_{G_{i+1}}(v) >d\}$.

Define the following:
	\begin{itemize}
		\item $\E_{i+1,0}=\{(S,j) \in \E_i: S \subseteq V(G_i) \cap V(G_{i+1})\}$.
		\item $\E_{i+1,1}=\{(\{v^*\} \cup (S \cap W) \cup T,j):$ there exists $(S,j) \in \E_i$ whose sink $v_S$ is in $N_{G_i[X]}^{\leq \ell_0-1}[z^*]$, $T \subseteq U_{i+1}-S$, there exists a matching in $G_i$ between $T$ and $(S-\{v_S\}) \cap N_{G_i[X]}^{\leq \ell_0-1}[z^*]$ with size $\lvert T \rvert\}$.
		\item $\E_{i+1,2}=\{(\{v^*\} \cup (S \cap W) \cup T,j+1):$ there exists $(S,j) \in \E_i$ whose sink $v_S$ is in $N_{G_i[X]}^{\leq \ell_0-1}[z^*]$, $T \subseteq U_{i+1}-S$, $(S \cap W) \cup T \neq \emptyset$, there exist $u \in (S-\{v_S\}) \cap N_{G_i[X]}^{\leq \ell_0-1}[z^*]$ and a matching in $G_i$ between $T$ and $(S-\{v_S,u\}) \cap N_{G_i[X]}^{\leq \ell_0-1}[z^*]$ with size $\lvert T \rvert\}$.
		\item $\E_{i+1,3}=\{(\{v^*\} \cup U_{i+1}, 1)\}$.
		\item $\E_{i+1}=\E_{i+1,0} \cup \E_{i+1,1} \cup \E_{i+1,2} \cup \E_{i+1,3}$.
	\end{itemize}

Note that every vertex in $V(G_{i+1})-\{v^*\}$ is in $V(G_i)$. 
And $v^*$ is obtained by contracting a subgraph of $G_i$ with $\lvert N_{G_i[X]}^{\leq \ell_0-1}[z^*] \rvert \leq \sum_{j=0}^{\ell_0-1}d^j \leq d^{\ell_0} \leq N^* \leq N$ vertices in $G_i$ by (i).
So $(G_{i+1},\M_{i+1},D_{i+1})$ satisfies (D6a).
And Claim 1 implies that $(G_{i+1},\M_{i+1},D_{i+1},\E_{i+1})$ satisfies (D6b).

\medskip

\noindent{\bf Claim 2:} If there exist $\A_{i+1}$ and $\A'_{i+1}$ such that $(\E_{i+1},\A_{i+1},\A'_{i+1})$ satisfies (D9) and (D10), then $\E_{i+1}$ satisfies (D5).

\noindent{\bf Proof of Claim 2:}
Since $\E_i$ satisfies (D5), every member of $\E_{i+1,0}$ satisfies (D5). 
And for every $(S,j) \in \E_{i+1}-\E_{i+1,0}$, $S \subseteq \{v^*\} \cup W$, so $\lvert S \rvert \leq \lvert W \rvert+1 \leq r$, and $v^*$ is the unique vertex in $S$ such that $(u,v^*)$ for every $u \in S-\{v^*\}$.
So to show that $\E_{i+1}$ satisfies (D5), it suffices to show that $j \in [h-2]$ for every $(S,j) \in \E_{i+1}-(\E_{i+1,0} \cup \E_{i+1,1})$.
Since $h \geq 3$, it suffices to show that $j \in [h-2]$ for every $(S,j) \in \E_{i+1,2}$.

Let $(S,j) \in \E_{i+1,2}$.
Note that $j \in [h-1]$ since $\E_i$ satisfies (D5).
And note that there exists a non-sink vertex for $(S,j)$ by the definition of $\E_{i+1,2}$. 
So if there exist $\A_{i+1}$ and $\A'_{i+1}$ such that $(\E_{i+1},\A_{i+1},\A'_{i+1})$ satisfies (D9) and (D10), then $\A'_{i+1,(S,j)} \neq \emptyset$, so contracting a member of $\A'_{i+1,(S,j)}$ into a vertex and each member of $\A_{i+1,(S,j)}$ into a vertex creates a $(K_1 \vee (k+h-j)\CT_{j,k})$-minor in $G$, so $G$ has a $\CT_{j+1,k}$-minor.
Since $G$ is $\CT_{h,k}$-minor free, $j+1 \leq h-1$.
That is, $j \in [h-2]$.
$\Box$

\medskip

\noindent{\bf Claim 3:} $\E_{i+1}$ satisfies (D7).

\noindent{\bf Proof of Claim 3:}
Let $(S,j)$ be a member of $\E_i$ with sink $v_S$ such that $S-\{v_S\} \subseteq V(G_{i+1})$ and $M_{i,v_S} \subseteq M_{i+1,v_{S'}}$ for some $v_{S'} \in V(G_{i+1})$. 
Let $S'=\{v_{S'}\} \cup (N_{G_{i+1}}(v_{S'}) \cap S-\{v_S\})$. 
To prove this claim, it suffices to prove $(S',j) \in \E_{i+1}$.

We first assume $v_S \not \in N_{G_i[X]}^{\leq \ell_0-1}[z^*]$.
So $v_{S'}=v_S$.
Recall that $v_S \not \in W$ by Claim 1.
Since every vertex in $S-\{v_S\}$ is adjacent in $G_i$ to $v_S$ by (D4) and (D5), (iii) implies that $S \subseteq V(G_i) \cap V(G_{i+1})$, so $S'=S$ and $(S',j) \in \E_{i+1,0} \subseteq \E_{i+1}$.
So we may assume $v_S \in N_{G_i[X]}^{\leq \ell_0-1}[z^*]$.
In particular, $v_{S'}=v^*$.
Since $S-\{v_S\} \subseteq V(G_{i+1})$, $S-\{v_S\} \subseteq W$ by (iii).
So $S'-\{v^*\}=S \cap W$ and $(S',j) \in \E_{i+1,1}$ (by taking $T=\emptyset$).
$\Box$

\medskip

\noindent{\bf Claim 4:} $(G_{i+1},\M_{i+1},D_{i+1},\E_{i+1})$ satisfies (D8a)-(D8h).

\noindent{\bf Proof of Claim 4:}
Recall that we defined $q_{i+1}=v^*$ and $U_{i+1}^+=W$.
So $q_{i+1} \in V(G_{i+1})$ and $U_{i+1} \subseteq U_{i+1}^+ \subseteq V(G_i) \cap V(G_{i+1}) \cap V(G)$ by Claim 1.

Every vertex $v$ in $V(G_{i})-V(G_{i+1})$ with $\lvert V(M_{i,v}) \rvert=1$ is in $X$, so (D8a) follows from (i).
And (D8b) follows from the definition of $U_{i+1}^+$ and $U_{i+1}$.
Note that $q_{i+1}=v^* \not \in W=U_{i+1}^+$, so (D8c) follows from (iii).
Similarly, (D8d) follows from (i) and (iii).
By (iii), there exists no $v' \in V(G_i)$ with $M_{i,v'} \subseteq M_{i+1,q_{i+1}}$ such that $N_{G_i}(v') -U_{i+1}^+ \neq \emptyset$, so (D8e) holds. 
And (iii) implies $U_{i+1} \cap N_{G_{i+1}}(q_{i+1})= U_{i+1} \cap W = U_{i+1}$, so $((U_{i+1} \cap N_{G_{i+1}}(q_{i+1})) \cup \{v^*\},1) = (U_{i+1} \cup \{v^*\},1) \in \E_{i+1,3} \subseteq \E_{i+1}$, so (D8f) holds.
Since $\lvert Z^* \rvert \geq k+h \geq 2$, $\bigcup_{M \in \M_i}V(M) \neq \bigcup_{M \in \M_{i+1}}V(M)$, so (D8g) holds.

So it suffices to prove (D8h).
(D8ha) clearly holds.
For every $M \in \M_{i}$ disjoint from all members of $\M_{i+1}$, $M=M_{i,v}$ for some $v \in N_{G_i[X]}^{\leq \ell_0-1}[z]$ with $z \in Z^*-\{z^*\}$, so (D8hb) follows from (iii).
And (D8hc) follows from Claim 1; (D8hd) follows from the definition of $G_{i+1}$.
By (i), $\lvert V(G_i)-V(G_{i+1}) \rvert = \sum_{z \in Z^*}\lvert N_{G_i[X]}^{\leq \ell_0-1}[z] \rvert \leq \lvert Z^* \rvert \cdot \sum_{j=0}^{\ell_0-1}d^j \leq (k+h) \cdot d^{\ell_0} \leq N^* \leq N$, so (D8he) holds.

Let $x \in V(G_i) \cap V(G) - \bigcup_{M \in \M_{i+1}}V(M)$.
So $x \in N_{G_i[X]}^{\leq \ell_0-1}[z]$ for some $z \in Z^*-\{z^*\}$, and hence $a_{z,N_G(x) \cap W}=1$.
By the definition of $Z_2$, $a_{z^*,N_G(x) \cap W}=a_{z,N_G(x) \cap W}=1$.
So there exists $x' \in N_{G_i[X]}^{\leq \ell_0-1}[z^*] \cap V(G) \subseteq V(M_{i+1,q_{i+1}}) \cap V(G_i) \cap V(G)$ such that $N_{G}(x') \cap W = N_{G}(x) \cap W$.
Since $U_{i+1} \subseteq W$, (D8hf) holds.

Let $(S,j) \in \E_i$ with $V(M_{i,v_S}) \subseteq \bigcup_{M \in \M_{i}}V(M)-\bigcup_{M \in \M_{i+1}}V(M)$, where $v_S$ is the sink for $(S,j)$.
So $v_S \in N_{G_i[X]}^{\leq \ell_0-1}[z]$ for some $z \in Z^*-\{z^*\}$.
Let $T=S \cap W$.
For every $\alpha \in [2^{\lvert W \rvert}]$, let $j_\alpha$ be the number of vertices $y$ in $S-\{v_S\}$ such that $N_G(y) \cap W$ equals the $\alpha$-th subset of $W$ based on $\sigma$.
So $a_{z,T,j,j_1,...,j_{2^{\lvert W \rvert}}}=1$.
By the definition of $Z_1$, $a_{z^*,T,j,j_1,...,j_{2^{\lvert W \rvert}}}=1$.
Hence there exists $(S',j) \in \E_i$ such that the sink $v_{S'}$ for $(S',j)$ is in $N_{G_i[X]}^{\leq \ell_0-1}[z^*]$, $S' \cap W=T$, and for every $\alpha \in [2^{\lvert W \rvert}]$, there are exactly $j_\alpha$ vertices $y$ in $S'-\{v_{S'}\}$ such that $N_{G}(y) \cap W$ equals the $\alpha$-th subset of $W$ based on $\sigma$.
In particular, $\lvert S' \rvert = 1+\sum_{\alpha \in [2^{\lvert W \rvert}]}j_\alpha = \lvert S \rvert$, $S \cap V(G_{i+1})-\{v_S\} = T = S' \cap V(G_{i+1})-\{v_{S'}\}$ by (iii), and there exists a bijection $\iota$ such that (D8hg) holds.
Therefore, (D8h) holds.
$\Box$

\medskip

\noindent{\bf Claim 5:} $(G_{i+1},\M_{i+1},D_{i+1},\E_{i+1})$ satisfies (D8).

\noindent{\bf Proof of Claim 5:}
By Claim 4, it suffices to show that $(G_{i+1},\M_{i+1},D_{i+1},\E_{i+1})$ satisfies (D8i).
Let $(S,j) \in \E_i$ with $M_{i,v} \cup M_{i,v_S} \subseteq M_{i+1,q_{i+1}}$ for some $v \in S-\{v_S\}$, where $v_S$ is the sink for $(S,j)$.

Let $\C$ be the multiset $\{N_G(x) \cap U_{i+1} \neq \emptyset: x \in S-(\{v_S\} \cup U_{i+1}^+)\}$.
So for every $T \in \C$, there exists $x_T \in S-(\{v_S\} \cup U_{i+1}^+) = S-(\{v_S\} \cup W)$ such that $T=N_G(x_T) \cap U_{i+1}$.
Note that we can choose those $x_T$ such that $x_{T_1} \neq x_{T_2}$ whenever $T_1$ and $T_2$ are distinct members of the multiset $\C$.
Let $f$ be a function that maps each member $T$ of $\C$ to a vertex in $T$.
Then there exists a matching in $G$ between $\{f(T): T \in \C\}$ (as a set) and $\{x_T: T \in \C\}$ with size $\lvert \{f(T): T \in \C\} \rvert$ (as a set).
So there exists a matching in $G$ between $\{f(T): T \in \C\}-S \subseteq U_{i+1}-S$ and $\{x_T: T \in \C\}$ with size $\lvert \{f(T): T \in \C\}-S \rvert$.
Since $U_{i+1} \cup (S-\{v_S\}) \subseteq V(G) \cap V(G_i)$, this matching is also in $G_i$ by (D2).
Note that $\{x_T: T \in \C\} \subseteq S-(\{v_S\} \cup W) \subseteq (S-\{v_S\}) \cap N_{G_i[X]}^{\leq \ell_0-1}[z^*]$ by (iii).
Hence $(\{v^*\} \cup (S \cap W) \cup \{f(T): T \in \C\}, j) = (\{v^*\} \cup (S \cap W) \cup (\{f(T): T \in \C\}-S), j) \in \E_{i+1,1} \subseteq \E_{i+1}$.
So (D8ia) holds.

Let $u \in S-(\{v_S\} \cup U_{i+1}^+)$.
So $u \in S-(\{v_S\} \cup W)$.
Let $\C_u$ be the multiset $\{N_G(x) \cap U_{i+1} \neq \emptyset: x \in S-(\{v_S,u\} \cup U_{i+1}^+)\}$.
So for every $T \in \C_u$, there exists $y_T \in S-(\{v_S,u\} \cup U_{i+1}^+) = S-(\{v_S,u\} \cup W)$ such that $T=N_G(y_T) \cap U_{i+1}$.
Note that we can choose those $y_T$ such that $y_{T_1} \neq y_{T_2}$ whenever $T_1$ and $T_2$ are distinct members of the multiset $\C_u$.
Assume either $S \cap U_{i+1}^+ \neq \emptyset$ or $\C_u \neq \emptyset$.
Let $f_u$ be a function that maps each member $T$ of $\C_u$ to a vertex in $T$.
Then there exists a matching in $G$ between $\{f_u(T): T \in \C_u\}$ (as a set) and $\{y_T: T \in \C_u\}$ with size $\lvert \{f_u(T): T \in \C_u\} \rvert$ (as a set).
Let $T^*$ be the set $\{f_u(T): T \in \C_u\}-(S \cap W)$.
So there exists a matching in $G$ between $T^* \subseteq U_{i+1}-S$ and $\{y_T: T \in \C_u\}$ with size $\lvert T^* \rvert$.
Since $U_{i+1} \cup (S-\{v_S\}) \subseteq V(G) \cap V(G_i)$, this matching is also in $G_i$ by (D2).
Note that $\{y_T: T \in \C_u\} \subseteq S-(\{v_S,u\} \cup W) \subseteq (S-\{v_S,u\}) \cap N_{G_i[X]}^{\leq \ell_0-1}[z^*]$ by (iii).
Since either $S \cap W = S \cap U_{i+1}^+ \neq \emptyset$ or $\C_u \neq \emptyset$, we know $(S \cap W) \cup T^* = (S \cap W) \cup \{f_u(T): T \in \C_u\} \neq \emptyset$.
Hence $(\{v^*\} \cup (S \cap W) \cup \{f_u(T): T \in \C_u\}, j+1) = (\{v^*\} \cup (S \cap W) \cup T^*, j+1) \in \E_{i+1,2} \subseteq \E_{i+1}$.
Therefore, (D8i) holds and hence (D8) holds.
$\Box$

\medskip

For every $Q \in \E_{i+1,0}$, we know $Q \in \E_i$, and we define $\A_{i+1,Q}$ and $\A'_{i+1,Q}$ to be $\A_{i,Q}$ and $\A'_{i,Q}$, respectively.
Since $(\A_i,\A'_i)$ satisfies (D10), (D10) holds for every member of $\E_{i+1,0}$.

For every $Q \in \E_{i+1,1}$, we define the following:
	\begin{itemize}
		\item Let $(S_Q,j_Q)$ be a member of $\E_i$ such that $Q=(\{v^*\} \cup (S_Q \cap W) \cup T_Q, j_Q)$, where $T_Q \subseteq U_{i+1}-S_Q$, $v_{S_Q}$ is the sink for $(S_Q,j_Q)$ and is in $N_{G_i[X]}^{\leq \ell_0-1}[z^*]$, and there exists a matching $M_{T_Q}$ in $G_i$ between $T_Q$ and $(S_Q-\{v_{S_Q}\}) \cap N_{G_i[X]}^{\leq \ell_0-1}[z^*]$ with size $\lvert T_Q \rvert$.
		\item Let $\A_{i+1,Q}=\A_{i,(S_Q,j_Q)}$.
		\item For every vertex $v \in S_Q-\{v_{S_Q}\}$, let $A_v$ be the member of $\A'_{i,(S_Q,j_Q)}$ such that $v \in V(A_v)$.
		\item For every vertex $u \in T_Q$, let $u'$ be the vertex in $(S_Q-\{v_{S_Q}\}) \cap N_{G_i[X]}^{\leq \ell_0-1}[z^*] \subseteq V(M_{i+1,v^*})$ matched with $u$ in $M_{T_Q}$, and let $A_{i+1,u}=G[V(A_{u'}) \cup \{u\}]$.
		\item Let $\A'_{i+1,Q}=\{A_v: v \in S_Q \cap W\} \cup \{A_{i+1,u}: u \in T_Q\}$.
	\end{itemize}
Clearly, (D10) holds for every member of $\E_{i+1,1}$.

For every $Q \in \E_{i+1,2}$, we define the following:
	\begin{itemize}
		\item Let $(S_Q,j_Q)$ be a member of $\E_i$ such that $Q=(\{v^*\} \cup (S_Q \cap W) \cup T_Q, j_Q+1)$, where $T_Q \subseteq U_{i+1}-S_Q$, $v_{S_Q}$ is the sink for $(S_Q,j_Q)$ and is in $N_{G_i[X]}^{\leq \ell_0-1}[z^*]$, and there exist $u_Q \in (S_Q-\{v_{S_Q}\}) \cap N_{G_i[X]}^{\leq \ell_0-1}[z^*]$ and a matching $M_{T_Q}$ in $G_i$ between $T_Q$ and $(S_Q-\{v_{S_Q},u_Q\}) \cap N_{G_i[X]}^{\leq \ell_0-1}[z^*]$ with size $\lvert T_Q \rvert$.
		\item For every $z \in Z^*-\{z^*\}$, by the definition of $Z_1$, there exists $(S_z,j_Q) \in \E_i$ with sink $v_{S_z} \in N_{G_i[X]}^{\leq \ell_0-1}[z]$ such that $S_z \cap W=S_Q \cap W$, and there exists a bijection $\iota_z$ between $(S_Q-\{v_{S_Q}\}) \cap N_{G_i[X]}^{\leq \ell_0-1}[z^*]$ and $(S_z-\{v_{S_z}\}) \cap N_{G_i[X]}^{\leq \ell_0-1}[z]$ such that $N_G(v) \cap W=N_G(\iota_z(v)) \cap W$ for every $v \in (S_Q-\{v_{S_Q}\}) \cap N_{G_i[X]}^{\leq \ell_0-1}[z^*]$. 
		\item Let $Z_Q$ be a subset of $Z^*-\{z^*\}$ with size $k+h-j_Q-1$. (Note that $Z_Q$ exists since $\lvert Z^* \rvert \geq k+h$ and $j_Q \geq 1$.)
		\item For every $z \in Z_Q$, define the following:
			\begin{itemize}
				\item For every $v \in S_z-\{v_{S_z}\}$, let $A_{z,v}$ be the member of $\A'_{i,(S_z,j_Q)}$ with $v \in V(A_{z,v})$.
				\item Let $\A^{(1)}_{i,(S_z,j_Q)}$ be a subset of $\A_{i,(S_z,j_Q)}$ such that contracting each member of $\A^{(1)}_{i,(S_z,j_Q)}$ into a vertex creates a $\CT_{j_Q,k}$-minor.

					(Note that $\A^{(1)}_{i,(S_z,j_Q)}$ exists since $\A_{i,(S_z,j_Q)}$ satisfies (D10).)
				\item Let $A_{i+1,z}=G[V(A_{z,\iota_z(u_Q)}) \cup \bigcup_{A \in \A^{(1)}_{i,(S_z,j_Q)}}V(A)]$.
				\item Let $\A^{(2)}_{i,(S_z,j_Q)}$ be a subset of $\A_{i,(S_z,j_Q)}-\A^{(1)}_{i,(S_z,j_Q)}$ such that contracting each member of $\A^{(2)}_{i,(S_z,j_Q)}$ into a vertex creates a $k\CT_{j_Q,k}$-minor.

					(Note that $\A^{(2)}_{i,(S_z,j_Q)}$ exists since $\A_{i,(S_z,j_Q)}$ satisfies (D10) and $j_Q \in [h-2]$ by (D5).)
				\item Let $\A_{i+1,z} = \{A_{i+1,z}\} \cup \A^{(2)}_{i,(S_z,j_Q)}$.

			(Note that since $\A_{i,(S_z,j_Q)}$ and $\A'_{i,(S_z,j_Q)}$ satisfy (D10), contracting each member of $\A_{i+1,z}$ into a vertex creates a $(K_1 \vee k\CT_{j_Q,k})$-minor and hence a $\CT_{j_Q+1,k}$-minor.)
			\end{itemize}
		\item Let $\A_{i+1,Q}=\bigcup_{z \in Z_Q}\A_{i+1,z}$.

			(Note that members of $\A_{i+1,Q}$ are pairwise disjoint since $(\A_i,\A'_i)$ satisfies (D11).
			Since $\lvert Z_Q \rvert=k+h-j_Q-1$, contracting each member of $\A_{i+1,Q}$ into a vertex creates a $(k+h-j_Q-1)\CT_{j_Q+1,k}$-minor.)
		\item For every vertex $w \in T_Q$, there exists $w_Q \in S_Q-\{v_{S_Q}\}$ matched with $w$ in $M_{T_Q}$, and we let $A_{i+1,w}=G[\{w\} \cup \bigcup_{z \in Z_Q}V(A_{z,\iota_z(w_Q)})]$. 

			(Note that the definition of $\iota_z$ implies that $\iota_z(w_Q)$ is adjacent in $G$ to $w$, so $A_{i+1,w}$ is connected.)
		\item Let $\A'_{i+1,Q}=\{G[V(\bigcup_{z \in Z_Q}A_{z,v})]: v \in S_Q \cap W\} \cup \{A_{i+1,w}: w \in T_Q\}$.
	\end{itemize}
Since $(\A_i,\A'_i)$ satisfies (D10) and (D11), we know that (D10) holds for every member of $\E_{i+1,2}$.

For the unique member $(\{v^*\} \cup U_{i+1}, 1)$ of $\E_{i+1,3}$, define $\A_{i+1,(\{v^*\} \cup U_{i+1},1)} = \linebreak \{G[\bigcup_{v \in N_{G_i[X]}^{\leq \ell_0-1}[z]}V(M_{i,v})]: z \in Z^*-\{z^*\}\}$, and $\A'_{i+1,(\{v^*\} \cup U_{i+1},1)} = \{G[\{v\}]: v \in U_{i+1}\}$.
Clearly, (D10) holds for $(\{v^*\} \cup U_{i+1},1)$.

Therefore, $(\E_{i+1},\A_{i+1},\A'_{i+1})$ satisfies (D9) and (D10).
Hence $\E_{i+1}$ satisfies (D5) by Claim 2.

\medskip

\noindent{\bf Claim 6:} $(\A_{i+1},\A'_{i+1})$ satisfies (D11).

\noindent{\bf Proof of Claim 6:}
Suppose to the contrary that there exist members $Q_1=(S_1,j_1)$ and $Q_2=(S_2,j_2)$ of $\E_{i+1}$ with distinct sinks violating (D11).
Since $(\A_i,\A'_i)$ satisfies (D11), we may assume $Q_1 \not \in \E_{i+1,0}$ by symmetry.
Hence $Q_1 \in \bigcup_{\alpha=1}^3\E_{i+1,\alpha}$, so the sink of $Q_1$ is $v^*$.
Since $Q_1$ and $Q_2$ have distinct sink, $Q_2 \in \E_{i+1,0}$.
In particular, by Claim 1, the sink for $Q_2$, denoted by $v_{Q_2}$, is in $V(G_i)-(W \cup \bigcup_{z \in Z^*}N_{G_i[X]}^{\leq \ell_0-1}[z])$, and $S_2 \cap \bigcup_{z \in Z^*}N_{G_i[X]}^{\leq \ell_0-1}[z]=\emptyset$.

Suppose to the contrary that $Q_1 \in \E_{i+1,3}$.
Then every member of $\A_{i+1,Q_1}$ is contained in $G[\bigcup_{z \in Z^*-\{z^*\}}\bigcup_{v \in N_{G_i[X]}^{\leq \ell_0-1}[z]}V(M_{i,v})] \subseteq G[V(G)-\bigcup_{v \in V(G_{i+1})}M_{i+1,v}]$ by definition.
Since $v_{Q_2} \in V(G_i)-(W \cup \bigcup_{z \in Z^*}N_{G_i[X]}^{\leq \ell_0-1}[z])$ and $Q_2 \in \E_{i+1,0}$, $V(M_{i+1,v_{Q_2}})=V(M_{i,v_{Q_2}})$ is disjoint from every member of $\A_{i+1,Q_1}$.
And $V(G)-\bigcup_{v \in V(G_{i})}V(M_{i,v})$ is disjoint from $\bigcup_{z \in Z^*-\{z^*\}}\bigcup_{v \in N_{G_i[X]}^{\leq \ell_0-1}[z]}V(M_{i,v})$.
By (D10), every member of $\A_{i+1,Q_2}=\A_{i,Q_2}$ is contained in $G[V(M_{i,v_{Q_2}}) \cup (V(G)-\bigcup_{v \in V(G_i)}M_{i,v})]$ and hence is disjoint from every member of $\A_{i+1,Q_1}$.
And every member of $\A'_{i+1,Q_2}=\A'_{i,Q_2}$ is contained in $G[\bigcup_{v \in S_2}V(M_{i,v}) \cup (V(G)-\bigcup_{v \in V(G_i)}V(M_{i,v})]$ (by (D10)) and hence is disjoint from every member of $\A_{i+1,Q_1}$ (since $S_2 \cap \bigcup_{z \in Z^*}N_{G_i[X]}^{\leq \ell_0-1}[z]=\emptyset$).
Moreover, every member of $\A'_{i+1,Q_1}$ consists of a vertex $a'$ in $U_{i+1} \subseteq W \subseteq V(G_i) \cap V(G_{i+1}) \cap V(G)$ by Claim 1, so it is disjoint from every member of $\A_{i+1,Q_2}=\A_{i,Q_2}$ by (D10) (since $v_{Q_2} \not \in W$); and it intersects an member of $\A'_{i+1,Q_2}=\A'_{i,Q_2}$ only possibly at $a' \in V(S_1)$, and if it happens, $a' \in V(S_1) \cap V(S_2)$ (since $v_{Q_2} \not \in W$ and $a' \in W \subseteq V(G) \cap V(G_i) \cap V(G_{i+1})$ by Claim 1).
So $Q_1$ and $Q_2$ do not violate (D11), a contradiction.

So $Q_1 \in \E_{i+1,1} \cup \E_{i+1,2}$.
Hence $(S_{Q_1},j_{Q_1}) \in \E_i$ is defined; when $Q_1 \in \E_{i+1,2}$, $Z_{Q_1}$ and $(S_z,j_{Q_1}) \in \E_i$ are defined for every $z \in Z_{Q_1}$.
Note that the sinks for $(S_{Q_1},j_{Q_1})$ and $(S_z,j_{Q_1})$ are in $\bigcup_{z \in Z^*}N_{G_i[X]}^{\leq \ell_0-1}[z]$ but the sink for $Q_2$ is not.
Since $(\A_i,\A_i')$ satisfies (D10) and (D11) and $S_2 \cap \bigcup_{z \in Z^*}N_{G_i[X]}^{\leq \ell_0-1}[z]=\emptyset$, every member of $\A_{i+1,Q_1}$ is disjoint from every member of $\A_{i+1,Q_2} \cup \A'_{i+1,Q_2}$.
Moreover, for every $A \in \A'_{i+1,Q_1}$, $V(A) \subseteq W \cup \bigcup_{z \in Z^*}\bigcup_{v \in N_{G_i[X]}^{\leq \ell_0-1}[z]}V(M_{i,v}) \cup (V(G)-\bigcup_{v \in V(G_i)}V(M_{i,v}))$ and $V(A) \cap W=V(A) \cap W \cap S_1$ is a set with size at most 1.
Hence every member of $\A'_{i+1,Q_1}$ is disjoint from every member of $\A_{i+1,Q_2}$, and if some member $A_1 \in \A'_{i+1,Q_1}$ intersects some member $A_2 \in \A'_{i+1,Q_2}$, then $V(A_1) \cap V(A_2) \subseteq V(A_1) \cap V(A_2) \cap W \subseteq (V(A_1) \cap S_{1}) \cap (V(A_2) \cap W \cap V(G_{i+1})) \subseteq V(A_1) \cap S_{1} \cap V(A_2) \cap S_{2}$ by Claim 1, contradicting that $Q_1$ and $Q_2$ are counterexamples.
$\Box$

\medskip

To prove this lemma, it suffices to show that $(G_{i+1},\M_{i+1},D_{i+1},\E_{i+1})$ satisfies (D12).
Let $v \in V(G_{i+1})$ such that either $\lvert V(M_{i+1,v}) \rvert \geq 2$, or $v$ is a head with respect to $D_{i+1}$ or the sink for some member of $\E_{i+1}$.
It suffices to show $\deg_{G_{i+1}}(v) \leq r$.

Suppose to the contrary that $\deg_{G_{i+1}}(v) > r$.
By (iii), $\deg_{G_{i+1}}(v^*)=\lvert W \rvert \leq r-1$.
So $v \neq v^*$.
Hence $v \in V(G_i) \cap V(G_{i+1})$.
In particular, $\deg_{G_{i+1}}(v) \leq \deg_{G_i}(v)$ and $M_{i,v}=M_{i+1,v}$.
So $\deg_{G_i}(v)>r$.
Since $(G_{i},\M_i,D_{i},\E_{i})$ satisfies (D12), $\lvert V(M_{i+1,v}) \rvert = \lvert V(M_{i,v}) \rvert=1$, $v$ is not a head with respect to $D_{i}$, and $v$ is not the sink for some member of $\E_{i}$.
Hence by the definition of $D_{i+1}$, since $v \neq v^*$, $v$ is not a head with respect to $D_{i+1}$.
So $v$ is a sink for some member of $\E_{i+1}$.
By the definition of $\E_{i+1}$, $v$ is the sink for some member of $\E_{i+1,0}$ since $v \neq v^*$.
However, $\E_{i+1,0} \subseteq \E_i$.
So $v$ is the sink for some member of $\E_i$, a contradiction.
This proves the lemma.
\end{pf}

\begin{lemma} \label{contract}
For any positive integers $h \geq 3,k$, there exist positive integers $r=r(h,k),d=d(h,k),N=N(h,k)$ such that the following hold.
Let $G$ be a graph with no $\CT_{h,k}$-minor.
Let $i$ be a positive integer.
Let $((G_\alpha,\M_\alpha,\E_\alpha,D_\alpha,\A_\alpha,\A_\alpha'): \alpha \in [i])$ be a strong $i$-$(G,h,k,r,d,N)$-defective elimination scheme.
If $\lvert V(G_i) \rvert>N$, then there exist a graph $G_{i+1}$ with $\lvert V(G_{i+1}) \rvert < \lvert V(G_i) \rvert$ and a tuple $(G_{i+1},\M_{i+1},\E_{i+1},D_{i+1},\A_{i+1}, \allowbreak \A_{i+1}')$ such that the sequence $((G_\alpha,\M_\alpha,\E_\alpha,D_\alpha,\A_\alpha,\A'_\alpha): \alpha \in [i+1])$ is a strong $(i+1)$-$(G,h,k,r,d,N)$-defective elimination scheme.
\end{lemma}

\begin{pf}
Let $h \geq 3,k$ be positive integers.
Let $\G$ be the class of $\CT_{h,k}$-minor-free graphs.
Define the following:
	\begin{itemize}
		\item Define $r=r_{\ref{homo}}(\G)$, where $r_{\ref{homo}}$ is the integer $r$ mentioned in Lemma \ref{homo}.
		\item Let $t_0=(h-2)(r+1)2^{r-1}r^{2^{r-1}}$. 
		\item Let $t_1=3^{r+t_0}$.
		\item Let $k_0=1+(h+k-1)(6t_1+1)$.
		\item Let $\ell_0 = n_{\ref{geodesic}}(t_1,k_0,3)+1$, where $n_{\ref{geodesic}}$ is the integer $n$ mentioned in Lemma \ref{geodesic}.
		\item Let $t=t_{\ref{del}}(h,k,r)$, where $t_{\ref{del}}$ is the integer $t$ mentioned in Lemma \ref{del}.
		\item Define $d=d_{\ref{homo}}(\G,t,\ell_0)+1$, where $d_{\ref{homo}}$ is the integer $d$ mentioned in Lemma \ref{homo}.
		\item Let $N_1=N_{\ref{homo}}(\G,t,\ell_0)$, where $N_{\ref{homo}}$ is the integer $N$ mentioned in Lemma \ref{homo}.
		\item Let $N_2 = N_{\ref{del}}(h,k,r,d,\ell_0)$, where $N_{\ref{del}}$ is the integer $N^*$ mentioned in Lemma \ref{del}.
		\item Define $N=d^{\ell_0}+N_1+N_2$.
	\end{itemize}

Let $G$ be a graph with no $\CT_{h,k}$-minor.
Let $i$ be a positive integer.
Let $((G_\alpha,\M_\alpha,\E_\alpha,D_\alpha, \allowbreak \A_\alpha,\A_\alpha'): \alpha \in [i])$ be a strong $i$-$(G,h,k,r,d,N)$-defective elimination scheme. 
Assume $\lvert V(G_i) \rvert >N$.

By (D1) and (D2), $G_i$ is a minor of $G$.
So $G_i$ is $\CT_{h,k}$-minor-free.
That is, $G_i \in \G$.
Since $\lvert V(G_i) \rvert > N \geq N_1$, by Lemma \ref{homo} (by taking $G=G_i$), there exist $X,Z,W \subseteq V(G_i)$ with $Z \subseteq X$, $\lvert Z \rvert=t$, $W \subseteq V(G_i)-X$ and $\lvert W \rvert \leq r-1$ such that
	\begin{itemize}
		\item[(i)] every vertex in $X$ has degree in $G_i$ at most $d$,
		\item[(ii)] for any distinct $z,z'$ in $Z$, the distance in $G_i[X]$ between $z,z'$ is at least $2\ell_0-1$, and
		\item[(iii)] $N_{G_i}(N_{G_i[X]}^{\leq \ell_0-1}[z])-X=W$ for every $z \in Z$.
	\end{itemize}
Since $\lvert V(G_i) \rvert > N \geq N_2$, if $N_{G_i}(N_{G_i[X]}^{\leq \ell_0-1}[z]) \subseteq W$ for every $z \in Z$, then (iii) implies that $N_{G_i}(N_{G_i[X]}^{\leq \ell_0-1}[z]) = W$ for every $z \in Z$, so we are done by Lemma \ref{del}, (i) and (ii). 

Hence we may assume that there exists $z^* \in Z$ such that $N_{G_i}(N_{G_i[X]}^{\leq \ell_0-1}[z]) \not \subseteq W$.
By (iii), $N_{G_i[X]}^{\leq \ell_0}[z^*]-N_{G_i[X]}^{\leq \ell_0-1}[z^*] \neq \emptyset$.
So there exists a geodesic in $G_i[X]$ starting from $z^*$ with length $\ell_0-1$.
Note that this geodesic contains a subpath $P_{z^*}$ starting from $z^*$ on $\ell_0-1-3t_1$ vertices, and $P_{z^*}$ is also a geodesic in $G_i[N_{G_i[X]}^{\leq \ell_0-1}[z^*]]$.

For every $(S,j) \in \E_i$, the {\it type} of $(S,j)$ is defined to be the sequence $$(j,\lvert S \rvert, S \cap W, \lvert \{v \in S-\{v_S\}: N_{G_i}(v) \cap W=T\} \rvert: T \subseteq W),$$ where $v_S$ is the sink for $(S,j)$.
By (D5), there are at most $(h-2) \cdot (r+1) \cdot 2^{\lvert W \rvert} \cdot r^{2^{\lvert W \rvert}} \leq (h-2)(r+1)2^{r-1}r^{2^{r-1}} \leq t_0$ different types for members of $\E_i$.
For every $v \in N_{G_i[X]}^{\leq \ell_0-1}[z^*]$, let $\phi(v)=(N_{G_i}(v) \cap W, a_\tau: \tau$ is a type for members of $\E_i)$, where for every $\tau$, 
	\begin{itemize}
		\item $a_\tau=0$ if $v$ is not in $S$ for every member $(S,j)$ of $\E_i$ with type $\tau$; 
		\item $a_\tau=1$ if $v$ is the sink for some member $(S,j)$ of $\E_i$ with type $\tau$; 
		\item $a_\tau=2$ otherwise.
	\end{itemize}
Let $Y_0$ be the image of $\phi$.
Note that $\lvert Y_0 \rvert \leq 2^{\lvert W \rvert} \cdot 3^{t_0} \leq 3^{r+t_0}=t_1$.
We call $\phi(v)$ the {\it type} of $v$.

Let $H=G_i[N_{G_i[X]}^{\leq \ell_0-1}[z^*]]$.
Note that $H$ is a graph such that $V(H) = N_{G_i[X]}^{\leq \ell_0-1}[z^*] = N_{H}^{\leq \ell_0-1}[z^*]$, and $\phi$ is a function from $N_{G_i[X]}^{\leq \ell_0-1}[z^*]=V(H)$ to the set $Y_0$ with size at most $t_1$, and $P_{z^*}$ is a geodesic in $G_i[N_{G_i[X]}^{\leq \ell_0-1}[z^*]]=H$ on $\ell_0-1-3t_1$ vertices.	
So by Lemma \ref{geodesic} (with taking $G=H$), there exist a subpath $P'_{z^*}$ of $P_{z^*}$ and a nonempty set $Y \subseteq Y_0$ such that
	\begin{itemize}
		\item[(iv)] for every $y \in Y_0-Y$, $\phi^{-1}(\{y\})$ is disjoint from $N_{H}^{\leq 3\lvert Y \rvert}[V(P'_{z^*})]= N_{G_i[N_{G_i[X]}^{\leq \ell_0-1}[z^*]]}^{\leq 3\lvert Y \rvert}[V(P'_{z^*})]$, and
		\item[(v)] $P'_{z^*}$ can be partitioned into $k_0$ disjoint subpaths $P'_{z^*,1},...,P'_{z^*,k_0}$ such that for every $\alpha \in [k_0]$ and $y \in Y$, $\phi^{-1}(\{y\})$  intersects $N_{H}^{\leq 3(\lvert Y \rvert-1)}[V(P'_{{z^*},\alpha})] = N_{G_i[N_{G_i[X]}^{\leq \ell_0-1}[z^*]]}^{\leq 3(\lvert Y \rvert-1)}[V(P'_{{z^*},\alpha})]$.
	\end{itemize}
Note that $P'_{z^*}$ is a subpath of $P_{z^*}$, and $P_{z^*}$ starts at $z^*$ and has at most $\ell_0-1-3t_1$ vertices.
So 
	\begin{itemize}
		\item[(vi)] for every $\alpha \in [3t_1]$ and for every subpath $P''$ of $P'_{z^*}$, we have $N_{H}^{\leq \alpha}[V(P'')] = \linebreak N_{G_i[N_{G_i[X]}^{\leq \ell_0-1}[z^*]]}^{\leq \alpha}[V(P'')] = N_{G_i[X]}^{\leq \alpha}[V(P'')] \subseteq N_{G_i[X]}^{\leq \ell_0-1}[z^*]$.
	\end{itemize}

Let 
$$O_0=N_{G_i[X]}^{\leq 3(\lvert Y \rvert-1)+1}[V(P'_{z^*})].$$
Let 
$$O=O_0 \cup \{x \in N_{G_i[X]}(O_0): \lvert V(M_{i,x}) \rvert \geq 2, \text{ or } x \text{ is the sink for some member of } \E_i\}.$$
Note that $O \subseteq N_{G_i[X]}^{\leq 3\lvert Y \rvert-1}[V(P'_{z^*})]$ and $N_{G_i[X]}[O] \subseteq N_{G_i[X]}^{\leq 3\lvert Y \rvert}[V(P'_{z^*})]$. 
So (iv), (v) and (vi) imply that 
	\begin{itemize}
		\item[(vii)] for every $y \in Y_0-Y$, $\phi^{-1}(\{y\})$ is disjoint from $N_{G_i[X]}[O]$, and
		\item[(viii)] $P'_{z^*}$ can be partitioned into $k_0$ disjoint subpaths $P'_{z^*,1},...,P'_{z^*,k_0}$ such that for every $\alpha \in [k_0]$ and $y \in Y$, $\phi^{-1}(\{y\})$ intersects $O \cap N_{G_i[X]}^{\leq 3\lvert Y \rvert}[V(P'_{{z^*},\alpha})]$.
	\end{itemize}

Define $G_{i+1}$ to be the graph obtained from $G_i$ by contracting $O$ into a new vertex $v^*$, deleting resulting parallel edges and loops, and deleting all edges between $v^*$ and $N_{G_i}(O) \cap X = N_{G_i[X]}(O)$.
Define $\M_{i+1} = \{G[\bigcup_{v \in O}V(M_{i,v})]\} \cup (\M_i-\{M_{i,v}: v \in O\})$.
Note that $M_{i+1,v^*}=G[\bigcup_{v \in O}V(M_{i,v})]$.

Since $Y \neq \emptyset$, $\lvert O \rvert \geq 2$, so $\lvert V(G_{i+1}) \rvert < \lvert V(G_i) \rvert$.
Hence $(G_{i+1},\M_{i+1})$ satisfies (D1) and (D3).
Since $(G_i,\M_i)$ satisfies (D2) and $\lvert V(M_{i+1,v^*}) \rvert \geq \lvert O \rvert \geq 2$, $(G_{i+1},\M_{i+1})$ satisfies (D2).

Define $q_{i+1}=v^*$.
Define $U_{i+1}^+=N_{G_i}(O)-X$.

\medskip

\noindent{\bf Claim 1:} The following statements hold.
	\begin{itemize}
		\item $U_{i+1}^+ = N_{G_i}(O) \cap W$ and $U_{i+1}^+ \subseteq W \cap V(G_i) \cap V(G_{i+1})$.
		\item For any $v \in O$ and $T \subseteq N_{G_i}(v)$, $T \cap W = T \cap U^+_{i+1}$.
		\item For every $w \in W$, $\lvert V(M_{i,w}) \rvert =1$ and $w$ is not a head (with respect to $D_i$) and is not the sink for a member of $\E_i$.
	\end{itemize}

\noindent{\bf Proof of Claim 1:}
Since $W \subseteq V(G_i)-X$, $U^+_{i+1}=N_{G_i}(O)-X \supseteq N_{G_i}(O) \cap W$.
And since $O \subseteq N_{G_i[X]}^{\leq \ell_0-1}[z^*]$, $U^+_{i+1}=N_{G_i}(O)-X \subseteq (W \cup X)-X \subseteq W$ by (iii).
Since $U^+_{i+1} \subseteq N_{G_i}(O)$, $U^+_{i+1} \subseteq N_{G_i}(O) \cap W$.
Hence $U^+_{i+1} = N_{G_i}(O) \cap W$.
In particular, $U^+_{i+1} \subseteq W \cap V(G_i)$.
Since $U_{i+1}^+ \cap O=\emptyset$, $U_{i+1}^+ \subseteq V(G_i) \cap V(G_{i+1})$.
This proves the first statement of this claim.

Let $v \in O$.
Let $T \subseteq N_{G_i}(v)$.
So $T \subseteq N_{G_i}[O]$.
Since $O \subseteq X$, $O \cap W=\emptyset$.
So $T \cap W = T \cap W-O \subseteq N_{G_i}(O) \cap W = U_{i+1}^+$ by the first statement of this claim.
Hence $T \cap W \subseteq T \cap U^+_{i+1}$.
By the first statement of this claim, $U^+_{i+1} \subseteq W$, so $T \cap U^+_{i+1} \subseteq T \cap W$.
So $T \cap U^+_{i+1} = T \cap W$.
This proves the second statement of this claim.

And by (ii) and (iii), for every $w \in W$, $\deg_{G_i}(w) \geq \lvert Z \rvert = t > r$, so $\lvert V(M_{i,w}) \rvert =1$ and $w$ is not a head (with respect to $D_i$) and is not the sink for a member of $\E_i$ since $(G_i,\M_i,D_i,\E_i)$ satisfies (D12).
$\Box$

\medskip

Note that $U_{i+1}^+ = N_{G_{i+1}}(v^*) = N_{G_{i+1}}(q_{i+1})$.
Define $D_{i+1} = \{(u,v) \in D_i: \{u,v\} \subseteq V(G_i) \cap V(G_{i+1})\} \cup \{(u,v^*): u \in U_{i+1}^+\}$.
Since $D_i$ satisfies (D4), Claim 1 implies that $D_{i+1}$ satisfies (D4).

By possibly changing the indices, we may assume that $P'_{z^*,1},...,P'_{z^*,k_0}$ appear in $P'_{z^*}$ in the order listed.
Since $P'_{z^*,1},...,P'_{z^*,k_0}$ are pairwise disjoint subpaths of a geodesic in $G_i[X]$, we know that for any $\alpha,\beta \in [k_0]$, the distance in $G_i[X]$ between $V(P'_{z^*,\alpha})$ and $V(P'_{z^*,\beta})$ is at least $\lvert \beta-\alpha \rvert$.

For every $\alpha \in [h+k]$, let $O_\alpha = O \cap N_{G_i[X]}^{\leq 3\lvert Y \rvert}[V(P'_{{z^*},1+(\alpha-1)(6t_1+1)})]$.
Since $k_0 \geq 1+(h+k-1)(6t_1+1)$, all $O_\alpha$'s are well-defined.
Since $\lvert Y \rvert \leq t_1$, we know that $O_1,O_2,..., \allowbreak O_{h+k}$ are pairwise disjoint.

\medskip

\noindent{\bf Claim 2:} If $(S,j)$ is a member of $\E_i$ with sink $v_S$ such that $v_S \in O$, then for every $\alpha \in [h+k]$, 
	\begin{itemize}
		\item there exists a member $(S_\alpha,j)$ of $\E_i$ with $S_\alpha \subseteq O_\alpha \cup (N_{G_i}(O_\alpha) \cap U^+_{i+1})$, $S_\alpha \cap U_{i+1}^+=S \cap U^+_{i+1}$, and $v_{S_\alpha} \in O_\alpha$, where $v_{S_\alpha}$ is the sink for $(S_\alpha,j)$, and 
		\item there exists a bijection $\iota_{(S,j),\alpha}: S-(\{v_S\} \cup U^+_{i+1}) \rightarrow S_\alpha-(\{v_{S_\alpha}\} \cup U^+_{i+1})$ such that for every $v \in S-(\{v_S\} \cup U^+_{i+1})$, $N_{G_i}(\iota_{(S,j),\alpha}(v)) \cap U^+_{i+1} = N_{G_i}(v) \cap U^+_{i+1}$. 
	\end{itemize}

\noindent{\bf Proof of Claim 2:}
Let $(S,j)$ be a member of $\E_i$ with sink $v_S$ such that $v_S \in O$.
By (vii), $\phi(v_S) \in Y$.
Let $\alpha \in [h+k]$.
By (v) and (vi), there exists $v_\alpha \in N_{G_i[X]}^{\leq 3(\lvert Y \rvert-1)}[V(P'_{z^*,1+(\alpha-1)(6t_1+1)})] \subseteq O_\alpha$ such that $\phi(v_\alpha)=\phi(v_S)$.
By the definition of $\phi$, $v_\alpha$ is the sink for some member $(S_\alpha,j_\alpha)$ of $\E_i$ such that the type of $(S_\alpha,j_\alpha)$ equals the type of $(S,j)$.
In particular, $j_\alpha=j$, $\lvert S_\alpha \rvert = \lvert S \rvert$, and $S_\alpha \cap W = S \cap W$.
Since every vertex in $S_\alpha-\{v_\alpha\}$ is adjacent in $G_i$ to $v_\alpha$, $S_\alpha \subseteq N_{G_i}[v_\alpha] \subseteq N_{G_i[X]}^{\leq 3(\lvert Y \rvert-1)+1}[V(P'_{z^*,1+(\alpha-1)(6t_1+1)})] \cup (N_{G_i}[v_\alpha]-X) \subseteq O_\alpha \cup (N_{G_i}(O_\alpha)-X) \subseteq O_\alpha \cup (N_{G_i}(O_\alpha) \cap U^+_{i+1})$.

Since $v_S \in O$ and $S-\{v_S\} \subseteq N_{G_i}(v_S)$, the second statement of Claim 1 implies that $(S-\{v_S\}) \cap W = (S-\{v_S\}) \cap U^+_{i+1}$. 
Since $v_S \in O \subseteq X$, $v_S \not \in W \cup U^+_{i+1}$, so $S \cap W = S \cap U^+_{i+1}$.
Similarly, since $v_\alpha \in O_\alpha \subseteq O$ and $S_\alpha-\{v_\alpha\} \subseteq N_{G_i}(v_\alpha)$, we have $S_\alpha \cap W = S_\alpha \cap U^+_{i+1}$.
Since $S \cap W = S_\alpha \cap W$, we have $S \cap U^+_{i+1} = S_\alpha \cap U^+_{i+1}$.

By the definition of the type of members of $\E_i$, for every $T \subseteq W$, $\lvert \{v \in S-\{v_S\}: N_{G_i}(v) \cap W=T\} \rvert = \lvert \{v \in S_\alpha-\{v_\alpha\}: N_{G_i}(v) \cap W=T\} \rvert$.
So this together with $S \cap U^+_{i+1} = S_\alpha \cap U^+_{i+1}$ imply that there exists a bijection $\iota_{(S,j),\alpha}: S-(\{v_S\} \cup U^+_{i+1}) \rightarrow S_\alpha-(\{v_\alpha\} \cup U^+_{i+1})$ such that for every $v \in S-(\{v_S\} \cup U^+_{i+1})$, $N_{G_i}(v) \cap W = N_{G_i}(\iota_{(S,j),\alpha}(v)) \cap W$.
Since $S_\alpha-U^+_{i+1} \subseteq O_\alpha \subseteq O$, for every $v \in S_\alpha-U^+_{i+1}$, $N_{G_i}(v) \cap W = N_{G_i}(v) \cap U^+_{i+1} \subseteq U^+_{i+1}$ by the second statement of Claim 1.
So for every $v \in S-(\{v_S\} \cup U^+_{i+1})$, $N_{G_i}(v) \cap W = N_{G_i}(\iota_{(S,j),\alpha}(v)) \cap W = N_{G_i}(\iota_{(S,j),\alpha}(v)) \cap U^+_{i+1} \subseteq U^+_{i+1}$, and hence $N_{G_i}(v) \cap W= N_{G_i}(v) \cap U^+_{i+1}$ (since $U^+_{i+1} \subseteq W$ by Claim 1).
Therefore, for every $v \in S-(\{v_S\} \cup U^+_{i+1})$, $N_{G_i}(v) \cap U^+_{i+1} = N_{G_i}(v) \cap W = N_{G_i}(\iota_{(S,j),\alpha}(v)) \cap W = N_{G_i}(\iota_{(S,j),\alpha}(v)) \cap U^+_{i+1}$.
$\Box$

\medskip

Let $U_{i+1} = \{v \in U_{i+1}^+: \deg_{G_{i+1}}(v) >d\}$.
Define the following:
	\begin{itemize}
		\item $\E_{i+1,0}=\{(S,j) \in \E_i: S \subseteq V(G_i) \cap V(G_{i+1})\}$.
		\item $\E_{i+1,1}=\{(\{v^*\} \cup (S \cap U^+_{i+1}) \cup T,j):$ there exists $(S,j) \in \E_i$ whose sink $v_S$ is in $O$, $T \subseteq U_{i+1}-S$, there exists a matching in $G_i$ between $T$ and $(S-\{v_S\}) \cap O$ with size $\lvert T \rvert\}$.
		\item $\E_{i+1,2}=\{(\{v^*\} \cup (S \cap U^+_{i+1}) \cup T,j+1):$ there exists $(S,j) \in \E_i$ whose sink $v_S$ is in $O$, $T \subseteq U_{i+1}-S$, $(S \cap U^+_{i+1}) \cup T \neq \emptyset$, there exist $u \in (S-\{v_S\}) \cap O$ and a matching in $G_i$ between $T$ and $(S-\{v_S,u\}) \cap O$ with size $\lvert T \rvert\}$.
		\item $\E_{i+1,3}=\{(\{v^*\} \cup U_{i+1}, 1)\}$.
		\item $\E_{i+1}=\E_{i+1,0} \cup \E_{i+1,1} \cup \E_{i+1,2} \cup \E_{i+1,3}$.
	\end{itemize}

\noindent{\bf Claim 3:} If there exist $\A_{i+1}$ and $\A'_{i+1}$ such that $(\E_{i+1},\A_{i+1},\A'_{i+1})$ satisfies (D9) and (D10), then $\E_{i+1}$ satisfies (D5).

\noindent{\bf Proof of Claim 3:}
Let $(S,j) \in \E_{i+1}$.
Since $\E_i$ satisfies (D5), every member of $\E_{i+1,0}$ satisfies (D5) by the definition of $D_{i+1}$.
So we may assume $(S,j) \in \E_{i+1}-\E_{i+1,0}$.
Then $S-\{v^*\} \subseteq U^+_{i+1} \subseteq W$.
So $\lvert S \rvert \leq 1+\lvert W \rvert \leq r$.
By Claim 1, $v^*$ is the unique vertex $v$ in $S$ such that $(u,v) \in D_{i+1}$ for every $u \in S-\{v\}$.
And if $(S,j) \not \in \E_{i+1,2}$, then $j \in [h-2]$ since $\E_i$ satisfies (D5).

So we may assume $(S,j) \in \E_{i+1,2}$ and it suffices to show $j \in [h-2]$.
Hence there exists $(S',j-1) \in \E_i$ such that $S=\{v^*\} \cup (S' \cap U^+_{i+1}) \cup T$ for some $T \subseteq U_{i+1}-S'$ with $(S' \cap U^+_{i+1}) \cup T \neq \emptyset$.
In particular, $S-\{v^*\} \neq \emptyset$ and $j \leq (h-2)+1=h-1$.
If there exist $\A_{i+1}$ and $\A'_{i+1}$ such that $(\A_{i+1},\A'_{i+1})$ satisfies (D9) and (D10), then since $S-\{v^*\} \neq \emptyset$ and $j \leq h-1$, contracting each member of $\A_{i+1,(S,j)} \cup \A'_{i+1,(S,j)}$ into a vertex creates a $(K_1 \vee (k+h-j)\CT_{j,k})$-minor in $G$, so there exists a $\CT_{j+1,k}$-minor in $G$.
Since $G$ is $\CT_{h,k}$-minor-free, $j+1 \leq h-1$.
So $j \in [h-2]$.
$\Box$

\medskip

\noindent{\bf Claim 4:} $(G_{i+1},\M_{i+1},D_{i+1},\E_{i+1})$ satisfies (D6) and (D7).

\noindent{\bf Proof of Claim 4:}
Let $v \in V(G_{i+1})$ such that either $\lvert V(M_{i+1,v}) \rvert \geq 2$ or $v$ is a head with respect to $D_{i+1}$ or the sink for some member of $\E_{i+1}$.
If $v \neq v^*$, then there exists at most $1 \leq N$ vertex $v'$ of $G_i$ such that $M_{i,v'} \subseteq M_{i+1,v}$.
If $v=v^*$, then $\{v' \in V(G_i): M_{i,v'} \subseteq M_{i+1,v^*}\} = O \subseteq N_{G_i[X]}[O_0] \subseteq N_{G_i[X]}^{\leq \ell_0-1}[z^*]$ which has size at most $\sum_{\alpha=0}^{\ell_0-1}d^\alpha \leq d^{\ell_0} \leq N$ by (i).
So (D6a) holds.
Suppose to the contrary that there exists $u \in V(G_{i+1})$ such that $uv \in E(G_{i+1})$, and either $\lvert V(M_{i+1,u}) \rvert \geq 2$ or $u$ is a head with respect to $D_{i+1}$ or the sink for some member of $\E_{i+1}$.
For $x \in \{u,v\}$, if $x \neq v^*$, then $x \in V(G_i) \cap V(G_{i+1})$, so either $\lvert V(M_{i,x}) \rvert = \lvert V(M_{i+1,x}) \rvert \geq 2$, or $x$ is a head with respect to $D_i$ (by the definition of $D_{i+1}$) or the sink for some member of $\E_{i+1,0} \subseteq \E_i$.
So if none of $u$ and $v$ equals $v^*$, then $\{u,v\} \subseteq V(G_i) \cap V(G_{i+1})$, and $uv \not \in E(G_i)$ (since $(G_i,\M_i,D_i,\E_i)$ satisfies (D6b)), so $uv \not \in E(G_{i+1})$, a contradiction. 
Hence one of $u$ and $v$ equals $v^*$.
Then the vertex $x \in \{u,v\}-\{v^*\}$ satisfies either $\lvert V(M_{i,x}) \rvert = \lvert V(M_{i+1,x}) \rvert \geq 2$, or $x$ is a head with respect to $D_i$ or the sink for some member of $\E_i$.
But $x \in N_{G_{i+1}}(v^*) \subseteq U^+_{i+1} \subseteq W$ by the first statement of Claim 1, contradicting the third statement of Claim 1.
This shows that $(G_{i+1},\M_{i+1},D_{i+1},\E_{i+1})$ satisfies (D6b).

Let $(S,j) \in \E_{i}$.
Let $v_S$ be the sink for $(S,j)$.
Let $v_{S'}$ be the vertex of $G_{i+1}$ such that $M_{i,v_S} \subseteq M_{i+1,v_{S'}}$.
Assume $S-\{v_S\} \subseteq V(G_{i+1})$.
Let $S'=\{v_{S'}\} \cup (N_{G_{i+1}}(v_{S'}) \cap (S-\{v_S\}))$.
To prove that $(G_{i+1},\M_{i+1},\E_{i+1})$ satisfies (D7), it suffices to show $(S',j) \in \E_{i+1}$.

We first assume $v_{S'} \neq v^*$.
Then $v_{S'}=v_S$.
Since $S-\{v_S\} \subseteq V(G_{i+1})$, $S \subseteq V(G_i) \cap V(G_{i+1})$ and $v^* \not \in S$.
So $S'=S$ and $(S',j)=(S,j) \in \E_{i+1,0} \subseteq \E_{i+1}$.

So we may assume $v_{S'}=v^*$.
Then $S'=\{v^*\} \cup (N_{G_{i+1}}(v^*) \cap (S-\{v_S\})) = \{v^*\} \cup (U^+_{i+1} \cap S)$ and $v_S \in O$.
So $(S',j) \in \E_{i+1,1} \subseteq \E_{i+1}$ (by taking $T=\emptyset$).
Hence $(G_{i+1},\M_{i+1},\E_{i+1})$ satisfies (D7).
$\Box$

\medskip

Recall that $q_{i+1}=v^*$.

\medskip

\noindent{\bf Claim 5:} $(G_{i+1},\M_{i+1},\E_{i+1},D_{i+1})$ satisfies (D8a)-(D8h).

\noindent{\bf Proof of Claim 5:}
By definition and Claim 1, $U_{i+1} \subseteq U_{i+1}^+ \subseteq W \cap V(G_i) \cap V(G_{i+1}) \subseteq V(G) \cap V(G_i) \cap V(G_{i+1})$.
Note that $V(G_i)-V(G_{i+1}) \subseteq O \subseteq X$.
So (D8a) follows from (i).
And (D8b) and (D8c) follow from the definition of $U_{i+1}$.
If $(S,j) \in \E_{i}$ with $M_{i,v_S} \subseteq M_{i+1,q_{i+1}}$, where $v_S$ is the sink for $(S,j)$, then (i) implies that $\{x \in S-\{v_S\}: \deg_{G_i}(x)>d\} \subseteq N_{G_i}(v_S)-X \subseteq N_{G_i}[O]-X = N_{G_i}(O)-X=U_{i+1}^+$, so (D8d) holds.

Now we prove (D8e).
Let $v \in V(G) \cap V(G_{i+1}) \cap N_{G_i}(v')-U^+_{i+1}$ with $\deg_{G_i}(v) \leq d$ for some $v' \in V(G_i)$ with $M_{i,v'} \subseteq M_{i+1,q_{i+1}}$.
Then $v \in V(G_{i+1}) \cap N_{G_i}[O]-U^+_{i+1} = (V(G_{i+1}) \cap N_{G_i}[O])-(N_{G_i}(O)-X) = V(G_{i+1}) \cap N_{G_i}(O) - (N_{G_i}(O)-X) = V(G_{i+1}) \cap N_{G_i}(O) \cap X \subseteq V(G_{i+1}) \cap N_{G_i[X]}(O)$.
Recall that $N_{G_i[X]}[O] \subseteq N_{G_i[X]}^{\leq 3\lvert Y \rvert}[V(P'_{z^*})] \subseteq N_{G_i[X]}^{\leq \ell_0-1}[z^*]$.
So $N_{G_i}(v)-X \subseteq N_{G_i}[N_{G_i[X]}[O]]-X \subseteq N_{G_i}[N_{G_i[X]}^{\leq \ell_0-1}[z^*]]-X =W$ by (iii).
In addition, by the definition of $G_{i+1}$, $N_{G_{i+1}}(v) \cap V(G) \subseteq N_{G_i}(v) \cap V(G_{i+1}) \cap V(G_i)$.
So $\{x \in N_{G_{i+1}}(v) \cap V(G): \deg_{G_{i+1}}(x)>d\} \subseteq \{x \in N_{G_i}(v): \deg_{G_i}(x)>d\} \subseteq N_{G_i}(v)-X$ by (i).
Therefore, $\{x \in N_{G_{i+1}}(v) \cap V(G): \deg_{G_{i+1}}(x)>d\} \subseteq N_{G_i}(v)-X \subseteq N_{G_i}(v) \cap W$.
Since $v \in N_{G_i[X]}(O)$, $\phi(v) \in Y$ by (vii).
So by (viii), there exists $v'' \in O$ with $\phi(v'')=\phi(v)$.
By the definition of $\phi$, we have $N_{G_i}(v) \cap W = N_{G_i}(v'') \cap W \subseteq N_{G_i}[O] \cap W = N_{G_i}(O) \cap W = U_{i+1}^+$ by Claim 1.
So $\{x \in N_{G_{i+1}}(v) \cap V(G): \deg_{G_{i+1}}(x)>d\} \subseteq \{x \in N_{G_i}(v) \cap W: \deg_{G_{i+1}}(x)>d\} \subseteq \{x \in U_{i+1}^+: \deg_{G_{i+1}}(x)>d\} = U_{i+1}$.
Hence (D8e) holds.

Since $N_{G_{i+1}}(q_{i+1})=U^+_{i+1}$, $U_{i+1} \cap N_{G_{i+1}}(q_{i+1}) = U_{i+1}$.
So $((U_{i+1} \cap N_{G_{i+1}}(q_{i+1})) \cup \{q_{i+1}\}, 1) = (U_{i+1} \cup \{q_{i+1}\},1) \in \E_{i+1,3} \subseteq \E_{i+1}$, and $q_{i+1}$ is its sink.
Hence (D8f) holds.

(D8ga) and (D8gb) clearly hold.
Now we prove (D8gc).
Let $u_{i}v_{i}$ be an edge of $G_i$ such that $u_{i+1} \neq v_{i+1}$ and $u_{i+1}v_{i+1} \not \in E(G_{i+1})$, where $u_{i+1},v_{i+1}$ are the vertices of $G_{i+1}$ with $M_{i,u_i} \subseteq M_{i+1,u_{i+1}}$ and $M_{i,v_i} \subseteq M_{i+1,v_{i+1}}$.
Then the definition of $G_{i+1}$ implies that one of $u_i,v_i$ is in $O$ and the other is in $N_{G_i}(O) \cap X$.
So $q_{i+1}=v^* \in \{u_{i+1},v_{i+1}\}$ and the vertex $x \in \{u_{i+1},v_{i+1}\}-\{q_{i+1}\}$ is in $N_{G_i}(O) \cap X$.
Hence $x \in V(G_i) \cap V(G_{i+1})$ and there exists $y \in O$ such that $x \in N_{G_i}(y)$.
Since $x \in X$, (i) implies that $\deg_{G_i}(x) \leq d$.
And $x \in X$ implies $x \not \in U^+_{i+1}$.
So it suffices to show $x \in V(G)$ and $x$ is not the sink for a member of $\E_i$.
Suppose to the contrary that $\lvert V(M_{i,x}) \rvert  = \lvert V(M_{i+1,x}) \rvert \geq 2$ or $x$ is the sink for some member of $\E_i$.
If $y \in O_0$, then $x \in O$, a contradiction.
So $y \in O-O_0$.
That is, $y \in N_{G_i[X]}(O_0)$, and either $\lvert V(M_{i,y}) \rvert \geq 2$, or $y$ is the sink for some member of $\E_i$.
Since $(G_i,\M_i,D_i,\E_i)$ satisfies (D6b), $xy \not \in E(G_i)$, a contradiction.
Hence (D8gc) holds.

Now we prove (D8gd).
Let $(S,j) \in \E_i$ with $M_{i,v_S} \subseteq M_{i+1,q_{i+1}}$, where $v_S$ is the sink for $(S,j)$.
So $v_S \in O$.
By Claim 2, there exists $(S',j) \in \E_i$ with $S' \subseteq O_1 \cup (N_{G_i}(O_1) \cap U^+_{i+1})$, $S' \cap U^+_{i+1} = S \cap U^+_{i+1}$, and $v_{S'} \in O_1 \subseteq O$, where $v_{S'}$ is the sink for $(S',j)$, and there exists a bijection $\iota: S-(\{v_S\} \cup U^+_{i+1}) \rightarrow S'-(\{v_{S'}\} \cup U^+_{i+1})$ such that $N_{G_i}(\iota(v)) \cap U^+_{i+1}=N_{G_i}(v) \cap U^+_{i+1}$ for every $v \in S-(\{v_S\} \cup U^+_{i+1})$.
So $S'-U^+_{i+1} \subseteq O_1 \subseteq O$.
Hence $\bigcup_{s \in S'-U^+_{i+1}}M_{i,s} \subseteq M_{i+1,q_{i+1}}$.
Therefore (D8gd) holds.

Since $\bigcup_{M \in \M_i}V(M) = \bigcup_{M \in \M_{i+1}}V(M)$, (D8h) holds.
$\Box$

\medskip

\noindent{\bf Claim 6:} $(G_{i+1},\M_{i+1},\E_{i+1},D_{i+1})$ satisfies (D8).

\noindent{\bf Proof of Claim 6:}
By Claim 5, it suffices to show (D8i).
Let $(S,j) \in \E_i$ with $M_{i,v} \cup M_{i,v_S} \subseteq M_{i+1,q_{i+1}}$ for some $v \in S-\{v_S\}$, where $v_S$ is the sink for $(S,j)$.
So $v_S \in O$.
By Claim 2, there exists $(S',j) \in \E_i$ with $S' \subseteq O \cup (N_{G_i}(O) \cap U^+_{i+1})$, $S' \cap U^+_{i+1}=S \cap U^+_{i+1}$ and $v_{S'} \in O$, where $v_{S'}$ is the sink for $(S',j)$, and there exists a bijection $\iota: S-(\{v_S\} \cup U^+_{i+1}) \rightarrow S'-(\{v_{S'}\} \cup U^+_{i+1})$ such that $N_{G_i}(\iota(v)) \cap U^+_{i+1} = N_{G_i}(v) \cap U^+_{i+1}$ for every $v \in S-(\{v_S\} \cup U^+_{i+1})$.
Since $U^+_{i+1} \cup (S-\{v_S\}) \cup (S'-\{v_{S'}\}) \subseteq V(G) \cap V(G_i)$ (by (D4) and (D5)), (D2) implies that $N_{G}(\iota(v)) \cap U^+_{i+1} = N_{G_i}(\iota(v)) \cap U^+_{i+1} = N_{G_i}(v) \cap U^+_{i+1} = N_{G}(v) \cap U^+_{i+1}$ for every $v \in S-(\{v_S\} \cup U^+_{i+1})$.

Let $\C$ be the multiset $\{N_G(x) \cap U_{i+1} \neq \emptyset: x \in S-(\{v_S\} \cup U_{i+1}^+)\}$.
So for every $T \in \C$, there exists $x_T \in S-(\{v_S\} \cup U_{i+1}^+)$ such that $T=N_G(x_T) \cap U_{i+1} = N_{G_i}(\iota(x_T)) \cap U_{i+1}$.
Note that we can choose those $x_T$ such that $x_{T_1} \neq x_{T_2}$ for any distinct members $T_1,T_2$ of the multiset $\C$.
Let $f$ be a function that maps each member $T$ of $\C$ to a vertex in $T$.
Then there exists a matching in $G_i$ between $\{f(T): T \in \C\}$ (as a set) and $\{\iota(x_T): T \in \C\}$ with size $\lvert \{f(T): T \in \C\} \rvert$ (as a set).
So there exists a matching in $G_i$ between $\{f(T): T \in \C\}-S' \subseteq U_{i+1}-S'$ and $\{\iota(x_T): T \in \C\}$ with size $\lvert \{f(T): T \in \C\}-S' \rvert$.
Note that $\{\iota(x_T): T \in \C\} \subseteq S'-(\{v_{S'}\} \cup U^+_{i+1}) \subseteq (S'-\{v_{S'}\}) \cap O$. 
Hence $(\{v^*\} \cup (S \cap U^+_{i+1}) \cup \{f(T): T \in \C\}, j) = (\{v^*\} \cup (S' \cap U^+_{i+1}) \cup (\{f(T): T \in \C\}-S'), j) \in \E_{i+1,1} \subseteq \E_{i+1}$.
So (D8ia) holds.

Let $u \in S-(\{v_S\} \cup U_{i+1}^+)$.
So $\iota(u) \in S'-(\{v_{S'}\} \cup U^+_{i+1}) \subseteq (S'-\{v_{S'}\}) \cap O$.
Let $\C_u$ be the multiset $\{N_G(x) \cap U_{i+1} \neq \emptyset: x \in S-(\{v_S,u\} \cup U_{i+1}^+)\}$.
So for every $T \in \C_u$, there exists $y_T \in S-(\{v_S,u\} \cup U_{i+1}^+)$ such that $T=N_G(y_T) \cap U_{i+1} = N_{G_i}(\iota(y_T)) \cap U_{i+1}$.
Note that we can choose those $y_T$ such that $y_{T_1} \neq y_{T_2}$ for any distinct members $T_1,T_2$ of the multiset $\C_u$.
Assume either $S \cap U_{i+1}^+ \neq \emptyset$ or $\C_u \neq \emptyset$.
Let $f_u$ be a function that maps each member $T$ of $\C_u$ to a vertex in $T$.
Then there exists a matching in $G_i$ between $\{f_u(T): T \in \C_u\}$ (as a set) and $\{\iota(y_T): T \in \C_u\}$ with size $\lvert \{f_u(T): T \in \C_u\} \rvert$ (as a set).
Let $T^*$ be the set $\{f_u(T): T \in \C_u\}-(S' \cap U^+_{i+1})$.
So there exists a matching in $G_i$ between $T^* \subseteq U_{i+1}-S'$ and $\{\iota(y_T): T \in \C_u\} \subseteq S'-(\{v_{S'},\iota(u)\} \cup U^+_{i+1}) \subseteq (S'-\{v_{S'},\iota(u)\}) \cap O$ with size $\lvert T^* \rvert$.
Since either $S \cap U_{i+1}^+ \neq \emptyset$ or $\C_u \neq \emptyset$, we know $(S' \cap U^+_{i+1}) \cup T^* = (S \cap U^+_{i+1}) \cup \{f_u(T): T \in \C_u\} \neq \emptyset$.
Hence $(\{v^*\} \cup (S \cap U^+_{i+1}) \cup \{f_u(T): T \in \C_u\}, j+1) = (\{v^*\} \cup (S' \cap U^+_{i+1}) \cup T^*, j+1) \in \E_{i+1,2} \subseteq \E_{i+1}$.
Therefore, (D8i) holds and hence (D8) holds.
$\Box$

\medskip

For every $Q \in \E_{i+1,0}$, we know $Q \in \E_i$, and we define $\A_{i+1,Q}$ and $\A'_{i+1,Q}$ to be $\A_{i,Q}$ and $\A'_{i,Q}$, respectively.
Since $(\A_i,\A'_i)$ satisfies (D10), (D10) holds for every member of $\E_{i+1,0}$.

For every $Q \in \E_{i+1,1}$, we define the following:
	\begin{itemize}
		\item Let $(S_Q,j_Q)$ be a member of $\E_i$ with $v_{S_Q} \in O$ such that $Q=(\{v^*\} \cup (S_Q \cap U^+_{i+1}) \cup T_Q, j_Q)$, where $T_Q \subseteq U_{i+1}-S_Q$, $v_{S_Q}$ is the sink for $(S_Q,j_Q)$, and there exists a matching in $G_i$ between $T_Q$ and $(S_Q-\{v_{S_Q}\}) \cap O$ with size $\lvert T_Q \rvert$.
		\item Let $(S'_Q,j_Q)$ be the member $(S_1,j)$ of $\E_i$ and $\iota_{(S_Q,j_Q),1}$ be the bijection mentioned in Claim 2 when taking $(S,j)=(S_Q,j_Q)$ and $\alpha=1$.
			Let $v_{S'_Q}$ be the sink for $(S'_Q,j_Q)$.

			(Note that Claim 2 implies that $S_Q \cap U^+_{i+1}=S'_Q \cap U^+_{i+1}$ and there exists a matching $M_{T_Q}$ in $G_i$ between $T_Q$ and $\iota_{(S_Q,j_Q),1}((S_Q-\{v_{S_Q}\}) \cap O)$ with size $\lvert T_Q \rvert$.)
		\item Let $\A_{i+1,Q}=\A_{i,(S'_Q,j_Q)}$.
		\item For every vertex $v \in S'_Q-\{v_{S'_Q}\}$, let $A_v$ be the member of $\A'_{i,(S'_Q,j_Q)}$ such that $v \in V(A_v)$.
		\item For every vertex $u \in T_Q$, let $u'$ be the vertex in $\iota_{(S_Q',j_Q),1}((S_Q-\{v_{S_Q}\}) \cap O) \subseteq O$ matched with $u$ in $M_{T_Q}$, and let $A_{i+1,u}=G[V(A_{u'}) \cup \{u\}]$.
		\item Let $\A'_{i+1,Q}=\{A_v: v \in S'_Q \cap U^+_{i+1}\} \cup \{A_{i+1,u}: u \in T_Q\}$.
	\end{itemize}
Clearly, (D10) holds for every member of $\E_{i+1,1}$.

For every $Q \in \E_{i+1,2}$, we define the following:
	\begin{itemize}
		\item Let $(S_Q,j_Q)$ be a member of $\E_i$ with $v_{S_Q} \in O$ such that $Q=(\{v^*\} \cup (S_Q \cap U^+_{i+1}) \cup T_Q, j_Q+1)$, where $T_Q \subseteq U_{i+1}-S_Q$, $v_{S_Q}$ is the sink for $(S_Q,j_Q)$, and there exist $u_Q \in (S_Q-\{v_{S_Q}\}) \cap O$ and a matching in $G_i$ between $T_Q$ and $(S_Q-\{v_{S_Q},u_Q\}) \cap O$ with size $\lvert T_Q \rvert$.
		\item For every $\alpha \in [h+k]$, let $(S_\alpha,j_Q)$ be the member $(S_\alpha,j)$ of $\E_i$ and $\iota_\alpha$ be the bijection $\iota_{(S,j),\alpha}$ mentioned in Claim 2 when taking $(S,j)=(S_Q,j_Q)$ and $\alpha=\alpha$. 
			Let $v_{S_\alpha}$ be the sink for $(S_\alpha,j_Q)$.

			(Note that Claim 2 implies that for every $\alpha \in [h+k]$, $S_Q \cap U^+_{i+1}=S_\alpha \cap U^+_{i+1}$, $v_{S_\alpha} \in O_\alpha$, and there exists a matching $M_{T_Q,\alpha}$ in $G_i$ between $T_Q$ and $\iota_\alpha((S_Q-\{v_{S_Q},u_Q\}) \cap O) = S_\alpha \cap O_\alpha - \{v_{S_\alpha}, \iota_\alpha(u_Q)\}$ with size $\lvert T_Q \rvert$.)
		\item For every $\alpha \in [h+k-j_Q-1]$, define the following:
			\begin{itemize}
				\item For every $v \in S_\alpha-\{v_{S_\alpha}\}$, let $A_{\alpha,v}$ be the member of $\A'_{i,(S_\alpha,j_Q)}$ with $v \in V(A_{\alpha,v})$.
				\item Let $\A^{(1)}_{i,(S_\alpha,j_Q)}$ be a subset of $\A_{i,(S_\alpha,j_Q)}$ such that contracting each member of $\A^{(1)}_{i,(S_\alpha,j_Q)}$ into a vertex creates a $\CT_{j_Q,k}$-minor.

					(Note that $\A^{(1)}_{i,(S_\alpha,j_Q)}$ exists since $\A_{i,(S_\alpha,j_Q)}$ satisfies (D10).)
				\item Let $A_{i+1,\alpha}=G[V(A_{\alpha,\iota_\alpha(u_Q)}) \cup \bigcup_{A \in \A^{(1)}_{i,(S_\alpha,j_Q)}}V(A)]$.
				\item Let $\A^{(2)}_{i,(S_\alpha,j_Q)}$ be a subset of $\A_{i,(S_\alpha,j_Q)}-\A^{(1)}_{i,(S_\alpha,j_Q)}$ such that contracting each member of $\A^{(2)}_{i,(S_\alpha,j_Q)}$ into a vertex creates a $k\CT_{j_Q,k}$-minor.

					(Note that $\A^{(2)}_{i,(S_\alpha,j_Q)}$ exists since $\A_{i,(S_\alpha,j_Q)}$ satisfies (D10) and $j_Q \in [h-2]$ by (D5).)
				\item Let $\A_{i+1,\alpha} = \{A_{i+1,\alpha}\} \cup \A^{(2)}_{i,(S_\alpha,j_Q)}$.

			(Note that since $\A_{i,(S_\alpha,j_Q)}$ and $\A'_{i,(S_\alpha,j_Q)}$ satisfy (D10), contracting each member of $\A_{i+1,\alpha}$ into a vertex creates a $(K_1 \vee k\CT_{j_Q,k})$-minor and hence a $\CT_{j_Q+1,k}$-minor.)
			\end{itemize}
		\item Let $\A_{i+1,Q}=\bigcup_{\alpha \in [h+k-j_Q-1]}\A_{i+1,\alpha}$.

			(Note that members of $\A_{i+1,Q}$ are pairwise disjoint since $(\A_i,\A'_i)$ satisfies (D11) and $O_1,O_2,..., \allowbreak O_{h+k}$ are pairwise disjoint.
			And contracting each member of $\A_{i+1,Q}$ into a vertex creates a $(k+h-j_Q-1)\CT_{j_Q+1,k}$-minor.)
		\item For every vertex $w \in T_Q$, let $A_{i+1,w}=G[\{w\} \cup \bigcup_{\alpha \in [h+k-j_Q-1]}V(A_{\alpha,w_\alpha})]$, where for every $\alpha \in [h+k-j_Q-1]$, $w_\alpha$ is the vertex in $S_\alpha \cap O_\alpha-\{v_{S_\alpha},\iota_\alpha(u_Q)\}$ matched with $w$ in $M_{T_Q,\alpha}$. 

			(Note that $A_{i+1,w}$ is connected.) 
		\item Let $\A'_{i+1,Q}=\{G[V(\bigcup_{\alpha \in [h+k-j_Q-1]}A_{\alpha,v})]: v \in S_Q \cap U^+_{i+1}\} \cup \{A_{i+1,w}: w \in T_Q\}$.
	\end{itemize}
Since $(\A_i,\A'_i)$ satisfies (D10) and (D11), we know that (D10) holds for every member of $\E_{i+1,2}$.

For the unique member $(\{v^*\} \cup U_{i+1}, 1)$ of $\E_{i+1,3}$, define $\A_{i+1,(\{v^*\} \cup U_{i+1},1)} = \{G[\bigcup_{v \in O_\alpha}V(M_{i,v})]: \alpha \in [h+k-1]\}$, and $\A'_{i+1,(\{v^*\} \cup U_{i+1},1)} = \{G[\{v\}]: v \in U_{i+1}\}$.
By (iii), (viii) and the definition of $\phi$, (D10) holds for $(\{v^*\} \cup U_{i+1},1)$.

Therefore, $(\E_{i+1},\A_{i+1},\A'_{i+1})$ satisfies (D9) and (D10).
Hence $\E_{i+1}$ satisfies (D5) by Claim 3.

\medskip

\noindent{\bf Claim 7:} $(\A_{i+1},\A'_{i+1})$ satisfies (D11).

\noindent{\bf Proof of Claim 7:}
Suppose to the contrary that there exist members $Q_1=(S_1,j_1)$ and $Q_2=(S_2,j_2)$ of $\E_{i+1}$ with distinct sinks violating (D11).
Since $(\A_i,\A'_i)$ satisfies (D11), we may assume $Q_1 \not \in \E_{i+1,0}$ by symmetry.
Hence $Q_1 \in \bigcup_{\alpha=1}^3\E_{i+1,\alpha}$, so the sink of $Q_1$ is $v^*$.
Since $Q_1$ and $Q_2$ have distinct sinks, $Q_2 \in \E_{i+1,0}$.
In particular, by the third statement of Claim 1, the sink for $Q_2$ is in $V(G_i)-(W \cup O)$.
This implies $S_2 \cap O=\emptyset$ by the definition of $O$ and (D6) (for $(G_i,\M_i,D_i,\E_i)$).
Since $(\A_i,\A'_i)$ satisfies (D10), it is straightforward to verify that $Q_1 \not \in \E_{i+1,3}$.
So $Q_1 \in \E_{i+1,1} \cup \E_{i+1,2}$.

Hence $(S_{Q_1},j_{Q_1}) \in \E_i$ is defined; when $Q_1 \in \E_{i+1,1}$, $(S'_{Q_1},j_{Q_1})$ and $\iota_{(S_{Q_1},j_{Q_1}),1}$ are defined; when $Q_1 \in \E_{i+1,2}$, $(S_\alpha, j_{Q_1})$ and $\iota_\alpha$ (for $\alpha \in [h+k]$) are defined. 
Note that the sinks for $(S'_{Q_1},j_{Q_1})$ and $(S_\alpha,j_{Q_1})$ are in $O$.
Since $(\A_i,\A_i')$ satisfies (D10) and (D11), every member of $\A_{i+1,Q_1}$ is disjoint from every member of $\A_{i+1,Q_2} \cup \A'_{i+1,Q_2}$.
Moreover, for every $A \in \A'_{i+1,Q_1}$, $V(A) \subseteq U^+_{i+1} \cup O \cup (V(G)-\bigcup_{v \in V(G_i)}V(M_{i,v}))$ and $V(A) \cap U^+_{i+1}=V(A) \cap U^+_{i+1} \cap S_1$ is a set with size at most 1.
Hence every member of $\A'_{i+1,Q_1}$ is disjoint from every member of $\A_{i+1,Q_2}$, and if some member $A_1 \in \A'_{i+1,Q_1}$ intersects some member $A_2 \in \A'_{i+1,Q_2}$, then $V(A_1) \cap V(A_2) \subseteq V(A_1) \cap V(A_2) \cap U^+_{i+1} \subseteq (V(A_1) \cap S_{1}) \cap (V(A_2) \cap U^+_{i+1} \cap V(G_{i+1})) \subseteq V(A_1) \cap S_{1} \cap V(A_2) \cap S_{2}$. 
So $Q_1$ and $Q_2$ do not form a counterexample, a contradiction.
$\Box$

\medskip

To prove this lemma, it suffices to show that $(G_{i+1},\M_{i+1},D_{i+1},\E_{i+1})$ satisfies (D12).
Let $v \in V(G_{i+1})$ such that either $\lvert V(M_{i+1,v}) \rvert \geq 2$, or $v$ is a head with respect to $D_{i+1}$ or the sink for some member of $\E_{i+1}$.
It suffices to show $\deg_{G_{i+1}}(v) \leq r$.

Suppose to the contrary that $\deg_{G_{i+1}}(v) > r$.
By the definition of $G_{i+1}$, $\deg_{G_{i+1}}(v^*) =\lvert U^+_{i+1} \rvert \leq \lvert W \rvert \leq r-1$.
So $v \neq v^*$.
Hence $v \in V(G_i) \cap V(G_{i+1})$.
In particular, $\deg_{G_{i+1}}(v) \leq \deg_{G_i}(v)$ and $M_{i,v}=M_{i+1,v}$.
So $\deg_{G_i}(v)>r$.
Since $(G_{i},\M_i,D_{i},\E_{i})$ satisfies (D12), $\lvert V(M_{i+1,v}) \rvert = \lvert V(M_{i,v}) \rvert=1$, $v$ is not a head with respect to $D_{i}$, and $v$ is not the sink for some member of $\E_{i}$.
Hence by the definition of $D_{i+1}$, since $v \neq v^*$, $v$ is not a head with respect to $D_{i+1}$.
So $v$ is a sink for some member of $\E_{i+1}$.
By the definition of $\E_{i+1}$, $v$ is the sink for some member of $\E_{i+1,0}$ since $v \neq v^*$.
However, $\E_{i+1,0} \subseteq \E_i$.
So $v$ is the sink for some member of $\E_i$, a contradiction.
This proves the lemma.
\end{pf}

\begin{lemma} \label{strong_scheme}
For any positive integers $h \geq 3$ and $k$, there exist positive integers $r=r(h,k),d=d(h,k),N=N(h,k)$ such that for every graph $G$ with no $\CT_{h,k}$-minor, there exists a strong $(G,h,k,r,d,N)$-defective elimination scheme.
\end{lemma}

\begin{pf}
Let $h \geq 3$ and $k$ be positive integers.
Define $r=r_{\ref{contract}}(h,k)$, $d=d_{\ref{contract}}(h,k)$ and $N_{\ref{contract}}(h,k)$, where $r_{\ref{contract}},d_{\ref{contract}},N_{\ref{contract}}$ are the integers $r,d,N$ mentioned in Lemma \ref{contract}, respectively.

Let $G$ be a graph with no $\CT_{h,k}$-minor.
Let $(G_1,\M_1,\E_1,D_1,\A_1,\A_1')=(G,\{G[\{v\}]: v \in V(G)\},\emptyset, \emptyset, \emptyset, \emptyset)$.
So the sequence $((G_\alpha,\M_\alpha,\E_\alpha,D_\alpha,\A_\alpha,\A'_\alpha): \alpha \in [1])$ is a strong 1-$(G,h,k,r,d,N)$-defective elimination scheme.

Suppose to the contrary that there exists no strong $(G,h,k,r,d,N)$-defective elimination scheme.
Then there exists the maximum positive integer $i$ such that the sequence $((G_\alpha,\M_\alpha,\E_\alpha,D_\alpha, \allowbreak \A_\alpha,\A'_\alpha): \alpha \in [i])$ is a strong $i$-$(G,h,k,r,d,N)$-defective elimination scheme.
By Lemma \ref{contract} and the maximality of $i$, $\lvert V(G_i) \rvert \leq N$.
For every $\alpha \in {\mathbb N}-[i]$, define $(G_\alpha,\M_\alpha,\E_\alpha,D_\alpha,\A_\alpha,\A'_\alpha)=(G_i,\M_i,\E_i,D_i,\A_i,\A'_i)$.
Then $((G_\alpha,\M_\alpha,\E_\alpha,D_\alpha,\A_\alpha,\A'_\alpha): \alpha \in {\mathbb N})$ is a strong $(G,h,k,r,d,N)$-elimination scheme, a contradiction.
\end{pf}

\begin{theorem} \label{main}
For any positive integers $h \geq 3$ and $k$, there exists a positive integer $d^*$ such that every graph with no $\CT_{h,k}$-minor has an $(h-1)$-coloring with defect $d^*$.
\end{theorem}

\begin{pf}
Let $h \geq 3$ and $k$ be positive integers.
Let $r=r_{\ref{strong_scheme}}(h,k)$, $d=d_{\ref{strong_scheme}}(h,k)$ and $N=N_{\ref{strong_scheme}}(h,k)$, where $r_{\ref{strong_scheme}},d_{\ref{strong_scheme}},N_{\ref{strong_scheme}}$ are the integers $r,d,N$ mentioned in Lemma \ref{strong_scheme}.
Define $d^*=d_{\ref{scheme_color}}(d,N)$, where $d_{\ref{scheme_color}}$ is the integer $d^*$ mentioned in Lemma \ref{scheme_color}.

Let $G$ be a graph with no $\CT_{h,k}$-minor.
By Lemma \ref{strong_scheme}, there exists a strong $(G,h,k,r,d,N)$-defective elimination scheme.
So there exists a $(G,h,k,r,d,N)$-defective elimination scheme. 
By Lemma \ref{scheme_color}, there exists an $(h-1)$-coloring of $G[V(G_1) \cap V(G)]=G$ with defect $d^*$.
\end{pf}

\bigskip

Now we are ready to prove Theorem \ref{main_intro}.

\bigskip

\noindent{\bf Proof of Theorem \ref{main_intro}:}
Let $\F$ be a minor-closed family.
If $\omega_\Delta(\F)=0$, then $K_1 \not \in \F$, so $\F=\emptyset$.
If $\omega_\Delta(\F)=1$, then $K_{1,k} \not \in \F$ for some integer $k$, so $K_{1,k}$ is not a minor (and hence a subgraph) of any graph in $\F$, so $\F$ has bounded maximum degree, and hence $\chi_\Delta(\F)=1$.
If $\omega_\Delta(\F)=\infty$, then $\chi_\Delta(\F)=\infty$.
So we may assume that $2 \leq \omega_\Delta(\F) <\infty$.
Hence there exist a positive integer $k$ such that $\CT_{\omega_\Delta(\F)+1,k} \not \in \F$.
So every graph in $\F$ is $\CT_{\omega_\Delta(\F)+1,k}$-minor free.
Hence $\chi_\Delta(\F) \leq \omega_\Delta(\F)$ by Theorem \ref{main}.
$\Box$

\bigskip

\bigskip

\noindent{\bf Acknowledgment:} This paper was partially completed when the author was visiting the Institute of Mathematics at Academia Sinica at Taiwan. The author thanks the Institute of Mathematics at Academia Sinica for its hospitality. The author thanks anonymous reviewers for their careful reading and suggestions.

\end{document}